\documentclass[11pt, reqno]{amsart}
 
\usepackage[cp1251]{inputenc}
\usepackage{amsmath,amsthm,mathrsfs}
\usepackage{amssymb}
\usepackage[english]{babel}
\usepackage{esint}
\usepackage[titletoc]{appendix}

\usepackage{graphicx}
\usepackage{float}
\usepackage{subfigure}

\usepackage{bbm}
\usepackage{amsfonts,amstext,amssymb,verbatim,epsfig}
\usepackage{pgf,tikz}
\usepackage{float}
\usetikzlibrary{arrows}
\usepackage{hyperref}
\usepackage{cite}
\usepackage{epstopdf}
\usepackage{color}
\usepackage{transparent}
\usepackage{subfigure}
\usepackage{pgfplots}
\usepgfplotslibrary{polar}
\usepgflibrary{shapes.geometric}
\usetikzlibrary{calc}

\usepackage{tikz}

\usepackage{xcolor}
\usepackage[normalem]{ulem}

\hypersetup{colorlinks=true,pdfborder=001,  citecolor  = blue, citebordercolor= {magenta}}
\hypersetup{linkcolor=magenta, linkbordercolor = blue}
\hypersetup{colorlinks,urlcolor=blue}
\makeatletter
\let\reftagform@=\tagform@
\def\tagform@#1{\maketag@@@{(\ignorespaces\textcolor{magenta}{#1}\unskip\@@italiccorr)}}
\renewcommand{\eqref}[1]{\textup{\reftagform@{\ref{#1}}}}
\makeatother

\makeatletter

\DeclareUrlCommand\ULurl@@{%
  \def\UrlLeft{\uline\bgroup}%
  \def\UrlRight{\egroup}}
\def\ULurl@#1{\hyper@linkurl{\ULurl@@{#1}}{#1}}
\DeclareRobustCommand*\ULurl{\hyper@normalise\ULurl@}
\makeatother

\def\lessim{\ \lower4pt\hbox{$
		\buildrel{\displaystyle <}\over\sim$}\ }
\def\gessim{\ \lower4pt\hbox{$\buildrel{\displaystyle >}
		\over\sim$}\ }

\newtheorem{definition}{\bf Definition}
\newtheorem{theorem}{\bf Theorem}
\newtheorem{lemma}[theorem]{\bf Lemma}
\newtheorem{corollary}[theorem]{\bf Corollary}
\newtheorem{example}{\bf Example}

\newtheorem{proposition}[theorem]{\bf Proposition}

\theoremstyle{remark}
\newtheorem{remark}{Remark}

\newenvironment{Proof of lemma}{\noindent{\bf Proof of Lemma}}{\hfill$\Box$\newline}
\newenvironment{Proof of theorem}{\noindent{\bf Proof of Theorem}}{\hfill{\footnotesize${\square}$}\newline}
\newenvironment{Proof of theorems}{\noindent{\bf Proof of Theorems}}{\hfill$\Box$\newline}
\newenvironment{Proof of proposition}{\noindent{\bf Proof of Proposition}}{\hfill$\Box$\newline}
\newenvironment{Proof of propositions}{\noindent{\bf Proof of Propositions}}{\hfill$\Box$\newline}
\newenvironment{Proof of exercise}{\noindent{\it Proof of Exercise:}}{\hfill$\Box$}
\begin{document}
\title{The spherical $p+s$ spin glass at zero temperature}

\author{Antonio Auffinger}
\address{Department of Mathematics, Northwestern Universty}
\email{tuca@northwestern.edu}


\author{Yuxin Zhou}
\address{Department of Mathematics, Northwestern Universty}
\email{yuxinzhou2023@u.northwestern.edu}

\begin{abstract}
We determine the structure of the Parisi measure at zero temperature for the spherical $p+s$ spin glass model. We show that depending on the values of $p$ and $s$,  four scenarios may emerge, including the existence of 1-FRSB and 2-FRSB phases as predicted by Crisanti and Leuzzi \cite{CL1,CL2}. We also provide consequences for the model at low temperature.
 \end{abstract}

\maketitle
\section{Introduction and main results}

Let $N, p,s$ be integers with $p,s \geq 2$, $N\geq1$ and let $\lambda \in [0,1]$.  The Hamiltonian of the spherical $p+s$ spin model is the Gaussian function defined on the $N$-dimensional sphere  $S_{N}:= \big \{ \sigma \in \mathbb{R}^N: \sum^N_{i=1} \sigma^2_i=N \big\}$ by
\begin{eqnarray*}
H_N(\sigma):=\frac{\sqrt{\lambda}}{N^{\frac{p-1}{2}}} \sum_{1 \leq i_1,\cdots,i_p \leq N} g_{i_1,\dots,i_p} \sigma_{i_1} \cdots \sigma_{i_p}+\frac{\sqrt{1-\lambda}}{N^{\frac{s-1}{2}}} \sum_{1 \leq i_1,\cdots,i_s \leq N} g_{i_1,\cdots,i_s} \sigma_{i_1} \cdots \sigma_{i_s},
\end{eqnarray*}
where all $(g_{i_1,\cdots, i_p})$ and $(g_{i_1,\cdots, i_s})$, $1 \leq i_1,\cdots,i_s, i_{p} \leq N$,  are independent, identically distributed standard Gaussian random variables.

The Gaussian field $H_{N}$ is centered with covariance given by 
\begin{eqnarray*}
\mathbb{E} H_N(\sigma^1) H_N(\sigma^2)=N \xi(R_{1,2})
\end{eqnarray*}
where  $R_{1,2}:=\frac{1}{N}\sum^N_{i=1} \sigma_i^1 \sigma^2_i$
is the normalized inner product between $\sigma^1$ and $\sigma^2$ and 
\begin{eqnarray}\label{eq:psxi}
\xi(x):=\lambda \cdot x^p+(1-\lambda) \cdot x^s.
\end{eqnarray}

The spherical $p+s$  spin glass model is an example of a mixed spherical $p$-spin glass model, where $\xi$ takes the form \eqref{eq:psxi}. The study of the partition function $Z_{N} = \int_{S_{N}} \exp(\beta H_{N}(\sigma)) d \sigma$ and the maximum 
\[
M_{N} = \max_{\sigma \in S_{N}} H_{N}(\sigma) 
\] 
has a long history in both physics and mathematics community. We refer the readers to the books of Mezard-Parisi-Virasoro \cite{mezard1987spin}, Talagrand \cite{TalagrandVolI}, the recent survey \cite{AMSurvey} and the numerous references therein.

Let $\mathcal{K}$ be the collection of all measures $\nu$ on $[0,1]$ with the form 
\begin{eqnarray*}
\nu(ds)=\mathbbm{1}_{[0,1)} \gamma(x)dx +\Delta \delta_{\{1\}}(dx),
\end{eqnarray*}
where $\gamma(x)$ is a nonnegative, right-continuous, and nondecreasing function on $[0,1)$, $\Delta>0,$ and $\delta_{\{1\}}$ is a Dirac measure at 1.  Define the Crisanti-Sommers \cite{CS} functional by 
\begin{eqnarray*}
\mathcal{Q}(\nu)=\frac{1}{2} \Big( \int^1_0  \xi'(x) \nu(dx) +\int^1_0 \frac{dx}{\nu \big( (x,1] \big)} \Big)
\end{eqnarray*}
for $\nu \in \mathcal{K}.$ The Parisi formula for the maximum energy (see \cite[Theorem 1]{ArnabChen15}) states that almost surely
\begin{eqnarray*}
\lim_{N \rightarrow \infty} N^{-1} M_N= \inf_{\nu \in \mathcal{K}} \mathcal{Q}(\nu).
\end{eqnarray*}
The functional $\mathcal{Q}$ is strictly convex on $\mathcal{K}$ and the right-hand side has a unique minimizer, denoted by 
\begin{eqnarray*}
\nu_P(dx)=\gamma_P(x) \mathbbm{1}_{[0,1)}(x)dx + \Delta_P \delta_{\{1\}}(dx).
\end{eqnarray*}
We denote by $\rho_P$ the measure on $[0,1)$ induced by $\gamma_P$, i.e.,
\begin{eqnarray*}
\gamma_P(x) = \rho_P([0,x]), \forall x \in [0,1),
\end{eqnarray*}
and we call $\nu_P$ the Parisi measure at zero temperature. At positive temperature, the Parisi formula reads as 
\begin{equation}\label{PFpt}
\lim_{N\to \infty} \frac{1}{N} \log Z_{N} = \inf_{\mu} \mathcal P_{\beta} (\mu),
\end{equation}
where the infimum in \eqref{PFpt} runs over the space of probability measures on $[0,1]$ and $\mathcal P_{\beta}$ is the Parisi functional described in \cite[Theorem 1.1]{Tal06} and \cite{CS}. The unique minimizer $\mu_{P}(\beta)$ of \eqref{PFpt} is called the Parisi measure at inverse temperature $\beta$.

In this paper we investigate the structure of the Parisi measure $\nu_{P}$ as a function of the parameters $p, s$ and $\lambda$ as well as the structure of the Parisi measure  $\mu_{P}(\beta)$ at low temperature. Our paper is inspired by the beautiful work of Crisanti and Leuzzi \cite{CL1, CL2} who predicted the results we describe.

Both in physics and mathematics, Parisi measures have been the main subject of study by several authors in the past decades. It is known that the structure of $\nu_{P}$  qualitatively describes the energy landscape, and it relates to efficiency of a large class of optimization algorithms \cite{el2021optimization, montanari2019optimization, subag2018following}.  Although the role and importance of the order parameter have been unveiled, the structure of the Parisi measure in many models remains very mysterious. Our main result in this paper is a complete characterization of the support of $\nu_{P}$ for all possible values of the parameters $p, s$ and $\lambda$. We describe our results now and put them into context at the end of this section. We start with a few definitions. 

\begin{definition}[k-RSB] For $k\geq0$, a Parisi measure $\nu_{P}$ has $k$ steps of replica symmetry breaking or, equivalently, the model is called $k$-RSB, if the  support of $\nu_{P}$ is discrete and has exactly $k+1$ points. When $k=0$ the model is also called Replica Symmetric (RS).

\end{definition}
\begin{definition}[FRSB] A Parisi measure $\nu_{P}$ has full steps of replica symmetry breaking or, equivalently, the model is FRSB, if $\nu_{P}$ is absolutely continuous on $[0,1)$ and $\nu_{P}([0,1)) >0$.
\end{definition}

\begin{definition}[1-FRSB] A Parisi measure $\nu_{P}$ is called 1-FRSB, if there exist $q>0, m\in(0,1)$ and a strictly increasing function $w(x)$ such that \begin{eqnarray*}
\nu(dx)=m \cdot \mathbbm{1}_{[0,q)}(x)dx+w(x)\cdot \mathbbm{1}_{[q,1)}(x)dx +\Delta \delta_{\{1\}}(dx),
\end{eqnarray*}
or
\begin{eqnarray*}
\nu(dx)=w(x)\cdot \mathbbm{1}_{[0,q)}(x)dx +m \cdot \mathbbm{1}_{[q,1)}(x)dx+\Delta \delta_{\{1\}}(dx).
\end{eqnarray*}
\end{definition}

\begin{definition}[2-FRSB] A Parisi measure $\nu_{P}$ is called 2-FRSB, if it there exist $0\leq q_{1}< q_{2}$, $m_{1}>0, m_{2}>0$ and a strictly increasing function $w(x)$ such that
\begin{eqnarray*}
\nu(dx)=m_1 \cdot \mathbbm{1}_{[q_{1},q_2)}(x)dx+w(x)\cdot \mathbbm{1}_{[q_1,q_2)} +m_2 \cdot \mathbbm{1}_{[q_2,1)}(x)dx+\Delta \delta_{\{1\}}(dx).
\end{eqnarray*}
\end{definition}

 Let $z$ be the unique solution to 
\begin{eqnarray}\label{eqnz}
\frac{1}{\xi'(1)}=\frac{1+z}{z^2} \log (1+z)-\frac{1}{z}
\end{eqnarray}
and define
\[ \bar{\lambda}_{1 \rightarrow 1F}:= \inf \{ \lambda \in [0,1] | 2 \lambda z \geq s(1-\lambda) \}.\] Note that when $\lambda=1,$ it holds that $2\lambda z \geq s(1-\lambda),$ which implies that $\bar{\lambda}_{1 \rightarrow 1F}$ is well-defined.

Last, set \[\bar{\lambda}_{1F \rightarrow F}: =\frac{s^2(s-1)}{(s-2)(s^2+s+6)}.\] Our first result classifies the Parisi measure at zero temperature in the case $p=2$.

\begin{theorem}\label{thm2+s}
For $s>2,$ the model $\xi(x)=\lambda x^2 +(1-\lambda)x^s$ is at zero temperature
\begin{eqnarray*}
&(i) &\text{ 1-RSB } \text{when } \lambda \in [0,\bar{\lambda}_{1 \rightarrow 1F}),\\
&(ii) &\text{ 1-FRSB }\text{when } \lambda \in [\bar{\lambda}_{1 \rightarrow 1F},\bar{\lambda}_{1F \rightarrow F}),\\
&(iii) &\text{ FRSB } \text{when } \lambda \in [\bar{\lambda}_{1F\rightarrow F} ,1), \\
&(iv) &\text{ RS }  \text{when } \lambda = 1.\
\end{eqnarray*}
\end{theorem}


The results for the case $p\neq 2$ are more elaborate to state and the phase diagram will depend on the values of $p$ and $s$. Define 
\begin{eqnarray*}
\Psi(\lambda):=-\frac{\xi'(1)-1}{\xi''(1)}-\log \Big (\frac{\xi''(1)}{\xi'(1)}  \Big ) +1-\frac{2}{\xi'(1)}+\frac{\xi''(1)}{\xi'(1)^2},
\end{eqnarray*}
and let
\begin{eqnarray}\label{lambda2f}
 Q(\lambda) &=&(s-p)^2(s^2-3s+2+p^2+3ps-3p)\lambda^2 \nonumber \\
&&-s(s-p)(2s^2-6s+4-p^2+3ps-3p)\lambda+s^2(s-1)(s-2).
\end{eqnarray}
The function $Q$ is a polynomial of degree two in $\lambda$ with discriminant given by 
\[
\Delta:=s^2-6(p-1)s+(p-1)(p+7).
\]
If $\Delta >0$ we let  $\lambda^*_1 < \lambda^*_2$ be the two real roots of $Q$. Our next theorem provides an exact condition for the mixed $p+s$ model to be 1RSB for all values of $\lambda \in (0,1)$.

\begin{theorem}[1RSB $p+s$ model]\label{coro1rsb}
Let $s\geq p>2$. The model $\xi(x)=\lambda x^p +(1-\lambda)x^s$ is 1RSB for all $\lambda\in[0,1]$ if and only if one of the following two conditions are satisfied:
\begin{eqnarray*}
&(i)& \Delta \leq  0, \\
&(ii)& \Delta > 0 \text{ and } \Psi(\lambda^*_1)<0.
\end{eqnarray*}
\end{theorem}

\begin{figure}[H] 
\centering 
\includegraphics[width=0.4\textwidth]{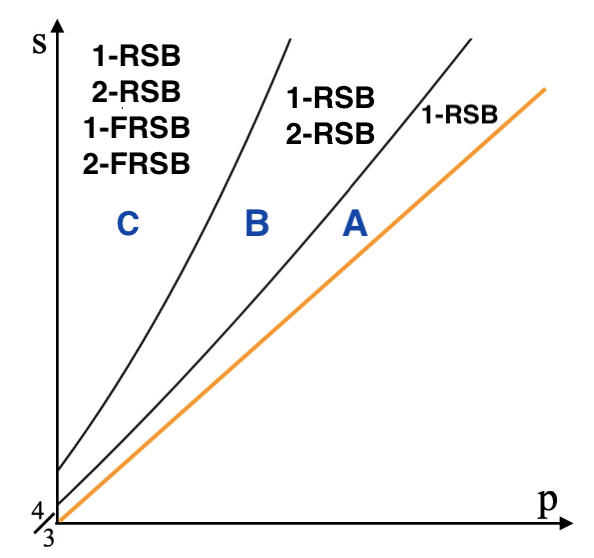} 
\caption{Different phases for the spherical $p+s$ model with $s\geq p > 2$ at zero temperature. In orange the line $p=s$.  Regions A, B, and C correspond to the values of $s$ and $p$ that satisfy Theorems \ref{coro1rsb}, \ref{coro12rsb},  and \ref{coro12frsb}, respectively. The curve plotted separating regions $A$ and $B$ is given by the equation $\Psi(\lambda^*_1)=0$. The curve that separates regions $B$ and $C$ is given by  $\Psi(\lambda_{2 \rightarrow 1F})=0$.} 
\label{Fig.main2}
\end{figure}

\begin{remark}
As seen in Figure \ref{Fig.main2} (see Region A), the conditions in Theorem \ref{coro1rsb} are satisfied and correspond to the cases when the difference $s-p$ is small.
\end{remark}

Now we turn to the $p+s$ spherical models that contain a phase that is not 1-RSB. Let \begin{eqnarray*}\label{lambda2f}
S(\lambda) &=& s^3(s-1)^2(s-2)(1-\lambda)^2+p^3(p-1)^2(p-2)\lambda^2 \nonumber \\
&&-2p(p-1)s(s-1)\Big[(p-2)(p-3)+(s-2)(s-3)-3(p-2)(s-2) \Big]\lambda(1-\lambda).
\end{eqnarray*}

If $\Delta>0,$ then straight-forward algebra shows that $S(\lambda)$ also has two real roots that we denote  by $\lambda_{2 \rightarrow 1F}<\lambda'_{2 \rightarrow 1F}.$  
Two possible scenarios may occur for the spherical  $p+s$ model depending on the sign of  $\Psi(\lambda_{2 \rightarrow 1F})$.

\begin{theorem}[1RSB/2RSB]\label{coro12rsb}
Let $s-1>p>2$. If
\begin{eqnarray*}
\Delta>0, \Psi(\lambda^*_1)>0 \text{ and } \Psi(\lambda_{2 \rightarrow 1F})\leq 0
\end{eqnarray*}
then there exist constants $0< \lambda_{1 \rightarrow 2} < \lambda_{2 \rightarrow 1} <1$ such that the model $\xi(x)=\lambda x^p +(1-\lambda)x^s$ is at zero temperature
\begin{eqnarray*}
&(i)&\text{ 1-RSB when } \lambda \in [0, \lambda_{1 \rightarrow 2}]\cup[\lambda_{2 \rightarrow 1},1],\\
&(ii)&\text{ 2-RSB when } \lambda \in (\lambda_{1 \rightarrow 2}, \lambda_{2 \rightarrow 1}).
\end{eqnarray*}

\end{theorem}

The next theorem describes the phase diagram when $\Psi(\lambda_{2 \rightarrow 1F})> 0$.

\begin{theorem}[1RSB/2RSB/1FRSB/2FRSB]\label{coro12frsb}
If
\begin{equation}\label{hyp34}
\Delta>0, \Psi(\lambda^*_1)>0 \text{ and } \Psi(\lambda_{2 \rightarrow 1F})> 0.
\end{equation}
then there exist constants $0< \lambda_{1 \rightarrow 2} < \lambda_{2 \rightarrow 2F} < \lambda _{2 \rightarrow 1F}< \lambda_{2 \rightarrow 1} <1$ such that the model $\xi(x)=\lambda x^p +(1-\lambda)x^s$ is at zero temperature
\begin{eqnarray*}
&(i)&\text{ 1-RSB when } \lambda \in [0, \lambda_{1 \rightarrow 2}]\cup[\lambda_{2 \rightarrow 1},1], \\
&(ii)&\text{ 2-RSB when } \lambda \in [\lambda_{1 \rightarrow 2}, \lambda_{2 \rightarrow 2F}],\\
&(iii)& \text{ 2-FRSB when } \lambda \in [\lambda_{2 \rightarrow 2F},\lambda _{2 \rightarrow 1F}], \\
&(iv)& \text{ 1-FRSB when } \lambda \in [\lambda_{2 \rightarrow 1F},\lambda _{2 \rightarrow 1}]. 
\end{eqnarray*}
\end{theorem}

\begin{remark} The constants $\lambda_{1 \rightarrow 2} < \lambda_{2 \rightarrow 2F} < \lambda _{2 \rightarrow 1F}< \lambda_{2 \rightarrow 1}$ in Theorems \ref{coro12rsb} and \ref{coro12frsb} are explicit and defined as solutions of a system of equations, see Section \ref{def:constant}, Equation \ref{eqnh1}.
\end{remark}
\begin{remark} Contrary to the mixed $p$-spin model on the hypercube $\{\pm 1\}^{N}$ where FRSB phases are expected \cite{auffinger2017sk}, there are no FRSB phases in the $p+s$ spherical model if $p, s>2$. When a 1-FRSB phase emerges, a 2-FRSB also appear. 
\end{remark}

The figure below summarizes our main results from Theorems \ref{coro1rsb}, \ref{coro12rsb},  amd \ref{coro12frsb} with the three possible scenarios for $p+s$ models with $s\geq p > 2$. 

A consequence of Theorem \ref{coro12frsb} we can show that once either $1$-FRSB or $2$-FRSB phases exist, the number of levels of replica symmetry breaking must diverge as $\beta$ goes to infinity.

\begin{corollary}\label{cor:finiteTemp} Assume $p,s$ satisfy the hypothesis of Theorem \ref{coro12frsb}. If $\lambda \in [\lambda_{2 \rightarrow 2F},\lambda _{2 \rightarrow 1F}]$ then for any $k\geq 1$  there exists $\beta_{k}$ such that the spherical $p+s$-spin model is at least $k$-RSB for all $\beta > \beta_{k}.$\end{corollary}

We conclude this section with some historical remarks and describing the main novelty of our approach. Parisi's formula \cite{parisi1979infinite} is a landmark result of Talagrand \cite{talagrand2006parisi} after major contribution from Guerra \cite{guerra2003broken}. At zero temperature, the formula was discovered in \cite{AC16}. The Parisi measure has been a topic of study for the past four decades. First rigorous results were obtained by Aizenman, Lebowitz, and Ruelle \cite{ALR} who showed that the SK model on $\{\pm 1\}^{N}$ is RS for $\beta <1$, while Toninelli \cite{Toni} showed that the Parisi measure is not RS in the low temperature region $\beta>1$. On $\{\pm 1\}^{N}$, for the SK model, the whole region $\beta >1$ is expected to be of FRSB. This prediction extends to all mixed $p$-spin models at sufficiently low temperature. After progress in \cite{talagrand2006parisi-b, auffinger2015parisi, AChen15PTRF}, the best result in this direction appeared in \cite{auffinger2017sk}, where for any mixed $p$-spin on $\{\pm 1\}^{N}$, it was shown that at zero temperature the model is not $k$-RSB for any $k\geq 1$ and an analogue of Corollary \ref{cor:finiteTemp}. Regularity properties of Parisi measures can be found in \cite{AChen15PTRF, JT, AC18Advances}. 

For the spherical model, the structure of the Parisi measure is different than that of $\{\pm 1\}^{N}$ at low temperature. The spherical pure $p$-spin is known to be $1$-RSB at low and zero temperature \cite{Tal06, JT, ArnabChen15}. Examples of $2$-RSB models at zero temperature first appeared in \cite{AZeng}, including some $p+s$ models. As far as we know, our main result is the first to characterize the full phase diagram for a large collection of models beyond the pure $p$-spin.

The main novelty of our approach is to explore the Parisi formula for the ground state energy directly at zero temperature and derive new and improved equivalent criteria to describe its support. As an example, Theorem \ref{2RSB} is a simple but powerful improvement of the results obtained in \cite{AZeng} while the conditions for 1-FRSB and 2-FRSB presented in Sections \ref{sec2frsb} and \ref{sec1frsb} are sharp and amenable to calculations. 

\subsection{Acknowledgments} A.A.'s research is partially supported by NSF Grant CAREER DMS-1653552 and NSF Grant DMS-1517894. A.A. would like to thank Andrea Crisanti and Luca Leuzzi for several helpful discussions and for sharing their own work that predicted the phase diagram for the $p+s$ spherical model \cite{CL2}. Y.Z.'s research is also partially supported by NSF Grant CAREER DMS-1653552.

\section{Boundaries and examples}\label{def:constant}

In this section we provide the definition of the constants $\lambda_{1 \rightarrow 2}, \lambda_{2 \rightarrow 2F}, \lambda_{2 \rightarrow 1}$ that appear in Theorems \ref{coro12rsb} and \ref{coro12frsb}, their existence will be shown in the proof of the aforementioned theorems. We also provide examples of mixtures that lie in each of the domains. 

Let 
\begin{equation}
h^1_1(x,\lambda):=\frac{x\xi'(1)-\xi'(x)}{\xi'(1)-\xi'(x)}+\frac{[\xi'(1)x-\xi'(x)]^2\xi(x)}{\xi'(1)\xi'(x)(\xi'(1)-\xi'(x))x(1-x)} -\log \Big(\frac{\xi'(1)-\xi'(x)}{\xi'(1)(1-x)} \Big)
\end{equation}
and
\begin{equation}
h^1_2(x,\lambda):=(1-x)[\xi'(1)-\xi'(x)] \Big ( \frac{\xi'(1)\xi'(x)x}{[\xi'(1)x-\xi'(x)]^2}\log\big(\frac{\xi'(1)x}{\xi'(x)}\big) -\frac{\xi'(x)}{\xi'(1)x-\xi'(x)} \Big) +\xi'(x)(1-x)-1+\xi(x).
\end{equation}

Consider the following system of equations with respect to $x$ and $\lambda:$
\begin{equation}\label{eqnh1}  \left\{
\begin{array}{lcl}
h^1_1(x, \lambda)&=&0 ,   \\
h^1_2(x, \lambda)&=&0.
\end{array} \right. \end{equation}
In the proof of Theorem \ref{coro12rsb} we will show that, under the assumptions oof $\Delta>0, \Psi(\lambda^*_1)>0$, a solution of the above system exists, and it is unique. We will call this solution $\lambda_{1 \rightarrow 2}.$

Similarly, let 
\begin{equation}h_1^2(x,\lambda):= \frac{x\xi''(x)-\xi'(x)}{x\xi''(x)} -\log\Big( \frac{x\xi''(x)}{\xi'(x)} \Big) +\frac{[\xi'(x)-x\xi''(x)]^2\xi(x)}{\xi'(x)^2 \xi''(x) x^2}
\end{equation}
and
\begin{align}
h^2_2(x,\lambda)=&- \xi(1)+\xi(x)+\xi'(x)(1-x)- \frac{\xi''(x)[\xi'(1)-\xi'(x)](1-x)^2}{\xi'(1)-\xi'(x)-\xi''(x)(1-x)} 
\\&+\frac{\xi''(x)[\xi'(1)-\xi'(x)]^2(1-x)^2}{[\xi'(1)-\xi'(x)-\xi''(x)(1-x)]^2}\log \big (\frac{\xi'(1)-\xi'(x)}{\xi''(x)(1-x)} \big).
\end{align}

In Section \ref{sec:rest}, we prove that if \eqref{hyp34} holds then the system of equations
\begin{equation}\label{eqnh2}  \left\{
\begin{array}{lcl}
h^2_1(x,\lambda)&=&0 ,   \\
h^2_2(x,\lambda)&=&0. \\
\end{array} \right. \end{equation}
has a unique solution given by $(x^{*}, \lambda_{2 \rightarrow 2F})$.

Last, if $\Delta>0$ and $\Psi(\lambda^*_1)>0$, we denote the zero of $\Psi(\lambda)$ between $\lambda^*_1$ and $\lambda^*_2$ by $\lambda_{2 \rightarrow 1}$.  It is not difficult to check that  $\lambda_{2 \rightarrow 1}$ exists and is unique under these assumptions.

We also give a concrete example for the phase characterization of the $p+s$ model. In the following example, we fix the value of $p$ and show the respective phase transitions of $\xi$ with different values of $s$.
\begin{example}
Let $p=4$ and $s=18.$ A computation shows that $\Delta=30>0$ and $\Psi(\lambda^*_1)<0$. Thus, by Theorem \ref{coro1rsb} the model $\xi = \lambda x^{4} + (1-\lambda)x^{18}$ is 1RSB for all values of $\lambda \in [0,1]$. 
\end{example}

\begin{example}
Let $p=4$ and $s=28$. In this case, $\Delta>0$,  $\Psi(\lambda^*_1)>0$ and $\Psi(\lambda_{2 \rightarrow 1F})<0$. Therefore the model $\xi$ is 1-RSB when $\lambda \in [0, \lambda_{1 \rightarrow 2}]\cup[\lambda_{2 \rightarrow 1},1],$ then 2RSB $\text{when } \lambda \in [\lambda_{1 \rightarrow 2}, \lambda_{2 \rightarrow 1}]$ and finally 1RSB when $\lambda \in [\lambda_{2 \rightarrow 1},1],$
\end{example}

\begin{example}
Fixing $p=4$ and $s=38,$ one can check that $p,s$ satisfy that $\Delta>0$, $\Psi(\lambda^*_1)>0$ and $\Psi(\lambda_{2 \rightarrow 1F})>0$. Therefore Theorem \ref{coro12frsb} applies.

\end{example}

The rest of the paper is organized as follows: Section \ref{secCr2+s} is devoted to study the $2+p$ case and prove Theorem \ref{thm2+s}. Section \ref{secCr} is devoted to prove Theorems \ref{coro1rsb}, \ref{coro12rsb} and \ref{coro12frsb}. To be more specific, Section \ref{sec1rsb} consists of the proof of Theorems \ref{coro1rsb}, Theorem \ref{coro12rsb}$(i)$ and Theorem \ref{coro12frsb}$(i)$. The proof of Theorem \ref{coro12rsb}$(i)$ and Theorem \ref{coro12frsb}$(i)$ is involved in Section \ref{sec2rsb}. Section \ref{sec2frsb} is devoted to prove Theorem \ref{coro12frsb}$(iii)$ while Section \ref{sec1frsb} is devoted to prove Theorem \ref{coro12frsb}$(iv)$.



\section{Proof of Theorem \ref{thm2+s}}\label{secCr2+s}

We now briefly recap previous results concerning Parisi measures that we will need. First, we state a general characterization of Parisi measures.
\begin{theorem}[Theorem 2 in \cite{ArnabChen15}]\label{criterion}
For $\nu \in \mathcal{K},$ let 
\[
g(u)=\int^1_u \Big( \xi'(t)-\int^t_0 \frac{dr}{\nu((r,1])^2} \Big)dt.\]
 Then $\nu$ is the Parisi measure for the model $\xi$ if and only if the following conditions are satisfied:
 \begin{eqnarray*}
 &(i)& \quad \xi'(1)=\int^1_0 \frac{dr}{\nu((r,1])^2} , \\
 &(ii)& \text{ the function g satisfies } \min_{u \in [0,1]} g(u) \geq 0 \\
 &(iii)& \text{ for } S:= \{u \in [0,1):g(u)=0 \}, \text{ one has } \rho(S)=\rho([0,1)). \text{ Here } \rho \text{ is the measure} \\
 && \text{ induced by } \gamma, i.e.  \rho([0,s])=\gamma(s).
\end{eqnarray*}
\end{theorem}

Second, we recall a characterization of 1-RSB measures at zero temperature.

\begin{theorem}[Theorem 5 in \cite{AC18Advances}]\label{Criterion1RSB}
The Parisi measure $\nu_P$ is 1RSB with supp $\rho_P=\{0\}$ if and only if $\zeta(x) \leq 0,\forall x \in [0,1].$
where 
\begin{eqnarray*}
\zeta(x):=\xi(x)+\xi'(x)(1-x)+\frac{\xi'(x)}{z}-\frac{(1+z)\xi'(1)}{z^2}\log\left(1+\frac{z\xi'(x)}{\xi'(1)}\right).
\end{eqnarray*}
In this case,
\begin{eqnarray*}
\nu_P(dx)&=&\frac{z}{\sqrt{(1+z)\xi'(1)}}\cdot \mathbbm{1}_{[0,1)}(x)dx+\frac{1}{\sqrt{(1+z)\xi'(1)}} \cdot \delta_{\{1\}}(dx) \\
ME&=&\frac{\xi'(1)+z}{\sqrt{(1+z)\xi'(1)}}
\end{eqnarray*}
\end{theorem}

\begin{proof}[Proof of Theorem \ref{thm2+s}$(i)$]
By Theorem \ref{Criterion1RSB}, it suffices to show that $\lambda$ and $z$ satisfies the relation $2 \lambda z \leq s(1-\lambda)$ if and only if $\zeta(x) \leq 0 , \forall x \in [0,1]$.

By computation, we obtain that 
\begin{eqnarray*}
\zeta'(x)&=&\frac{\xi''(x)}{\xi'(1)+z\xi'(x)}  [ -z \cdot x\xi'(x)+(z+1)\xi'(x)-\xi'(1)x] \\
&=&\frac{x\xi''(x)}{\xi'(1)+z\xi'(x)} [-sz(1-\lambda)x^{s-1}+s(z+1)(1-\lambda)x^{s-2}-2z\lambda x+2\lambda z-s(1-\lambda)]
\end{eqnarray*}
Note that $\zeta(0)=\zeta(1)=0$ and $\zeta'(0)=\zeta'(1)=0.$

If $\zeta(x) \leq 0 , \forall x \in [0,1],$ combing with the fact that $\zeta(0)=\zeta'(0)=0$, we conclude that
\begin{eqnarray*}
\zeta''(0)=\frac{\xi''(0)}{\xi'(1)} \cdot [(z+1) \xi''(0)-\xi'(1)]=\frac{\xi''(0)}{\xi'(1)} \cdot [2z  \lambda-s(1-\lambda)] \leq 0
\end{eqnarray*}
which implies that $2z\lambda \leq s(1-\lambda).$ Now we prove the opposite direction.

By Descartes' rule of signs and the assumption $2\lambda z-s(1-\lambda)<0$, $\zeta'(x)=0$ has either 2 or no strictly positive roots. Since $\zeta'(1)=0,$ $\zeta'(x)$ has 2 strictly positive roots. Also since $\zeta(0)=\zeta(1)=0$, by mean value theorem, there exists $c \in (0,1)$ such that $\zeta'(c)=0.$ Hence we conclude that $\zeta(x)$ has exactly one critical point in $(0,1).$

Then by the assumption $2z  \lambda<s(1-\lambda)$, we obtain that $\zeta''(0) \leq 0$, which yields the desired conclusion that $\zeta(x)<0, \forall x \in (0,1). $\end{proof}

We now turn into the proof of Theorem \ref{thm2+s}$(ii)$.

\begin{proposition}\label{bar}
For $s \geq 4$ and $\lambda \in (\bar{\lambda}_{1 \rightarrow 1F},1), $ $\lambda$ and $z$ satisfy that $2\lambda z>s(1-\lambda).$ Moreover,  it holds that $\bar{\lambda}_{1 \rightarrow 1F}<\bar{\lambda}_{1F \rightarrow F}.$ 
\end{proposition}

\begin{proof}

Recall the definition $\bar{\lambda}_{1 \rightarrow 1F}:= \inf \{ \lambda \in [0,1] | 2 \lambda z \geq s(1-\lambda) \}$. In order to show that $2 \lambda z>s(1-\lambda)$ for $\lambda \in [\lambda_{1 \rightarrow 1F},1]$, it's equivalent for us to show that 
\begin{eqnarray*}
\frac{d}{d\lambda} \Big[ 2z\lambda-s(1-\lambda) \Big] |_{2z\lambda=s(1-\lambda)}>0.
\end{eqnarray*}

Since $z$ satisfies that $\frac{1}{\xi'(1)}=\frac{1+z}{z^2} \log(1+z)-\frac{1}{z}$, it yields that
\begin{eqnarray*}
\frac{s-p}{\xi'(1)^2}=\Big[ -\frac{z+2}{z^3} \log(1+z)+\frac{2}{z^2} \Big] \cdot z'_{\lambda},
\end{eqnarray*}
which implies that
\begin{eqnarray}\label{x^*3}
z'_{\lambda}=-\frac{(s-p)z(1+z)}{\xi'(1)[z+2-\xi'(1)]}.
\end{eqnarray}

Plugging the relation $z=\frac{s(1-\lambda)}{2\lambda}$ and $p=2$ in $\eqref{x^*3}$, we obtain that
\begin{eqnarray*}
z'_\lambda &=& -\frac{s-2}{\xi'(1)} \cdot \frac{s(1-\lambda)}{2\lambda}\cdot \frac{\xi'(1)}{2\lambda}\cdot \frac{2\lambda}{s-2(s-2)\lambda}=-\frac{s(s-2)}{2\lambda[s-2(s-2)\lambda]}
\end{eqnarray*}
which implies that
\begin{eqnarray*}
\frac{d}{d\lambda} \Big[ 2z\lambda-s(1-\lambda) \Big] |_{2z\lambda=s(1-\lambda)}&=&2z+2\lambda z'_\lambda+s\\
&=&\frac{s(1-\lambda)}{\lambda}-{s(s-2)}{s-2(s-2)\lambda}\\
&=&\frac{s \big [s-3(s-2)\lambda\big]}{\lambda \big[s-2(s-2)\lambda\big]}
\end{eqnarray*}
Now we claim that $\lambda<\frac{s}{3(s-2)}$ for $s >3$, which yields that $s-2(s-2)\lambda>0$ and then $\frac{d}{d\lambda} \Big[ 2z\lambda-s(1-\lambda) \Big] |_{2z\lambda=s(1-\lambda)}>0.$

By \eqref{eqnz} and an elementary inequality $\log(1+z) >\frac{3z(2+z)}{z^2+6z+6}$ for $z>0$, we obtain that
\begin{eqnarray}\label{2+sineq1}
\frac{1}{\xi'(1)}>\frac{1+z}{z^2}\cdot \frac{3z(2+z)}{z^2+6z+6}-\frac{1}{z}=\frac{2z+3}{z^2+6z+6}
\end{eqnarray}
By computation, the following inequalities are equivalent:
\begin{eqnarray*}
&&\eqref{2+sineq1}\\
&\Longleftrightarrow&z^2+6z-2\xi'(1)z+6-3\xi'(1)>0\\
&\Longleftrightarrow& \frac{(1-\lambda)\big[ 4(s-2)(s-3)\lambda^2-s(5s-12)\lambda+s^2 \big]}{4 \lambda^2}>0\\
&\Longleftrightarrow&4(s-2)(s-3)\lambda^2-s(5s-12)\lambda+s^2>0
\end{eqnarray*}

In order to show that $\lambda<\frac{s}{3(s-2)},$ it suffices for us to show that $4(s-2)(s-3)\lambda^2-s(5s-12)\lambda+s^2<0$ for $\lambda =\frac{s}{3(s-2)}$ and 1.
Indeed, by computation, for $\lambda=\frac{s}{3(s-2)},$
\begin{eqnarray*}
4(s-2)(s-3)\lambda^2-s(5s-12)\lambda+s^2=-\frac{2s^2(s-3)}{9(s-2)}<0
\end{eqnarray*}
and for $\lambda=1,$
\begin{eqnarray*}
4(s-2)(s-3)\lambda^2-s(5s-12)\lambda+s^2=-8(s-3)<0.
\end{eqnarray*}

Now we prove that $\bar{\lambda}_{1 \rightarrow 1F}< \bar{\lambda}_{1F \rightarrow F}$. By the definition of $\bar{\lambda}_{1 \rightarrow 1F},$ it suffices to show that for $\lambda=\bar{\lambda}_{1 \rightarrow 1F},$
\begin{eqnarray}
2z \cdot \frac{s^2(s-1)}{(s-2)(s^2+s+6)} \geq s \Big( 1- \frac{s^2(s-1)}{(s-2)(s^2+s+6)} \Big)
\end{eqnarray}
which can be simplified as
\begin{eqnarray*}
z \geq \frac{2(s-3)}{s(s-1)}.
\end{eqnarray*}

For $\lambda=\bar{\lambda}_{1F \rightarrow F},$ we obtain that
\begin{eqnarray*}
\xi'(1)=2 \cdot \frac{s^2(s-1)}{(s-2)(s^2+s+6)} +s \cdot \frac{s^2(s-1)}{(s-2)(s^2+s+6)}=\frac{2s(s+3)}{s^2+s+6}.
\end{eqnarray*}

Since $c(x):=\frac{1+x}{x^2}\log(1+x)-\frac{1}{x}$ is decreasing on $[0,+\infty)$ and $c(z)=\frac{1}{\xi'(1)}$, it is equivalent for us to show that 
\begin{eqnarray*}
\frac{s^2+s+6}{2s(s+3)}=\frac{1}{\xi'(1)} \leq c\Big(\frac{2(s-3)}{s(s-1)}\Big)
\end{eqnarray*}

Because of the elementary inequality 
\begin{eqnarray}\label{logz}
\log(1+z) \geq \frac{2z}{2+z},
\end{eqnarray}
it holds that $c(x) \geq \frac{1+x}{x^2} \cdot \frac{2x}{2+x} -\frac{1}{x}=\frac{1}{2+x}$. Thus for $x=\frac{2(s-3)}{s(s-1)},$ we obtain that
\begin{eqnarray*}
c\Big( \frac{2(s-3)}{s(s-1)} \Big) \geq \frac{s(s-1)}{2(s^2-3)}.
\end{eqnarray*}
Moreover, for $s\geq5,$
\begin{eqnarray*}
\frac{s^2+s+6}{2s(s+3)}-\frac{s(s-1)}{2(s^2-3)} \leq \frac{-s^3+6s^2-3s-18}{2s(s+3)(s^2-3)}=\frac{-(s-3)(s^2-3s-6)}{s(s+3)(s^2-3)}<0
\end{eqnarray*}
which implies that $c\Big( \frac{2(s-3)}{s(s-1)} \Big) > \frac{s^2+s+6}{2s(s+3)}$. Therefore for $s \geq 5,$ $\bar{\lambda}_{1 \rightarrow 1F}<\frac{s^2(s-1)}{(s-2)(s^2+s+6)}.$ 

Now for $s=4,$ we claim that when $\lambda=\frac{s^2(s-1)}{(s-2)(s^2+s+6)}=\frac{12}{13},$ $z \geq \frac{2(s-3)}{s(s-1)}=\frac{1}{6}$. Indeed, it holds that $z > \frac{1}{5}$. Thus for $s=4,$ it also holds that $\bar{\lambda}_{1 \rightarrow 1F}<\bar{\lambda}_{1F \rightarrow F}.$

\end{proof}

\begin{lemma}\label{ms}
Define 
\begin{eqnarray*}
m(x)=\xi'''(x)[\xi'(1)-\xi'(x)](1-x)-2\xi''(x)[\xi'(1)-\xi'(x)-\xi''(x)(1-x)]
\end{eqnarray*}
For $0<\lambda< \bar{\lambda}_{1F \rightarrow F},$ there exists $c_P\in (0,1)$ such that $m(x)<0$ for $x \in (0,c_P)$ and $m(x)>0$ for $x \in (c_P,1).$ Also, for $\bar{\lambda}_{1F \rightarrow F} \leq \lambda \leq 1$, it holds that $m(x)<0$ for $x \in (0,1).$
\end{lemma}
\begin{proof}
By computation,
\begin{eqnarray*}
m(x)
&=&-s^2(s-1)(s-2)(1-\lambda)^2x^{2s-3}+s^3(s-1)(1-\lambda)^2x^{2s-4}+2s(s-2)(s-3)\lambda(1-\lambda)x^{s-1}\\
&&-s(s-1)(1-\lambda)[s^2(1-\lambda)+4(s-3)\lambda]x^{s-2} \\
&&-4s\lambda(1-\lambda) +s(s-1)(s-2)(1-\lambda)[2\lambda+s(1-\lambda)]x^{s-3}.
\end{eqnarray*}
By Descartes' rule of signs, $m(s)$ has either 0, 2, or 4 strictly positive roots counting multiplicity. Note that $m(1)=m'(1)=m''(1)=0$ and \begin{eqnarray*}
m'''(1)&=&2\xi''(1) \xi''''(1)-3 \xi'''(1)^2\\
&=&s(s-1)(s-2)^2(s^2+s+6)(1-\lambda) \Big [\lambda-\frac{s^2(s-1)}{(s-2)(s^2+s+6)} \Big].
\end{eqnarray*}
Then for $\lambda \neq \bar{\lambda}_{1F \rightarrow F},$ $x=1$ is a zero of $m(x)$ with multiplicity 3, which implies that $m(x)$ has exactly 1 strictly positive root except $x=1.$

When $0< \lambda < \bar{\lambda}_{1F \rightarrow F}$, it holds that $m'''(1)<0,$ which implies that there exists $\delta>0$ such that for $x \in (1-\delta,1), m(x)>0.$ Combining with the fact that $m(0)<0,$ there exists $c_P \in (0,1)$ such that $m(c_P)=0.$ Since $m(x)$ has only 1 strictly positive positive root except $x=1,$ $c_P$ is the only zero of $m(x)$ in $(0,1).$

Thus for $0< \lambda < \bar{\lambda}_{1F \rightarrow F}$, we conclude that when $0<x<c_P,$ $m(x)<0$ and when $c_P<x<1,$ $m(x)>0.$

When $\bar{\lambda}_{1F \rightarrow F}<  \lambda<1$, it holds that $m'''(1)>0,$ which implies that there exists $\delta>0$ such that for $x \in (1,1+\delta), m(x)>0.$ Since $m(x) \rightarrow -\infty$ as $x \rightarrow + \infty,$ the only strictly positive zero except 1 is larger than 1. Since $m(0)<0$ and $m(1)=0,$ we conclude that $m(x)<0$ for $x \in (0,1).$\end{proof}

\begin{lemma}\label{hs}
Define $h(x):=-h^2_2(x,\lambda),$ then for  $\bar{\lambda}_{1 \rightarrow 1F}<\lambda< \bar{\lambda}_{1F \rightarrow F},$ h(x) has exactly 1 root $q_P$ in $(0,1)$ with $q_P \in (0,c_P)$. Also for $ \bar{\lambda}_{1F \rightarrow F}<\lambda<1,$ h(x) has no solution in $(0,1).$
\end{lemma}

\begin{proof}
By computation,
\begin{eqnarray*}
h'(x)&=&\frac{[\xi'(1)-\xi'(x)]^2(1-x)^2}{[\xi'(1)-\xi'(x)-\xi''(x)(1-x)]^3}\cdot  \Big[ \xi'''(x)-\frac{2\xi''(x) \Big( \xi'(1)-\xi'(x)-\xi''(x)(1-x) \Big) }{\big( \xi'(1)-\xi'(x) \big)(1-x)} \Big ] \cdot \\
&&\Big[ 2 \Big( \xi'(1)-\xi'(x)-\xi''(x)(1-x) \Big) - \Big (  \xi'(1)-\xi'(x)+\xi''(x)(1-x) \Big) \log \Big( \frac{\xi'(1)-\xi'(x)}{\xi''(x)(1-x)} \Big) \Big]\\
&=&\frac{[\xi'(1)-\xi'(x)](1-x)}{[\xi'(1)-\xi'(x)-\xi''(x)(1-x)]^3} \cdot m(x) \cdot \\
&&\Big[ 2 \Big( \xi'(1)-\xi'(x)-\xi''(x)(1-x) \Big) - \Big (  \xi'(1)-\xi'(x)+\xi''(x)(1-x) \Big) \log \Big( \frac{\xi'(1)-\xi'(x)}{\xi''(x)(1-x)} \Big) \Big].
\end{eqnarray*}

By the inequality \eqref{logz} with setting $z=\frac{\xi'(1)-\xi'(x)-\xi''(x)(1-x)}{\xi''(x)(1-x)}$, we obtain that
\begin{eqnarray*}
\log \Big( \frac{\xi'(1)-\xi'(x)}{\xi''(x)(1-x)} \Big)>\frac{2 \Big[ \xi'(1)-\xi'(x)-\xi''(x)(1-x) \Big]}{\xi'(1)-\xi'(x)+\xi''(x)(1-x)} 
\end{eqnarray*}
which implies that 
\begin{eqnarray*}
 2 &\Big(& \xi'(1)-\xi'(x)-\xi''(x)(1-x) \Big) - \Big (  \xi'(1)-\xi'(x)+\xi''(x)(1-x) \Big) \log \Big( \frac{\xi'(1)-\xi'(x)}{\xi''(x)(1-x)} \Big) \\
&<&2 \left( \xi'(1)-\xi'(x)-\xi''(x)(1-x) \right) -\left (  \xi'(1)-\xi'(x)+\xi''(x)(1-x) \right)  \frac{2 \Big[ \xi'(1)-\xi'(x)-\xi''(x)(1-x) \Big]}{\xi'(1)-\xi'(x)+\xi''(x)(1-x)} \\
&=&0.
\end{eqnarray*}

Thus $h'(x)$ and $m(x)$ have opposite signs for $x \in (0,1).$

Note that the assumption $2 \lambda z \geq s(1-\lambda)$ is equivalent to $z \geq \frac{\xi'(1)}{\xi''(0)}-1.$ Since $x \mapsto \frac{1}{\xi'(1)}+\frac{1}{x}-\frac{1+x}{x^2}\log(1+x)$ is strictly increasing on $(0,+\infty),$ it implies that
\begin{eqnarray*}
h(0)&=& 1+\frac{\xi''(0)\xi'(1)}{\xi'(1)-\xi''(0)}-\frac{\xi''(0)\xi'(1)^2}{[\xi'(1)-\xi''(0)]^2}  \log \Big ( \frac{\xi'(1)}{\xi''(0)} \Big) \\
&=& \xi'(1) \Big[ \frac{1}{\xi'(1)}+\frac{1}{\frac{\xi'(1)}{\xi''(0)}-1}-\frac{(\frac{\xi'(1)}{\xi''(0)}-1)+1}{[\frac{\xi'(1)}{\xi''(0)}-1]^2} \log \Big (1+(\frac{\xi'(1)}{\xi''(0)}-1) \Big) \Big] \\
&\leq& \xi'(1) \Big[ \frac{1}{\xi'(1)}+\frac{1}{z}-\frac{1+z}{z^2} \log(1+z) \Big]=0.
\end{eqnarray*}

By Lemma \ref{ms}, for $0<\lambda< \bar{\lambda}_{1F \rightarrow F} ,$ there exists $c_P\in (0,1)$ such that $h'(x)>0$ for $x \in (0,c_P)$ and $h'(x)<0$ for $x \in (c_P,1).$ Since $h(0) \leq 0$ and $h(1)=0,$ $h(x)$ has only 1 zero $q_P $ in $(0,1)$ and $q_P$ lies in $(0,c_P)$ when $2\lambda z > s(1-\lambda).$ Also, for $\bar{\lambda}_{1F \rightarrow F} \leq \lambda \leq 1$, it holds that $h'(x)>0$ for $x \in (0,1).$ Since when $2\lambda z > s(1-\lambda)$,$h(0)<0$ and $h(1)=0$, we obtain that $h(x)<0$ for $x \in (0,1)$.

\end{proof}

\begin{lemma}\label{gs}
Define 
\begin{eqnarray*}
\bar{g}(x)&=&\xi(1)-\xi(x)-\xi'(q_P)(1-x)+\frac{\xi''(q_P)[\xi'(1)-\xi'(q_P)](1-q_P)}{\xi'(1)-\xi'(q_P)-\xi''(q_P)(1-q_P)} \cdot (1-x) \\
&&-\frac{\xi''(q_P)[\xi'(1)-\xi'(q_P)]^2(1-q_P)^2}{[\xi'(1)-\xi'(q_P)-\xi''(q_P)(1-q_P)]^2} \log \Big (1+\frac{\xi'(1)-\xi'(q_P)-\xi''(q_P)(1-q_P)}{\xi''(q_P)(1-q_P)^2}\cdot(1-x) \Big).
\end{eqnarray*}
Then for $\bar{\lambda}_{1 \rightarrow 1F}<\lambda< \bar{\lambda}_{1F \rightarrow F},$ it holds that $\bar{g}(x)>0$, when $x \in (q_P,1).$
\end{lemma}

\begin{proof}
By computation,
\begin{eqnarray*}
\bar{g}'(x)&=&-\xi'(x)+\xi'(q_P)-\frac{\xi''(q_P)[\xi'(1)-\xi'(q_P)](1-q_P)}{\xi'(1)-\xi'(q_P)-\xi''(q_P)(1-q_P)} \\
&&+\frac{1}{\Big[\xi''(q_P)(1-q_P)^2+[\xi'(1)-\xi'(q_P)-\xi''(q_P)(1-q_P)](1-x)\Big]}  \cdot \\
&& \qquad \frac{\xi''(q_P)[\xi'(1)-\xi'(q_P)]^2(1-q_P)^2}{\Big[\xi'(1)-\xi'(q_P)-\xi''(q_P)(1-q_P) \Big]  } \\
&&=\frac{1}{\Big[ \xi''(q_P)(1-q_P)^2+\big[ \xi'(1)-\xi'(q_P)-\xi''(q_P)(1-q_P) \big](1-x) \Big]} \cdot \\
&&\qquad \frac{a_s x^s  -a_{s-1} x^{s-1} +a_2x^2+a_1 x +a_0}{\big[ \xi'(1)-\xi'(q_P)-\xi''(q_P)(1-q_P)\big]},
\end{eqnarray*}
where 
\begin{eqnarray*}
a_p&=&s(1-\lambda) \big [  \xi'(1)-\xi'(q_P)-\xi''(q_P)(1-q_P) \big]^2,\\
a_{p-1}&=&s(1-\lambda) \big[ \xi'(1)-\xi'(q_P)-\xi''(q_P)(1-q_P) \big] \times \\
&&\Big[  \xi''(q_P)(1-q_P)^2+ \big[ \xi'(1)-\xi'(q_P)-\xi''(q_P)(1-q_P) \big] \Big],\\
a_2&=&2\lambda \big[ \xi'(1)-\xi'(q_P)-\xi''(q_P)(1-q_P) \big]^2\\
a_1 &=& \big[ \xi'(1)-\xi'(q_P)-\xi''(q_P)(1-q_P)\big] \times\\
&&\Big[ \xi''(q_P) \big(\xi'(1)-\xi'(q_P) \big)(1-q_P)-\xi'(q_P) \big[ \xi'(1)-\xi'(q_P)-\xi''(q_P)(1-q_P) \big] \\
&& -2\lambda \big[  \xi'(1)-\xi'(q_P)-\xi''(q_P)(1-q_P)+\xi''(q_P)(1-q_P)^2 \big] \Big] \\
a_0 &=& \big[ \xi'(1)-\xi'(q_P)-\xi''(q_P)(1-q_P)\big]\xi'(q_P) \times \\
&&\big[ \xi'(1)-\xi'(q_P)-\xi''(q_P)(1-q_P)+\xi''(q_P)(1-q_P)^2\big]\\
&&-\xi''(q_P)[\xi'(1)-\xi'(q_P)](1-q_P)\big[ \xi'(1)-\xi'(q_P)-\xi''(q_P)(1-q_P)+\xi''(q_P)(1-q_P)^2 \big]\\
&&+\xi''(q_P)[\xi'(1)-\xi'(q_P)]^2(1-q_P)^2.
\end{eqnarray*}
Also
\begin{eqnarray*}
\bar{g}''(x)=-\xi''(x)+\frac{\xi''(q_P)[\xi'(1)-\xi'(q_P)]^2(1-q_P)^2}{\Big[\xi''(q_P)(1-q_P)^2+\big[\xi'(1)-\xi'(q_P)-\xi''(q_P)(1-q_P)\big](1-x) \Big]^2}
\end{eqnarray*}
and
\begin{eqnarray*}
\bar{g}'''(x)=-\xi'''(x)+\frac{2\xi''(q_P)[\xi'(1)-\xi'(q_P)]^2(1-q_P)^2\big[\xi'(1)-\xi'(q_P)-\xi''(q_P)(1-q_P)\big]}{\Big[\xi''(q_P)(1-q_P)^2+\big[\xi'(1)-\xi'(q_P)-\xi''(q_P)(1-q_P)\big](1-x) \Big]^3}.
\end{eqnarray*}

Note that $\bar{g}(q_P)=\bar{g}(1)=\bar{g}'(q_P)=\bar{g}'(1)=0$, $\bar{g}''(q_P)=0$ and 
\begin{eqnarray*}
\bar{g}'''(q_P)=-\xi'''(q_P)+\frac{2\xi''(q_P) \big[ \xi'(1)-\xi'(q_P)-\xi''(q_P)(1-q_P) \big]}{\big[\xi'(1)-\xi'(q_P)\big](1-q_P)}=-m(q_P).
\end{eqnarray*}
Thus, by the mean value theorem, there exists $t_0 \in (q_P,1) $ such that $\bar{g}'(t_0)=0.$

Since $2 \lambda z \geq s(1-\lambda)$ with $\lambda< \bar{\lambda}_{1F \rightarrow F},$ by Lemma \ref{ms}, there exists $c_P\in (0,1)$ such that $m(x)<0$ for $x \in (0,c_P)$ and $m(x)>0$ for $x \in (c_P,1).$ Also by Lemma \ref{hs}, $0<q_P<c_P$, which implies that $m(q_P)<0$ and $\bar{g}'''(q_P)>0.$

By Descartes' rule of signs, $\bar{g}'(x)$ has at most 4 strictly positive roots counting multiplicity. Since $\bar{g}'(q_P)=\bar{g}''(q_P)=0,$ $q_P$ is a zero of $\bar{g}'(x)=0$ with multiplicity at least 2. Combining with the fact that $\bar{g}'(t_0)=\bar{g}'(1)=0,$ $\bar{g}'(x)$ has 4 strictly positive roots $x=q_P,t_0,1$  counting multiplicity.

Since $\bar{g}(q_P)=\bar{g}'(q_p)=\bar{g}''(q_P)=\bar{g}(1)=0$ and $\bar{g}'''(q_P)>0$, combining with the reasoning above that there is only 1 critical point between $q_P$ and 1, we obtain that $\bar{g}(x)>0$ for $x \in (q_P,1).$

\end{proof}

\begin{lemma}\label{increasing}
If $\lambda>\bar{\lambda}_{1F \rightarrow F}$, $\frac{\xi'''(x)}{\xi''(x)^{\frac{3}{2}}}$ is strictly increasing on $[0,1).$ If $\lambda<\bar{\lambda}_{1F \rightarrow F}$, $\frac{\xi'''(x)}{\xi''(x)^{\frac{3}{2}}}$ is strictly increasing on $[0,q_P].$
\end{lemma}

\begin{proof}
By computation,
\begin{eqnarray*}
\frac{d}{dx} \Big( \frac{\xi'''(x)}{\xi''(x)^{\frac{3}{2}}} \Big) &=&\xi''(x)^{\frac{5}{2}} \cdot \Big [ 2\xi''(x)\xi''''(x)-3\xi'''(x)^2 \Big] \\
&=& \xi''(x)^{\frac{5}{2}} \cdot \Big[4s(s-1)(s-2)(s-3) \lambda (1-\lambda)x^{s-4}-s^3(s-1)^2(s-2)(1-\lambda)^2x^{2s-6}      \Big]\\
&=& -s^3(s-1)^2(s-2)(1-\lambda)^2x^{s-4}  \xi''(x)^{\frac{5}{2}} \cdot \Big[ x^{s-2}-\frac{4(s-3)\lambda}{s^2(s-1)(1-\lambda)} \Big].
\end{eqnarray*}

Note that for $\lambda>\bar{\lambda}_{1F \rightarrow F},$
\begin{eqnarray*}
\frac{4(s-3)\lambda}{s^2(s-1)(1-\lambda)} >\frac{4(s-3)}{s^2(s-1)}\cdot \frac{\frac{s^2(s-1)}{(s-2)(s^2+s+6)}}{1-(\frac{s^2(s-1)}{(s-2)(s^2+s+6)})}=1.
\end{eqnarray*}
Therefore for $\lambda>\bar{\lambda}_{1F \rightarrow F}$ and $x \in [0,1),$ we obtain that $\frac{d}{dx} \Big( \frac{\xi'''(x)}{\xi''(x)^{\frac{3}{2}}} \Big)>0$.

Now we assume that $0<\lambda<\bar{\lambda}_{1F \rightarrow F}.$  In order to prove that $\frac{\xi'''(x)}{\xi''(x)^{\frac{3}{2}}}$ is strictly increasing on $[0,q_P],$ it is equivalent to show that $q_P<\Big( \frac{4(s-3)\lambda}{s^2(s-1)(1-\lambda)}   \Big)^{\frac{1}{s-2}}$. Recall $m(x)$ in Lemma \ref{ms} and now regard $m$ as a function with respect to $\lambda$ and $x$. We denote it by $m(\lambda,x).$ For $0<\lambda< \bar{\lambda}_{1F \rightarrow F}$, it holds that $m(\lambda,x)<0$ when $x \in (0,c_P)$ and $m(\lambda, x)>0,$  when $x \in (c_P,1).$ It suffices to show that $m \Big( \lambda, \big( \frac{4(s-3)\lambda}{s^2(s-1)(1-\lambda)} \big)^{\frac{1}{s-2}} \Big)>0$ for $0<\lambda< \bar{\lambda}_{1F \rightarrow F}.$ 

Now we set $x=\big( \frac{4(s-3)\lambda}{s^2(s-1)(1-\lambda)} \big)^{\frac{1}{s-2}} $ and then $\lambda=\frac{s^2(s-1)x^{s-2}}{s^2(s-1)x^{s-2}+4(s-3)}.$ It's enough for us to show that $m(\frac{s^2(s-1)x^{s-2}}{s^2(s-1)x^{s-2}+4(s-3)},x)>0$, for $x \in (0,+\infty).$

By plugging the expression of $\lambda$ in $m(x),$ we obtain that
\begin{eqnarray*}
m(x)&=&\frac{1}{[s^2(s-1)x^{s-2}+4(s-3)]^2} \cdot \Big \{ 8s^2(s-1)(s-2)^2(s-3)^2x^{2s-3}\\
&&-16s^3(s-1)(s-2)(s-3)^2x^{2s-4}+8s^3(s-1)^2(s-2)(s-3)x^{2s-5}\\
&&+16s^3(s-1)(s-2)(s-3)x^{s-2}+16s^2(s-1)(s-2)(s-3)^2x^{s-3} \Big \}\\
&=&\frac{8s^2(s-1)(s-2)(s-3)x^{s-3}}{[s^2(s-1)x^{s-2}+4(s-3)]^2} \cdot \\
&&\Big \{ (s-2)(s-3)s^s-2s(s-3)x^{s-1}+s(s-1)x^{s-2}-2sx+2(s-3) \Big \}
\end{eqnarray*}
By Descartes' rule of signs, $m(x)$ has at most 4 strictly positive zeros counting multiplicity.

Notice that for $n(x):=(s-2)(s-3)x^s-2s(s-3)x^{s-1}+s(s-1)x^{s-2}-2sx+2(s-3) ,$
\begin{eqnarray*}
n(1)&=&n'(1)=n''(1)=n'''(1)=0 \\
n''''(1)&=&2s(s-1)(s-2)(s-3)>0,
\end{eqnarray*}
which implies that $x=1$ is a zero of $m(x)$ with multiplicity 4.
 Therefore $x=1$ is the only strictly positive zero  of $m(x).$ Since $n(x) \mapsto + \infty$, as $x \mapsto \infty,$ we conclude that $m(x)\geq0$ for $x \in (0,+\infty).$ Therefore for $\lambda< \bar{\lambda}_{1F \rightarrow F}$, we obtain that $q_P<\Big( \frac{4(s-3)\lambda}{s^2(s-1)(1-\lambda)}   \Big)^{\frac{1}{s-2}}$ and $\frac{\xi'''(x)}{\xi''(x)^{\frac{3}{2}}}$ is strictly increasing on $[0,q_P].$

\end{proof}

Now we start the proof of Theorem \ref{thm2+s}$(ii)$ and $(iii)$ as follows:
\begin{proof}[Proof of Theorem \ref{thm2+s}$(ii)$,$(iii)$]

We first deal with the case when $ \bar{\lambda}_{1F \rightarrow F}< \lambda<1.$

Assume that $\nu$ takes the form
\begin{eqnarray*}
\nu(dx)=\frac{1}{2} \xi''(x)^{-\frac{3}{2}} \xi'''(x) \cdot \mathbbm{1}_{[0,1)}(x) dx+ \xi''(1)^{-\frac{1}{2}} \cdot  \delta_1(dx).
\end{eqnarray*}

By Lemma \ref{increasing}, when $\lambda>\bar{\lambda}_{1F \rightarrow F},$ $\xi''^{-\frac{3}{2}}\cdot\xi'''$ is strictly increasing on $[0,1)$, which implies that $\nu \in \mathcal{K}.$

By Theorem \ref{criterion}, it suffices for us to show that
\begin{eqnarray}\label{supp}
\int^1_u \Big( \xi'(x)-\int^x_0 \frac{dr}{\nu([r,1])^2} \Big) dx=0 \text{ for all } u \in [0,1),
\end{eqnarray}
and
\begin{eqnarray}\label{supp1}
 \int^1_0 \frac{dx}{\nu([x,1])^2}=\xi'(1)
\end{eqnarray}

By the definition of $\nu$, it holds that, for all $x \in [0,1),$
\begin{eqnarray*}
\nu([x,1])=\xi''(1)^{-\frac{1}{2}}+\int^1_x x(r) dr=\xi''(1)^{-\frac{1}{2}}+\int^1_x  \frac{\xi'''(r) }{2 \xi''(r)^{\frac{3}{2}}}  dr=\xi''(x)^{-\frac{1}{2}}.
\end{eqnarray*}
which implies that
\begin{eqnarray*}
\xi'(x)=\int^x_0 \frac{dr}{\nu([r,1])^2}.
\end{eqnarray*}
Then \eqref{supp} and \eqref{supp1} are verified. By the uniqueness of the Parisi measure, $\nu$ is the minimizer of $\mathcal{Q}.$

We now deal with the case when $\bar{\lambda}_{1 \rightarrow 1F}< \lambda< \bar{\lambda}_{1F \rightarrow F}$. 
Assume that $\nu$ takes the form
\begin{eqnarray*}
\nu(dx)=\frac{1}{2} \xi''(x)^{-\frac{3}{2}} \xi'''(x)  \mathbbm{1}_{[0,q_P)}(x)dx + a_P \mathbbm{1}_{[q_P,1)}(x)dx+\Delta_P \delta_{\{1\}}(dx)
\end{eqnarray*}
where
\begin{eqnarray*}
\Delta_P=\frac{\xi''(q_P)^{\frac{1}{2}} \cdot (1-q_P)}{ \xi'(1)-\xi'(q_P)} \text{ and } a_P=\frac{\xi'(1)-\xi'(q_P)-\xi''(q_P)(1-q_P)}{\xi''(q_P)^{\frac{1}{2}}[\xi'(1)-\xi'(q_P)](1-q_P)}.
\end{eqnarray*}
and $q_P$ is the unique zero $q_P$ of $h(x)$ in $(0,1)$. 

Note that
\begin{eqnarray*}
&&\frac{1}{2} \cdot \frac{\xi'''(q_P)}{\xi''(q_P)^{\frac{3}{2}} }-a_P\\
&=&\frac{1}{2} \cdot \frac{\xi'''(q_P)}{\xi''(q_P)^{\frac{3}{2}} }-\frac{\xi'(1)-\xi'(q_P)-\xi''(q_P)(1-q_P)}{\xi''(q_P)^{\frac{1}{2}}[\xi'(1)-\xi'(q_P)](1-q_P)}\\
&=&\frac{\xi'''(q_P)[\xi'(1)-\xi'(q_P)](1-q_P)-2\xi''(q_P)^{\frac{1}{2}} [\xi'(1)-\xi'(q_P)-\xi''(q_P)(1-q_P)] }{2\xi''(q_P)^{\frac{3}{2}}[\xi'(1)-\xi'(q_P)](1-q_P)}\\
&=&\frac{m(q_P)}{2\xi''(q_P)^{\frac{3}{2}}[\xi'(1)-\xi'(q_P)](1-q_P)}.
\end{eqnarray*}
Since Lemma \ref{ms} and Lemma \ref{hs}, we obtain that $q_P<c_P$ and $m(q_P)<0$, which implies that $a_P>\frac{\xi'''(q_P)}{2\xi''(q_P)^{\frac{3}{2}} }.$  Also by Lemma \ref{increasing}, $\xi''^{-\frac{3}{2}}\cdot \xi'''$ is strictly increasing on $[0,q_P].$\

Similarly, it suffices for us to show that $g(u)=0$ for $u \in [0,q_P),$ $g(u)>0$ for $u \in [q_P,1)$,
 and $ \int^1_0 \frac{dx}{\nu([x,1])^2}=\xi'(1).$

By computation, 
 \begin{equation*} \nu((r,1])= \left\{
\begin{array}{lcl}
\Delta_P+a_P(1-r)=\frac{\xi''(q_P)^{\frac{1}{2}}(r-q_P)}{\xi'(1)-\xi'(q_P)} + \frac{1-r}{\xi''(q_P)^{\frac{1}{2}}(1-q)}&\text{ for }& q_P \leq r \leq 1,    \\
\Delta_P+a_P(1-q_P) + \int^{q_P}_r x(s) ds=\xi''(r)^{-\frac{1}{2}}  &\text{ for }& 0 \leq r \leq q_P, \\
\end{array} \right. \end{equation*} 
and for $u \in [0,q_P],$
\begin{eqnarray*}
\int^{q_P}_u \Big( \xi'(x)-\int^s_0 \frac{dr}{\nu([r,1])^2} \Big) ds =\int^{q_P}_u \Big( \xi'(x)-\xi'(x) \Big) dx=0.
\end{eqnarray*}

Moreover, 
\begin{eqnarray*}
 \int^1_0 \frac{ds}{\nu([x,1])^2}&=& \int^{q_P}_0 \frac{dx}{\nu([x,1])^2}+ \int^1_{q_P} \frac{dx}{\nu([x,1])^2} \\
 &=& \int^{q_P}_0 \xi''(r)dr+ \int^1_{q_P} \frac{dr}{[\Delta_P+a_P(1-r)]^2}\\
 &=& \xi'(q_P)+\frac{1-q_P}{\Delta_P[\Delta_P+a_P(1-q_P)]}\\
 &=&\xi'(q_P)+\xi'(1)-\xi'(q_P)=\xi'(1).
\end{eqnarray*}

By computation, for $u \in [q_P,1),$
\begin{eqnarray*}
\int^1_u &\Big(& \xi'(x)-\int^x_0 \frac{dr}{\nu([r,1])^2} \Big) dx\\
&=&\xi(1)-\xi(u)-\xi'(q_P)(1-u)+\frac{1-u}{a_P[\Delta_P+a_P(1-q_P)]}-\frac{1}{a^2_P} \log(1+\frac{a_P}{\Delta_P}(1-u)) \\
&=&\xi(1)-\xi(u)-\xi'(q_P)(1-u)+\frac{\xi''(q_P)[\xi'(1)-\xi'(q_P)](1-q_P)(1-u)}{\xi'(1)-\xi'(q_P)-\xi''(q_P)(1-q_P)} \\
&&-\frac{\xi''(q_P)[\xi'(1)-\xi'(q_P)]^2(1-q_P)^2}{[\xi'(1)-\xi'(q_P)-\xi''(q_P)(1-q_P)]^2} \log \Big (1+ \frac{[\xi'(1)-\xi'(q_P)-\xi''(q_P)(1-q_P)](1-u)}{\xi''(q_P)(1-q_P)^2} \Big) \\
&=&\bar{g}(u).
\end{eqnarray*}
By Lemma \ref{gs}, for $\bar{\lambda}_{1 \rightarrow 1F} < \lambda< \bar{\lambda}_{1F \rightarrow F},$ it holds that $\bar{g}(x)>0$, when $x \in (q_P,1),$ and $\bar{g}(q_P)=0$, which verifies the condition $\int^1_u \Big( \xi'(x)-\int^x_0 \frac{dr}{\nu([r,1])^2} \Big) dx>0,$ for $u \in [q_P,1).$
Therefore $\nu$ is the minimizer of $\mathcal{Q}.$

\end{proof}

\section{The spherical $p+s$ model with $p,s>2$}\label{secCr}\label{sec:rest}

In this section, we prove the phase transition of the model $\xi$ when $s>p>2.$

\subsection{The 1-RSB phase at zero temperature }\label{sec1rsb}

We characterize the 1-RSB phase by the following proposition.
\begin{proposition}\label{Cr1rsb}
The model $\xi(x)$ is 1-RSB if and only if $h^1_1(q,\lambda)<0$, where $q$ satisfies $h^1_1(q)=\frac{\xi'(1)\xi'(q)[\xi'(1)-\xi'(q)]q(1-q)}{[\xi'(1)q-\xi'(q)]^2} \cdot  h^1_2(q,\lambda).$
\end{proposition}

\begin{proof}[Proof of Proposition \ref{Cr1rsb}]

By computation,
\begin{eqnarray*}
\zeta'(x)&=&-\frac{x\xi''(x)}{\xi'(1)+z\xi'(x)} \Big( \xi'(1)+z\xi'(x)-\frac{\xi'(x)}{x}(1+z)\Big) \\
&:=&-\frac{x\xi''(x)}{\xi'(1)+z\xi'(x)} \psi(x)
\end{eqnarray*}
and $$\psi'(1)=-  \xi''(1)+\xi'(1)(1+z)=\frac{1}{\xi'(1)} \Big (1+z-\frac{\xi''(1)}{\xi'(1) } \Big),$$

If the mode $\xi$ is 1-RSB, it holds that $\zeta(x) \leq 0$ for $x \in [0,1].$ Since $\zeta(1)=\psi(1)=0$, it yields that $\psi'(1)>0.$ Then we obtain that 
\begin{eqnarray}\label{x^*4}
1+z>\frac{\xi''(1)}{\xi'(1)}.
\end{eqnarray} 
By equation \eqref{eqnz} and the fact that $\frac{1+z}{z^2} \log (1+z)-\frac{1}{z}$ is strictly decreasing with respect to $z$, we have that
\begin{eqnarray}\label{h233}
\frac{\xi'(1)\xi''(1)}{[\xi''(1)-\xi'(1)]^2} \log \Big( \frac{\xi''(1)}{\xi'(1)} \Big) -\frac{\xi'(1)}{\xi''(1)-\xi'(1)}-\frac{1}{\xi'(1)}>0
\end{eqnarray}

Now if $x=q$ is a critical point of $\zeta(x),$ it yields that $\psi(q)=0$, which is equivalent to
\begin{eqnarray*}
z=\frac{q\xi'(1)-\xi'(q)}{\xi'(q)(1-q)}.
\end{eqnarray*}
Since $\zeta(0)=\zeta(1)=0,$ then the mean value theorem guarantees the existence of the critical point $q$ of $\zeta(x).$
Then the critical points $x=q$ of $\zeta(x)$ are the points satisfying that 
\begin{eqnarray}\label{criticalh3}
\frac{1+z(q)}{z(q)^2} \log(1+z(q)) -\frac{1}{z(q)}-\frac{1}{\xi'(1)}=0.
\end{eqnarray}
Note that $q$ is not necessarily unique.We denote all its possible solutions by $q$ in the following proof.

Then by setting $z=\frac{[\xi'(1)-\xi'(x)]x}{\xi'(x)(1-x)}-1,$  the equation \eqref{criticalh3} can be rewritten as
\begin{eqnarray*}
h_3(x):=\frac{\xi'(x)[\xi'(1)-\xi'(x)]x(1-x)}{[x\xi'(1)-\xi'(x)]^2} \log \Big( \frac{[\xi'(1)-\xi'(x)]x}{\xi'(x)(1-x)}  \Big) - \frac{\xi'(x)(1-x)} {x\xi'(1)-\xi'(x)}-\frac{1}{\xi'(1)}=0.
\end{eqnarray*}
Here by algebraic calculation, we notice that 
\begin{eqnarray*}
h_3(x)=\frac{\xi'(x)[\xi'(1)-\xi'(x)]x(1-x)}{[x\xi'(1)-\xi'(x)]^2} \cdot \Big(\frac{\xi'(1)\xi'(x)[\xi'(1)-\xi'(x)]x(1-x)}{[\xi'(1)x-\xi'(x)]^2} h^1_2(x,\lambda) -h^1_1(x)  \Big)
\end{eqnarray*}
and define $\tilde{h}^1_2(x)=\frac{\xi'(1)\xi'(x)[\xi'(1)-\xi'(x)]x(1-x)}{[\xi'(1)x-\xi'(x)]^2} h^1_2(x,\lambda) $.
Thus the critical points $x=q$ of $\zeta(x)$ satisfy that $\tilde{h}^1_2(q)=\frac{\xi'(1)\xi'(x)[\xi'(1)-\xi'(x)]x(1-x)}{[\xi'(1)x-\xi'(x)]^2} h^1_2(q,\lambda). $

We rewrite $\zeta(x)$ as follows:
\begin{eqnarray*}
\zeta(x)&=&\xi(x)+\xi'(x)(1-x)+\frac{\xi'(q)(1-q)}{\xi'(1)q-\xi'(q)}\xi'(x)\\
&&-\frac{\xi'(1)\xi'(q)[\xi'(1)-\xi'(q)]q(1-q)}{[q\xi'(1)-\xi'(q)]^2} \log\Big (1+\frac{[q\xi'(1)-\xi'(q)]\xi'(x)}{\xi'(1)\xi'(q)(1-q)} \Big )\\
&=&\frac{\xi'(1)\xi'(q)[\xi'(1)-\xi'(q)]q(1-q)}{[q\xi'(1)-\xi'(q)]^2} \cdot \\
&&\Big \{ \frac{[q\xi'(1)-\xi'(q)]^2\xi(x)} {\xi'(1)\xi'(q)[\xi'(1)-\xi'(q)]q(1-q)} +\frac{[q\xi'(1)-\xi'(q)]^2\xi'(x)(1-x)} {\xi'(1)\xi'(q)[\xi'(1)-\xi'(q)]q(1-q)}  \\
&&+\frac{[q\xi'(1)-\xi'(q)]}{\xi'(1)[\xi'(1)-\xi'(q)]q} \xi'(x) -\log\Big (1+\frac{[q\xi'(1)-\xi'(q)]\xi'(x)}{\xi'(1)\xi'(q)(1-q)} \Big )  \Big \}.
\end{eqnarray*}
In particular for $x=q$, we obtain that
\begin{eqnarray*}
\zeta(q)&=&\frac{\xi'(1)\xi'(q)[\xi'(1)-\xi'(q)]q(1-q)}{[q\xi'(1)-\xi'(q)]^2} \cdot \\
&&\Big \{ \frac{[q\xi'(1)-\xi'(q)]^2\xi(q)} {\xi'(1)\xi'(q)[\xi'(1)-\xi'(q)]q(1-q)} +\frac{q\xi'(1)-\xi'(q)} {\xi'(1)-\xi'(q)}    -\log\Big (\frac{[\xi'(1)-\xi'(q)]}{\xi'(1)(1-q)} \Big ) \Big\}\\
&=&\frac{\xi'(1)\xi'(q)[\xi'(1)-\xi'(q)]q(1-q)}{[q\xi'(1)-\xi'(q)]^2} \cdot h^1_1(q,\lambda)
\end{eqnarray*}


Now since the critical points $x=q$ of $\zeta(x)$ are the points satisfing that $\tilde{h}^1_2(q)=\frac{\xi'(1)\xi'(q)[\xi'(1)-\xi'(q)]q(1-q)}{[\xi'(1)q-\xi'(q)]^2} h^1_2(q,\lambda) $, the inequality $\zeta(x)<0, \forall x \in (0,1)$ is equivalent to $h^1_1(q,\lambda)<0,$ where $q$ satisfies that $h^1_1(q,\lambda)=\tilde{h}^1_2(q).$

\end{proof}

\begin{lemma}\label{tx}
Consider the function 
\begin{eqnarray}
t(x)=\xi'(1)\xi''(x)x(1-x)-\xi'(x)[\xi'(1)-\xi'(x)].
\end{eqnarray}
If $s^2+6s-7+p^2+6p-6ps<0,$ then for any $\lambda \in (0,1)$, it holds that $t(x)$ has either 0 or 2 roots in $(0,1).$

If $s^2+6s-7+p^2+6p-6ps>0,$ when $\lambda \in (\lambda^*_1,\lambda^*_2),$ $t(x)$ has exactly 1 root in $(0,1)$. For $0< \lambda<\lambda^*_1$, $t(x)$ has either zero or two roots in $(0,1)$.  For $\lambda> \lambda^*_2,$ $t(x)$ has no zeros in $(0,1).$
\end{lemma}

\begin{proof}
By computation, we obtain that
\begin{eqnarray*}
t(x)&=&\xi'(1)p(p-2)\lambda x^{p-1}-\xi'(1)p(p-1)\lambda x^p+p^2\lambda^2x^{2p-2}+\xi'(1)s(s-2)(1-\lambda) x^{s-1}\\
&&-\xi'(1)s(s-1)(1-\lambda) x^s+2ps \lambda(1-\lambda)x^{p+s-2}+s^2(1-\lambda)^2x^{2s-2}
\end{eqnarray*}

When $s \geq 2p,$ the coefficients of $t(x)$ have signs $(+,-,+,+,-,+,+)$ and when $s<2p$, they have signs $(+,-,+,-,+,+,+).$ For both of two cases, $t(x)$ has 4 changes of signs. Then by Descartes' rule of signs, $t(x)$ has at most four strictly positive roots.

Now we claim $x=1$ is a zero of $t(x)$ with multiplicity at least 2. Indeed, by computation, it yields that
\begin{eqnarray*}
t'(x)=\xi'(1)\xi'''(x)x(1-x)-2\xi''(x)[\xi'(1)x-\xi'(x)],
\end{eqnarray*}
and
\begin{eqnarray*}
t''(1)=-\xi'(1)[\xi'''(1)-\frac{2\xi''(1)[\xi''(1)-\xi'(1)]}{\xi'(1)}].
\end{eqnarray*}
Regard $t''(1)$ as a quadratic equation with respect to $\lambda$, which can be expressed explicitly as $(s-p)^2(s^2-3s+2+p^2+3ps-3p)\lambda^2-s(s-p)(2s^2-6s+4-p^2+3ps-3p)\lambda+s^2(s-1)(s-2)=0$
If $s^2+6s-7+p^2+6p-6ps>0,$ it has two roots $\lambda^*_1$ and $\lambda^*_2.$

Thus when $\lambda \neq \lambda^*_1,\lambda^*_2$, $t(1)=0$, $t'(1)=0$ and $t''(1) \neq 0,$ which implies that $x=1$ is a zero of $t(x)$ with multiplicity 2.

Moreover, when $\lambda \in (\lambda^*_1,\lambda^*_2),$ $\xi'''(1)>\frac{2\xi''(1)[\xi''(1)-\xi'(1)]}{\xi'(1)}$ and $t''(1)<0$, which implies that  $x=1$ is a local maximum, Then there exists $\delta_1>0$ such that $t(x)<0$ for $x \in (1-\delta_1,1+\delta_1).$ Moreover, since the coefficient of the term with smallest index in $t(x)$ is strictly positive, there exists $\delta_0>0$ such that $t(x)>0$ for $x \in (0,\delta)$. Then by intermediate value theorem, there exists $x_1 \in(0,1)$ such that $t(x_1)=0.$ Also, since the leading coefficient of $t(x)$ is also strictly positive, $t(x)$ will tend to $+\infty$ as $x \rightarrow \infty.$ By similar reasoning, there exists $x_2 \in (1,+\infty),$ such that $t(x_2)=0.$ Since $t(x)$ can have at most 2 roots besides $x=1,$ $x_1$ and $x_2$ are the exact remaining 2 roots of $t(x).$

 When $\lambda>\lambda^*_2$ or $\lambda<\lambda_1^*$, it holds that $\xi'''(1)<\frac{2\xi''(1)[\xi''(1)-\xi'(1)]}{\xi'(1)}$ and $t''(1)>0$, which implies that $x=1$ is a local minimum of $t(x).$ Since $t(x)$ is strictly positive in some neighborhood near 0, $t(x)$ has either 0 or 2 zeros in $(0,1).$

Now we prove that for $\lambda \in (\lambda^*_2,1)$, $t(x)$ has no roots in $(0,1).$ We prove by contradiction and assume there exists $\lambda \in (\lambda^*_2,1)$ such that $t(x)$ has 2 roots in $ (0,1).$ By continuity, there exists $\lambda' \in  (\lambda^*_2,1)$ and $x' \in (0,1)$ such that $x'$ is a zero of $t(x)$ with multiplicity 2.


Since $x'$ and $\lambda'$ satisfies that
$t(x')=t'(x')=0, $
it yields that for $i=1,2,$
\begin{eqnarray*}
t'(x')&=&\xi'(1)x'(1-x') \Big( \xi'''(x')-\frac{2\xi''(x')[\xi'(1)x'-\xi'(x')]}{\xi'(1)x'(1-x')} \Big)\\
&=&\xi'(1)x'(1-x') \Big( \xi'''(x')-\frac{2\xi''(x')[x'\xi''(x')-\xi'(x')]}{x'\xi'(x')} \Big)=0
\end{eqnarray*}
which implies that $x'\xi'(x') \xi'''(x')-2\xi''(x')[x'\xi''(x')-\xi'(x')] =0.$

Recall that
\begin{eqnarray*}
&&x\xi'(x)\xi'''(x)-2\xi''(x)[x\xi''(x)-\xi'(x)]\\
&=&-p^2(p-1)(p-2)\lambda^2x^{2p-3}+ps\lambda(1-\lambda)x^{p+s-3}\cdot[s^2+3s-4+p^2+3p-4ps]\\
&&-s^2(s-1)(s-2)(1-\lambda)^2x^{2s-3},
\end{eqnarray*}
which is a quadratic equation with respect to $x^{s-p}$ and then $(x')^{s-p}=\frac{p (s^2+3s-4+p^2+3p-4ps) \mp p \sqrt{\Delta}}{2s(s-1)(s-2)} \cdot \frac{\lambda}{1-\lambda},$ where $\Delta= (s^2+3s-4+p^2+3p-4ps)^2-4(p-1)(p-2)(s-1)(s-2)$.

By computation, the following inequalities are equivalent to each other:
\begin{eqnarray*}
x'<1
&\Longleftrightarrow&\frac{p (s^2+3s-4+p^2+3p-4ps) \mp p \sqrt{\Delta}}{2s(s-1)(s-2)} \cdot \frac{\lambda}{1-\lambda} <1\\
&\Longleftrightarrow&\lambda <\frac{2s(s-1)(s-2)}{ \Big( p (s^2+3s-4+p^2+3p-4ps) \mp p \sqrt{\Delta}2s(s-1)(s-2) +2s(s-1)(s-2) \Big)} \\
&\Longleftrightarrow&\lambda< \frac{2s(s-1)(s-2)}{(s-p) \Big( 2s^2-6s+4-p^2+3ps-3p \mp p\sqrt{s^2+6s-7+p^2+6p-6ps}\Big)} \\
&\Longleftrightarrow&\lambda<\lambda^*_{2,1}.
\end{eqnarray*}
The last equivalence can be checked directly by plugging the expression above in the quadratic equation that defines $\lambda^{*}$.
Then it contradicts with the assumption $\lambda > \lambda^*_2,$ which proves the conclusion.

Finally we consider the case when $s^2+6s-7+p^2+6p-6ps<0.$ Under this assumption, we obtain that $t''(1)>0$ for $\lambda \in (0,1)$, which implies that $x=1$ is a local minimum. Then by similar reasoning as above, it yields that $t(x)$ has either 2 or 0 zeros in $(0,1)$ for any $\lambda \in (0,1)$, which ends the proof.\end{proof}

We are now ready to provide the proof of Theorem \ref{coro1rsb}, \ref{coro12rsb}$(i)$ and Theorem \ref{coro12frsb}$(i)$. 

\begin{proof}[Proof of Theorem \ref{coro1rsb}]

We first prove that if $s^2+6s-7+p^2+6p-6ps>0$ and $\Psi(\lambda^*_1)>0,$ then it holds that  the Non-1RSB phase appears for some $\lambda \in [0,1]$.  Indeed, in the proof of Proposition \ref{Cr1rsb}, for $\lambda$ satisfying that $\Psi(\lambda)>0,$ it holds that $\zeta''(1)>0,$ which implies that the model $\xi$ cannot  be 1RSB.

Now we prove that if the Non-1RSB phase appears, then it yields that $\Psi(\lambda^*_1)>0.$ We prove its contrapositive and assume $\Psi(\lambda^*_1)<0.$
We first compute the derivatives of $h^1_1(x,\lambda)$ with respect to $x$ as follows:
\begin{eqnarray*} 
\frac{d}{dx}h^1_1(x,\lambda)
&=&-\frac{\xi'(1)-\xi'(x)-\xi''(x)(1-x)}{[\xi'(1)-\xi'(x)](1-x)}-\frac{2}{x}+\frac{2\xi(x)[\xi'(x)+x\xi''(x)}{x^2\xi'(x)^2}\\
&&+\frac{\xi'(x)[\xi'(1)-\xi'(x)][x\xi'(x)(1-x)-\xi(x)(1-2x)]-\xi'(1)\xi(x)\xi''(x)x(1-x)}{\xi'(1)\xi'(x)^2(1-x)^2x^2}\\
&&-\frac{\xi'(1)\Big[\xi'(x)[\xi'(1)-\xi'(x)][\xi(x)-x\xi'(x)]+x(1-x) \big [ \xi'(x)[x\xi'(x)-2\xi(x)]+\xi'(1)\xi(x)] \big] \Big]}{x^2\xi'(x)^2[\xi'(1)-\xi'(x)]^2} \\
&&=\left [\frac{\xi''(x)(1-x)}{\xi'(1)-\xi'(x)} -\frac{\xi'(x)}{\xi'(1)x} \right ]\cdot \frac{\xi'(1)-\xi'(x)}{1-x} \\
&& \left \{ \frac{\xi(x)(1-2x)}{\xi'(x)^2x(1-x)} +\frac{1}{\xi'(1)-\xi'(x)}-\frac{\xi'(1)(1-x)[x\xi'(x)-\xi(x)]}{x\xi'(x)[\xi'(1)-\xi'(x)]^2}-\frac{\xi'(1)(1-x)\xi(x)}{x\xi'(x)^2[\xi'(1)-\xi'(x)]} \right \} \\
&&=\Big [\frac{\xi''(x)(1-x)}{\xi'(1)-\xi'(x)} -\frac{\xi'(x)}{\xi'(1)x} \Big ]\cdot  \frac{\xi'(x)-\xi'(1)x}{x(1-x)^2\xi'(x)^2[\xi'(1)-\xi'(x)]} \times\\
&&\Big[ \xi'(x)(1-x)[\xi(x)-x\xi'(x)]+x\xi(x)[\xi'(1)-\xi'(x)]\Big].
\end{eqnarray*}

We then claim for $x \in (0,1),$ $\xi'(x)(1-x)[\xi(x)-x\xi'(x)]+x\xi(x)[\xi'(1)-\xi'(x)]$ is always strictly positive. If the claim holds, then the sign of $\frac{d}{dx} h^1_1(x)$ is the opposite as $t(x).$

Indeed by computation, it yields that for $x \in (0,1),$
\begin{eqnarray*}
\xi(x)&\cdot& \frac{\xi'(1)-\xi'(x)}{1-x}-\xi'(x)\cdot \Big [ \xi'(x)-\frac{\xi(x)}{x} \Big] \\
&=&\xi(x)[p\lambda\sum^{p-2}_{k=0}x^k+s(1-\lambda)\sum^{s-2}_{k=0}x^k]-\xi'(x)[(p-1) \lambda x^{p-1}+(s-1)(1-\lambda)x^{s-1}]\\
&=& p\lambda^2 x^p \sum^{p-2}_{k=0}x^k-p(p-1)\lambda^2x^{2p-2}+s(1-\lambda)^2x^s\sum^{s-2}_{k=0}x^k-s(s-1)(1-\lambda^2)^2x^{2s-2}\\
&&+s\lambda(1-\lambda)x^p\sum^{s-2}_{k=0}x^k+p\lambda(1-\lambda)x^s\sum^{p-2}_{k=0}x^k \\
&&-p(s-1)\lambda(1-\lambda)x^{p+s-2}-s(p-1)\lambda(1-\lambda)x^{p+s-2}\\
&\geq& x^p \lambda(1-\lambda) \Big\{ s\sum^{s-2}_{k=0}x^k-s(s-1)x^{s-2}+px^{s-p}\sum^{p-2}_{k=0}x^k-p(p-1)x^{s-2}\\
&&+x^{s-2} \Big [s(s-1)+p(p-1)-p(s-1)-(p-1)s \Big ]x^{s-2} \Big \} \\
&\geq& x^{p+s-2}\lambda(1-\lambda)(s-p)^2>0,
\end{eqnarray*}
which implies that $\xi'(x)(1-x)[\xi(x)-x\xi'(x)]+x\xi(x)[\xi'(1)-\xi'(x)]>0.$

By Lemma \ref{tx}, for $\lambda \in (0,\lambda^*_1),$ $t(x)$ has either 0 or 2 roots in $ (0,1).$ If $t(x)$ has no roots in $ (0,1),$ combining with the fact that $\lim_{x \rightarrow 0+} \frac{d^2}{dx^2} h^1_1(x) =-\frac{1}{2}<0$, we obtain that $h^1_1(x)<0$ for $x \in (0,1),$ which implies that the model $\xi$ is then 1RSB by Proposition \ref{Cr1rsb}.

We then compute the derivative of $h^2_1(x,\lambda)$ to obtain
\begin{align*}
\frac{d}{dx} h_1^2(x,\lambda)&=\frac{[\xi'(x)-x\xi''(x)] }{x^3\xi'(x)^3\xi''(x)^2} \cdot [x\xi'(x)^2-x\xi(x)\xi''(x)-\xi(x)\xi'(x)]
\\&\quad \times \left [x\xi'(x)\xi'''(x)-2\xi''(x)[x\xi''(x)-\xi'(x)] \right].
\end{align*}
Notice that $x\xi'(x)^2-x\xi(x)\xi''(x)-\xi(x)\xi'(x)=-(s-p)^2\lambda(1-\lambda)x^{p+s-1}<0$ and $\xi'(x)-x\xi''(x)<0$. Thus the sign of $\frac{d}{dx} h_1^2(x,\lambda)$ is the same as $t^2_1(x):=x\xi'(x)\xi'''(x)-2\xi''(x)[x\xi''(x)-\xi'(x)].$ 

By algebraic calculation, we obtain that
\begin{eqnarray*}
x\xi'(x)\xi'''(x)&-&2\xi''(x)[x\xi''(x)-\xi'(x)]\\
&=&-p^2(p-1)(p-2)\lambda^2x^{2p-3}+ps\lambda(1-\lambda)x^{p+s-3}\cdot[s^2+3s-4+p^2+3p-4ps]\\
&&-s^2(s-1)(s-2)(1-\lambda)^2x^{2s-3}.
\end{eqnarray*}
Since the right-hand side is a quadratic equation with respect to $x^{s-p},$ it holds that $t^2_1(x)$ has either 0 or 2 roots.
If $t^2_1(x)$ has no zeros, by the fact that $h^2_1(x,\lambda)$ is decreasing near $0$, we obtain that $h^2_1(x,\lambda)<0$ for $x \in (0,1),$ which implies that the model $\xi$ is 1RSB by Proposition \ref{Cr1rsb}.

Thus we consider the case when $\lambda \in (0,\lambda^*_1)$ and both $t(x),t^1_2(x)$ have two zeros in $(0,1).$ We still denote the two zeros of $t(x)$ by $\bar{q}_1$ and $\bar{q}_2.$ Also, we denote the zeros of $t^1_2(x)$ by $\tilde{q}_1<\tilde{q}_2.$

Since $\tilde{h}^2_1(\bar{q}_2)=h^1_1(\bar{q}_2),$ we conclude that the maximum of $\tilde{h}^2_1(x)$ is always larger than the maximum of $h^1_1(x)$ in $[0,1].$

Also, we consider $h^2_1(x,\lambda)$ as a function also with respect to $\lambda$ and compute 
\begin{eqnarray*}
\frac{\partial }{\partial \lambda} h^2_1(x,\lambda) &=&\frac{(p-s)(p-1)(s-1)x^{p+s-4}}{\xi''(x)^2}-\frac{2(s-p)}{\xi'(x)^2}x^{p+s-2} + \frac{ps(s-p)x^{p+s-3}}{\xi'(x)\xi''(x)} \\
&&+\frac{1}{\xi'(x)^3} \Big( p\lambda \big[ -p(p-1)+2(p-1)s-s(s-1) \Big]x^{2p+s-3}\\
&&+s(1-\lambda) \Big[ p(p-1)+s(s-1)-2p(s-1) \big] x^{p+2s-3} \Big)\\
&=& \frac{(p-s)x^{p+s-4}}{\xi'(x)^3\xi''(x)^2} \cdot \Big \{ (p-1)(s-1)\xi'(x)^3 +2x^2\xi'(x)\xi''(x)^2-ps\xi'(x)^2\xi''(x)x \\
&&+p\lambda (1+s-p)x^{p+1} \xi''(x)^2 +s(1-\lambda)(1+p-s)x^{s+1}\xi''(x)^2 \Big \} \\
&=&  \frac{(p-s)x^{p+s-4}[x\xi''(x)-\xi'(x)]}{\xi'(x)^3\xi''(x)^2} \cdot \Big \{ x\xi'(x)\xi'''(x)-2\xi''(x)[x\xi''(x)-\xi'(x)] \Big \}
\end{eqnarray*}

Then it yields that $\frac{\partial }{\partial \lambda} h^2_1(x,\lambda)|_{x= \tilde{q}_1}=\frac{\partial }{\partial \lambda} h^2_1(x,\lambda)|_{x= \tilde{q}_2}=0$ which implies that for fixed $s-1>p>2$, the maximum and minimum of $h^1_2(x,\lambda)$ are constants, regardless of the values of $\lambda \in (0,1)$. 

Thus if $\lim_{x \rightarrow \infty} h^2_1(x,\lambda)=\Psi(\lambda)<0$ for $\lambda \in (0,1),$ it yields that $h^2_1(x,\lambda)<0$ for $x \in (0,1).$ Thus by Proposition \ref{Cr1rsb}, the model  is 1RSB for all values of $\lambda \in [0,1]$.

Now it suffices for us to show that $\Psi(\lambda)$ achieves its maximum at $\lambda=\lambda^*_1.$
Indeed by computing the critical points of $\Psi(\lambda)$, it yields that  
\begin{eqnarray*}
0&=&\frac{d}{d \lambda}\Big (\Psi(\lambda) \Big)\\
&=&(p+s)\xi''(1) \xi'(1)^3-(p+s-1)\xi'(1)^3 [\xi'(1)-1]\\
&&-\xi''(1)^2\xi'(1)^2-(p+s+1) \xi''(1)^2\xi'(1)+2\xi''(1)^3\\
&=&\big( (s-p)(s+p-2)\lambda-(s-2)s \big) \Big [\xi'(1)\xi'''(1)-2\xi''(1)[\xi''(1)-\xi'(1)]\Big],
\end{eqnarray*}
which implies that the equation has altogether 3 zeros, which are $\lambda^*_1,\lambda^*_2 \text{ and } \frac{s(s-2)}{(s-p)(s+p-2)}.$ Here $\lambda^*_1$ and $\frac{s(s-2)}{(s-p)(s+p-2)}$ are local maximum of $\Psi(\lambda)$ and $\lambda^*_2$ is a local minimum of $\Psi(\lambda),$ which proves the claim.

\end{proof}

Now we turn to the proof of Theorem \ref{coro12rsb}$(i)$ and Theorem \ref{coro12frsb}$(i)$ which are similar.

\begin{proof}[Proof of Theorem \ref{coro12rsb}$(i)$ and Theorem \ref{coro12frsb}$(i)$]

We first prove that the model $\xi(x)=\lambda x^p +(1-\lambda)x^s$ is 1RSB if $\lambda \in [0,\lambda_{1 \rightarrow 2}]\cup [\lambda_{2 \rightarrow 1},1].$

We first consider the case when $\lambda \in [\lambda_{2 \rightarrow 1},1]$ and split it into 2 intervals as $ [\lambda_{2 \rightarrow 1}, \lambda^*_2]$ and $[\lambda^*_2,1].$

By Lemma \ref{tx}, for $\lambda \in [\lambda_{2 \rightarrow 1},\lambda^*_2],$ $t(x)$ has exactly 1 root in $(0,1).$ 
Also, by L'H\^opital's rule, we notice that 
\begin{eqnarray*}
\lim_{x \rightarrow 1}h^1_1(x,\lambda)&=&-\frac{\xi'(1)-\xi''(1)}{\xi''(1)} + \frac{\big(\xi'(1)-\xi''(1)\big)^2}{\xi'(1)^2\xi''(1)}-\log \Big( \frac{\xi''(1)}{\xi'(1)} \Big) \\
&=&\Psi(\lambda)<0.
\end{eqnarray*}
Combining these two facts, it yields that $h^1_1(x,\lambda)<0$ for $x \in (0,1).$ Then by Proposition \ref{Cr1rsb}, the model $\xi$ is 1RSB.

Now we consider the case that $\lambda \in (\lambda^*_2,1).$  Since $t(x)$ has no zeros in $(0,1)$ when $\lambda \in (\lambda^*_2,1),$ then by the fact that $h^2_1(0,\lambda)=-1<0$ and $h^2_1(0,\lambda)=0$, it holds that $h^2_1(x,\lambda)<0$ for $x \in (0,1).$ Then we conclude that for $\lambda \in (\lambda_{2 \rightarrow 1},1),$ the model is 1RSB by Proposition \ref{Cr1rsb}.

Now we turn to the case when $\lambda$ is near 0. To be more specific, we restrict the range of $\lambda$ to be $[0, \lambda_{2 \rightarrow 1}].$

Since $\psi(x)=zs(1-\lambda)x^{s-1}-(1+z)s(1-\lambda)x^{s-2}+zp\lambda x^{p-1}-(1+z)p\lambda x^{p-2}+\xi'(1),$ by Descartes' rule, $\psi(x)$ has at most 3 zeros except $x=1$. Recall that when the model $\xi(x)$ is 1RSB, it always holds that 
\begin{eqnarray*}
\zeta''(1)=\frac{\xi''(1)}{(1+z)} \Big [ (1+z)-\frac{\xi''(1)}{\xi'(1)} \Big]>0, 
\end{eqnarray*}
which implies that $\zeta(x)$ has either 1 or 3 critical points in $(0,1)$. In particular, when $\zeta(x)$ has 1 critical point in $(0,1),$ the model $\xi(x)$ must be 1RSB. When $\zeta(x)$ has 3 critical points in $(0,1),$ then two of them are local minimums and the rest one is a local maximum.

Then when $\lambda< \lambda_{2 \rightarrow 1},$ the boundary of 1RSB and Non-1RSB phase occurs if and only if $\zeta(x^*)=0,$ where $x=x^*$ is the local maximum of $\zeta(x).$ We denote the corresponding model parameter by $\lambda^*.$ By the proof of Proposition \ref{Cr1rsb}, we obtain that $h^1_1(x^*,\lambda^{*})=\tilde{h}^1_2(x^*)=0.$ Here we notice that $(\lambda^*,x^*)$ is a solution of the system of equations \eqref{eqnh1}. 

Now we prove that $(\lambda^*,x^*)$ is the unique solution of the system \eqref{eqnh1}. Since when $\lambda <  \lambda_{2 \rightarrow 1},$ the solution of the system \eqref{eqnh1} is exactly the boundary of 1RSB and Non-1RSB phase, it is equivalent for us to show that the boundary of 1RSB and Non-1RSB phase appears only once.

Now it suffices for us to show that $\frac{d}{d\lambda} \zeta(x^*)|_{\lambda=\lambda^*}>0.$ By computation, we obtain that
\begin{eqnarray}\label{dzeta}
\frac{d}{d\lambda} \zeta(x^*)|_{\lambda=\lambda^*}&=&(x^*)^p-(x^*)^s+\big[p(x^*)^{p-1}-s(x^*)^{s-1}\big](1-x^*)+\frac{p(x^*)^{p-1}-s(x^*)^{s-1}}{z} \nonumber \\
&+&\frac{(s-p)(1+z)}{z^2} \log \Big( 1+\frac{z\xi'(x)}{\xi'(1)} \Big) -\frac{ps \big[ (x^*)^{p-1}-(x^*)^{s-1} \big] (1+z)}{z \big[ \xi'(1)+z\xi'(x^*)\big]} \nonumber \\
&+&z'_\lambda \cdot \Big \{ -\frac{\xi'(x^*)}{z^2}+\frac{\xi'(1)(z+2)}{z^3} \log \Big( 1+\frac{z\xi'(x^*)}{\xi'(1)} \Big)-\frac{\xi'(1)\xi'(x^*)(1+z)}{z^2\Big[\xi'(1)+z\xi'(x^*)\Big]}\Big\}.
\end{eqnarray}

Since $x^*$ and $\lambda^*$ satisfy that $\zeta(x^*)=\zeta'(x^*)=0,$ we obtain that 
\begin{eqnarray}\label{x^*1}
\log \Big (1+\frac{z\xi'(x^*)}{\xi'(1)} \Big)= \frac{z^2}{\xi'(1)(1+z)} \Big[ \xi(x^*)+\xi'(x^*)(1-x^*)+\frac{\xi'(x^*)}{z} \Big]
\end{eqnarray}
and
\begin{eqnarray}\label{x^*2}
z=\frac{\xi'(1)x^*-\xi'(x^*)}{\xi'(x^*)(1-x^*)}.
\end{eqnarray}

Moreover, since for $\lambda=\lambda^*$, $\zeta(x)$ has 3 critical points in $(0,1),$ then it holds that $\zeta''(1)>0$ and then
$1+z>\frac{\xi''(1)}{\xi'(1)}.$

Now we plug relations \eqref{x^*1},\eqref{x^*2}, \eqref{x^*3} and \eqref{x^*4}  in \eqref{dzeta} and obtain that
\begin{eqnarray*}
\frac{d}{d\lambda} \zeta(x^*)|_{\lambda=\lambda^*}&=&\frac{s(x^*)^p-p(x^*)^s}{\xi'(1)}-\frac{(s-p)(z+2)\xi(x^*)}{\xi'(1)[z+2-\xi'(1)]}+\frac{(s-p)x^*\xi'(x^*)}{\xi'(1)[z+2-\xi'(1)]} \\
&=&\frac{1}{\xi'(1)[z+2-\xi'(1)]} \Big\{ (z+2) \big[ s(x^*)^p-p(x^*)^s \big] -\xi'(1) \big[ s(x^*)^p-p(x^*)^s \big] \\
&&-(z+2)(s-p)\xi(x^*)+(s-p)x^*\xi'(x^*) \Big \} \\
&=&\frac{1}{\xi'(1)[z+2-\xi'(1)]} \Big\{ (z+2)\xi'(1) \big[ (x^*)^p-(x^*)^s \big]-[p^2 \lambda+s^2(1-\lambda)]\big[ (x^*)^p-(x^*)^s\big] \Big\} \\
&=&\frac{(x^*)^p-(x^*)^s}{[z+2-\xi'(1)]} \cdot \Big( 1+z-\frac{\xi''(1)}{\xi'(1)} \Big)>0,
\end{eqnarray*}
which implies that the uniqueness of the boundary of 1RSB and Non-1RSB phase when $\lambda \in [0,\lambda_{2 \rightarrow 1})$. Thus when $\lambda \in [0,\lambda_{2 \rightarrow 1})$, the model $\xi(x)$ is 1RSB.

Now we consider the case that $\lambda \in (\lambda_{1 \rightarrow 2},  \lambda_{2 \rightarrow 1}).$
We denote the zero of $\Psi(\lambda)$ in $(0,\lambda^*_1)$ by $\lambda'_{1 \rightarrow 2RSB}.$ When $\lambda \in (\lambda_{1 \rightarrow 2}, \lambda'_{1 \rightarrow 2RSB}),$ it holds that $\zeta(x^*)>0,$ which implies that $\xi(x)$ is not 1RSB. When $\lambda \in (\lambda'_{1 \rightarrow 2RSB}, \lambda_{2 \rightarrow 1}),$ it holds that $\zeta''(1)>0,$ and therefore the model $\xi(x)$ is not 1RSB.

\end{proof}

\subsection{The 2RSB Criterion}\label{sec2rsb}

Assume that the Parisi measure $\nu_P$ has the form 
\begin{eqnarray*}
\nu(ds)=k_1 \cdot \mathbbm{1}_{[0,q)}(s)ds+k_2 \cdot \mathbbm{1}_{[q,1)}(s)ds+\Delta \delta_{\{1\}}(ds),
\end{eqnarray*}

By direct computation, we obtain that for $0 \leq x \leq q$,
\begin{eqnarray*}
g(x)&&=\xi(1)-\xi(x)-\frac{1}{k_1^2} \log(1+\frac{k_1(q-x)}{k_2(1-q)+\Delta})+\frac{q-x}{k_1[k_1x+k_2(1-x)+\Delta]} \\
&&-\frac{1}{k_2^2}\log(1+\frac{k_2(1-x)}{\Delta})+\frac{1-q}{k_2(k_2(1-q)+\Delta)}-\frac{(1-q)q}{[k_2(1-q)+\Delta][k_1q+k_2(1-q)+\Delta]}
\end{eqnarray*}
and for $q \leq x \leq 1,$
\begin{eqnarray*}
g(x)&=&\xi(1)-\xi(x)-\frac{q(1-x)}{[k_2(1-q)+\Delta][k_1q+k_2(1-q)+\Delta]}\\
&&+\frac{1-x}{k_2[k_2(1-q)+\Delta]}-\frac{1}{k^2_2}\log(1+\frac{k_2(1-x)}{\Delta})
\end{eqnarray*}

Now we set $z_1=\frac{k_1q}{\Delta}$ and $z_2=\frac{k_2(1-q)}{\Delta}.$ The following was shown in \cite{AZeng}.

\begin{theorem}[Theorem 10 in \cite{AZeng}]\label{2RSBCr}
Let $s-1>p>2,s,p\in \mathbb{N}$, the model $\xi(x)=\lambda x^p+(1-\lambda)x^s$ is 2RSB provided the following conditions hold:
\begin{enumerate}
\item The following equation holds
\begin{eqnarray}\label{eq4}
1+z_1+z_2=\frac{q[\xi'(1)-\xi'(q)]}{\xi'(q)(1-q)},
\end{eqnarray}
\item $z_1>0,z_2>\frac{(1-q)z_1}{q}$ and $q\in (0,1)$ satisfy the following two equations
\begin{eqnarray}\label{f1}
f_1(q,z_2)&:=&-\frac{[q\xi'(q)-\xi(q)](1-q)(1+z_2)}{\xi'(1)-\xi'(q)}-q^2\log \frac{q[\xi'(1)-\xi'(q)]}{(1+z_2)\xi'(q)(1-q)}+q^2  \nonumber \\
&&-\frac{2\xi(q)q}{\xi'(q)}+\frac{\xi(q)q^2[\xi'(1)-\xi'(q)]}{(1+z_2)\xi'(q)^2(1-q)}=0
\end{eqnarray}
and 
\begin{eqnarray}\label{f2}
f_2(q,z_2):= (1-q)[\xi'(1)-\xi'(q)](\frac{1+z_2}{z^2_2}\log(1+z_2)-\frac{1}{z_2})+\xi'(q)(1-q)-1+\xi(q)=0.
\end{eqnarray}

\item the following inequality holds
\begin{eqnarray}\label{eq15}
\xi''(q)(1+z_2)(1-q) \leq \xi'(1)-\xi'(q),
\end{eqnarray}

\item $h_1(1)>0$ and $h_2(0)<0$, where 
\begin{eqnarray*}
h_1(x)=\xi'(x)(q+qz_1+qz_2-z_1x)-(1+z_2)\xi'(q)x,
\end{eqnarray*}
and
\begin{eqnarray*}
h_2(x)=[\xi'(x)-\xi'(q)](1+z_2-q-z_2x)-(x-q)[\xi'(1)-\xi'(q)].
\end{eqnarray*}
\end{enumerate}
\end{theorem}

\begin{remark}
Based on the definition of $f_1(q,z_2)$ and $f_2(q,z_2)$, we notice that
\begin{eqnarray*}
h^1_1(x,\lambda)= x^{-2} \cdot f_1\left(x, \frac{\xi'(1)x}{\xi'(x)}-1\right) \text{ , } h^1_2(x,\lambda)=  f_2\left(x,\frac{\xi'(1)x}{\xi'(x)}-1\right),
\end{eqnarray*}
and
\begin{eqnarray*}
h_1^2(x,\lambda)=x^{-2}\cdot f_1\left(x, \frac{\xi'(1)-\xi'(x)}{\xi''(x)\cdot(1-x)}-1\right) \text{ , } h^2_2(x,\lambda)=    f_2\left(x, \frac{\xi'(1)-\xi'(x)}{\xi''(x)\cdot(1-x)}-1\right).
\end{eqnarray*}

\end{remark}

Now we claim that condition $(4)$ can be deduced by the other three conditions in Theorem \ref{2RSBCr}. Thus we simplify Theorem \ref{2RSBCr} and derive an exact condition for the model $\xi(x)$ to be 2RSB as follows:
\begin{theorem}\label{2RSB}
For some $s-1>p>2,s,p\in \mathbb{N}$, the model $\xi(x)=\lambda x^p+(1-\lambda)x^s$ if and only if the following conditions hold:
\begin{eqnarray}\label{Condition2rsb}
&(i)& \text{ there exists }  z_2>0 \text{ and }q\in (0,1) \text{ satisfying} \nonumber \\
&&\frac{\xi'(1)q}{\xi'(q)}<1+z_2<\frac{\xi'(1)-\xi'(q)}{\xi''(q)(1-q)},\\
&(ii)& q \text{ and } z_2 \text{ satisfy that } f_1(q,z_2)=0 \text{ and } f_2(q,z_2)=0. 
\end{eqnarray}

\end{theorem}

\begin{proof}[Proof of Theorem \ref{2RSB}]
It suffices for us to show that condition $(iv)$ in Theorem \ref{2RSBCr} can be deduced by condition $(i), (ii)$ and $(iii)$ in Theorem \ref{2RSBCr}.
 
 By computation, we obtain that 
 \begin{eqnarray*}
 h_1(1)&=&\xi'(1) \cdot [q+qz_1+qz_2-z_1]-(1+z_2)\xi'(q) \\
 &=&q(1+z_2)\xi'(1)-(1-q)\xi'(1)z_1-(1+z_2)\xi'(q)\\
 &=& [\xi'(1)-\xi'(q)]\cdot[1+z_2-\frac{\xi'(1)q}{\xi'(q)}]
 \end{eqnarray*}
 and 
 \begin{eqnarray*}
 h_2(0)&=&-\xi'(q)\cdot(1+z_2-q)+q[\xi'(1)-\xi'(q)]\\
 &=&-\xi'(q)\cdot (1+z_2)+q\xi'(1)\\
 &=&-\xi'(q)\cdot[1+z_2-\frac{\xi'(1)q}{\xi'(q)}].
 \end{eqnarray*}
  Therefore the condition $h_1(1)>0$ and $h_2(0)<0$ hold if and only if $1+z_2-\frac{\xi'(1)q}{\xi'(q)}>0$.
  
  Since $z_2>\frac{(1-q)z_1}{q}$, we know that $z_1 < \frac{q}{1-q}z_2$. Combining with \eqref{eq4}, it yields that
\begin{eqnarray*}
\frac{q[\xi'(1)-\xi'(q)]}{\xi'(q)(1-q)}&=&1+z_1+z_2 \\
&<&1+\frac{q}{1-q}z_2+z_2 = 1+\frac{1}{1-q}z_2,
\end{eqnarray*}
which implies that $z_2 >\frac{\xi'(1)q}{\xi'(q)}-1.$
Therefore the condition $h_1(1)>0$ and $h_2(0)<0$ are derived.

Moreover if $z_2 >\frac{\xi'(1)q}{\xi'(q)}-1,$ then 
\begin{eqnarray*}
\frac{q[\xi'(1)-\xi'(q)]}{\xi('q)(1-q)}-(1+z_2)&<& \frac{q[\xi'(1)-\xi'(q)]}{\xi'(q)(1-q)} -\frac{\xi'(1)q}{\xi'(q)}\\
&=&\frac{q}{1-q}\cdot \frac{q\xi'(1)-\xi'(q)}{\xi'(q)}\\
&<&\frac{q}{1-q} z_2.
\end{eqnarray*}

Also, notice that since $q \xi''(q)>\xi'(q)$, we obtain that
\begin{eqnarray*}
\frac{q[\xi'(1)-\xi'(q)]}{\xi'(q)(1-q)} > \frac{\xi'(1)-\xi'(q)}{\xi''(q)(1-q)}
\end{eqnarray*}
which implies that $1+z_2<\frac{q[\xi'(1)-\xi'(q)]}{\xi'(q)(1-q)} .$ Then we set $z_1:=\frac{q[\xi'(1)-\xi'(q)]}{\xi'(q)(1-q)} -(1+z_2)$ and $z_1$ satisfies the conditions in Theorem \ref{2RSBCr}.
\end{proof}

The next lemma provides some qualitatively information about the functions $h^1_1, h^1_2, h^2_1$ and $h^2_2$.
\begin{lemma}\label{h4}
The following holds.
\begin{enumerate}
\item For all $x, \lambda$ \[
\text{sgn} \left( \partial_{x} h^1_1(x,\lambda)\right) = - \text{sgn} \left( \partial_{x} h^1_2(x,\lambda)\right).\]

\item If $\partial_{x} h^1_2(x,\lambda) <0$ then $h^1_2(x,\lambda)> h^2_2(x,\lambda)$and $h_1^2(x,\lambda) < h^1_1(x,\lambda)$.

\item Exactly one of the following cases can hold: (i) there exists $0\leq \bar{q}_1<\bar{q}_2\leq 1$ such that the strictly decreasing interval of $h^1_2(x,\lambda)$ is $ [\bar{q}_1,\bar{q}_2];$ (ii) $h^1_2(x,\lambda)$ is strictly increasing for $x\in (0,1).$  

\item $h_1^2(x,\lambda)$ has at most 2 roots in $(0,1)$ and $\lim_{x \rightarrow 1}h^2_1(x,\lambda)=\Psi(\lambda).$

\item  $\text{sgn}\left(\partial_{x} h^2_2(x,\lambda)\right) = \text{sgn}\left(m(x)\right)$. Moreover, $\frac{d^2}{dx^2}h^2_2(1,\lambda)>0$ if and only if $\lambda \in (\lambda_{2 \rightarrow 1F},\lambda'_{2 \rightarrow 1F}).$

\end{enumerate}
\end{lemma}

Before we prove Lemma \ref{h4}, we recall  that $\lambda^*_{1,2}$ and $\lambda_{2 \rightarrow 1F}$,$\lambda'_{2 \rightarrow 1F}$ are solutions of  quadratic equations. We can express $\lambda^*_{1,2}$ explicitly as follows:
\begin{eqnarray*}
\lambda^*_{1,2}=\frac{2s(s-1)(s-2)}{(s-p)\big(a\pm p\sqrt{\Delta}\big)},
\end{eqnarray*}
where $a=2s^2+3(p-2)s-(p-1)(p+4)$.

We also express $\lambda_{2 \rightarrow 1F}$,$\lambda'_{2 \rightarrow 1F}$ as
\begin{eqnarray*}
\lambda_{2 \rightarrow 1F},\lambda'_{2 \rightarrow 1F} =\frac{s^2(s-1)(s-2)}{(s-p)\big( a_F\pm p(p-1)\sqrt{\Delta+2(p-2)(s-2)} \big)},  
\end{eqnarray*}
where $a_F=s^3+(p-3)s^2+2[p(p-2)+1]-p(p-1)(p+1).$

Here we notice that for $s>p>2,$ if $\Delta>0,$ then $\lambda_{2 \rightarrow 1F}$ and $\lambda'_{2 \rightarrow 1F}$ also exist. Also, for $s>p>2,$ it holds that $a>0$ and $a_F>0.$
\begin{lemma}\label{lemmalambda}
If $\Delta>0,$ then $0<\lambda^*_1<\lambda^*_2<1$ and $\lambda^*_1<\lambda_{2 \rightarrow 1F}.$
\end{lemma}

\begin{proof}[Proof of Lemma \ref{lemmalambda}]
Since $s^2+6s-7+p^2+6p-6ps>0,$ by elementary computation, it holds that
\begin{eqnarray*}
0<s^2-6(p-1)s+(p-1)(p-7) <[s-3(p-1)]^2-8(p-1)(p-2)
\end{eqnarray*}
which implies that $s>\sqrt{8(p-1)(p-2)}+3(p-1)>5(p-2).$

Therefore
\begin{eqnarray*}
2s^2-6s+4-p^2+3ps-3p&=&2s^2+3(p-2)s-(p-1)(p-4)\\
&=&2[s+\frac{3}{4}(p-2)]^2-\Big( \frac{9}{8}(p-2)^2+(p-1)(p+4) \Big) \\
&\geq& 65(p-2)^2-(p-1)(p+4)>0,
\end{eqnarray*}
and it also holds that $s^2-3s+2+p^2+3ps-3p=s^2+3(p-1)s+(p-1)(p-2)>0$.
Thus plugging in back in we obtain that  $\lambda^*_1>0.$

Now it suffices for us to show that $\lambda^*_2<1.$
Indeed, we notice that the following inequalities are equivalent:
\begin{eqnarray*}
\lambda^*_2<1
&\Longleftrightarrow&ps  \sqrt{s^2+6s-7+p^2+6p-6ps}< p(-3ps+s^2+3s-2p^2+6p-4) \\
&\Longleftrightarrow&s  \sqrt{s^2-6(p-1)s+(p-1)(p-7)}<s^2-3(p-1)s-2(p-1)(p-2)\\
&\Longleftrightarrow&2(2p+3)(p-1)s^2+12(p-1)^2(p-2)s+4(p-1)^2(p-2)^2>0.
\end{eqnarray*}
Since the last inequality holds for any $p,s>0,$ we obtain that $\lambda^*_2<1.$

Now we turn to the proof of $\lambda^*_1<\lambda_{2 \rightarrow 1F}.$ Since for $s>p>2$, it holds that $a,a_F>0$. Then if $\Delta>0,$ the following inequalities are equivalent:
\begin{eqnarray*} 
\lambda^*_1<\lambda_{2 \rightarrow 1F}
&\Longleftrightarrow&\frac{s}{a_F+p(p-1)\sqrt{\Delta+2(p-2)(s-2)}}>\frac{2}{a+p\sqrt{\Delta}}, \\
&\Longleftrightarrow&s[a+p\sqrt{\Delta}]-2[a_F+p(p-1)\sqrt{\Delta+2(p-2)(s-2)}]>0.
\end{eqnarray*}

By the inequality $\sqrt{x+y}<\sqrt{x}+\sqrt{y}$ for $x,y>0$, we set $x=\Delta,y=2(p-2)(s-2)$ and then obtain that $\sqrt{\Delta+2(p-2)(s-2)}<\sqrt{\Delta}+\sqrt{2(p-2)(s-2)}$. In order to prove $\lambda^*_1<\lambda_{2 \rightarrow 1F},$ it suffices for us to show that 
\begin{eqnarray}\label{ps1}
s\cdot a-2a_F-2p(p-1)\sqrt{2(p-2)(s-2)}>-p \big[s-2(p-1)\big] \sqrt{\Delta}.
\end{eqnarray} 
Note that the right-hand side of the inequality \eqref{ps1} is always negative for $s>p>2$ and $\Delta>0.$ We then claim that the left-hand side of \eqref{ps1} for $s>p>2$, which then implies that $\lambda^*_1<\lambda_{2 \rightarrow 1F}.$

Indeed, by algebraic computation, we obtain that
\begin{eqnarray}\label{ps2}
s \cdot a-2a_F=p(s^2-5ps+5s+2p^2-2),
\end{eqnarray}
which implies that for $s>0>2$ and $\Delta>0,$
\begin{eqnarray*}
&&s\cdot a-2a_F-2p(p-1)\sqrt{2(p-2)(s-2)}\\
&=&p \big[(s^2-5ps+5s+2p^2-2)-2(p-1)\sqrt{2(p-2)(s-2)} \big]\\
&=&\frac{p[s-2(p-1)]^2\Delta}{(s^2-5ps+5s+2p^2-2)+2(p-1)\sqrt{2(p-2)(s-2)}}>0.
\end{eqnarray*}

\end{proof}

Now we come to the proof of Lemma $\ref{h4}.$
\begin{proof}[Proof of Lemma \ref{h4}]

We first compute the derivative of $h^1_2(x,\lambda)$ with respect to $x$ to obtain
\begin{eqnarray*}
\frac{d}{dx}h^1_2(x,\lambda)
&=& \frac{\xi'(1)\xi''(x)x(1-x)-\xi'(x)[\xi'(1)-\xi'(x)]}{[\xi'(1)x-\xi'(x)]^3}  \\
&\times&\qquad \Big \{ \xi'(1)\Big[ x[\xi'(1)-\xi'(x)]+\xi'(x)(1-x) \Big] \log \Big( \frac{\xi'(1)x}{\xi'(x)}\Big) \\
&&\qquad \qquad \qquad - [\xi'(1)-\xi'(x)+\xi'(1)(1-x)][\xi'(1)x-\xi'(x)] \Big \}.
\end{eqnarray*}

Note that by \eqref{logz}, we set $z=\frac{\xi'(1)x}{\xi'(x)}-1$ and obtain that
\begin{eqnarray*}
\xi'(1)&\Big[& x[\xi'(1)-\xi'(x)]+\xi'(x)(1-x) \Big] \log \Big( \frac{\xi'(1)x}{\xi'(x)}\Big) \\
&& \qquad \qquad - [\xi'(1)-\xi'(x)+\xi'(1)(1-x)][\xi'(1)x-\xi'(x)] \\
&\geq&\frac{2\xi'(1)[\xi'(1)x-\xi'(x)]}{\xi'(1)x+\xi'(x)} \Big[ x[\xi'(1)-\xi'(x)]+\xi'(x)(1-x) \Big]  \\
&&-[\xi'(1)-\xi'(x)+\xi'(1)(1-x)][\xi'(1)x-\xi'(x)] \\
&=& \frac{[\xi'(1)x-\xi'(x)]^3}{\xi'(1)x+\xi'(x)}>0.
\end{eqnarray*}

Thus $\frac{d}{dx}h^1_2(x,\lambda)$ has the same sign as $t(x)$. Recall that $\frac{d}{dx}h^1_2(x,\lambda)$ has the opposite sign to $t(x)$, which implies the first part of Lemma \ref{h4}.

Now we turn to the proof of item (2). We fix $x \in (0,1)$ and regard $f_1(x,z_2)$ as a function with respect to $w:=z_2+1$, denoted by $\phi(w).$ Then taking derivative in $w$, it yields that
 \begin{eqnarray*}
 \phi'(w)=-\frac{[x\xi'(x)-\xi(x)](1-x)}{\xi'(1)-\xi'(x)}+\frac{x^2}{w}-\frac{\xi(x)x^2[\xi'(1)-\xi'(x)]}{\xi'(x)^2(1-x)w^2}
 \end{eqnarray*}
 Solving the quadratic equation given by $\phi'(w)=0,$ we obtain the two critical points of $\phi(w)$ as follows:
 \begin{eqnarray*}
 w_{1}=\frac{x\xi(x)[\xi'(1)-\xi'(x)]}{\xi'(x)[x\xi'(x)-\xi(x)](1-x)} \text{ and } w_{2}=\frac{x[\xi'(1)-\xi'(x)]}{\xi'(x)(1-x)}
 \end{eqnarray*}
Since $\phi(w)$ approaches $+\infty$ as $w \rightarrow 0,$ we obtain that $\phi(w)$ is strictly decreasing on $(0,w_1)$. Since $z=w-1,$ it is equivalent to say for fixed $x \in (0,1),$ $f_1(x,z_2)$ is strictly decreasing on the interval $\big(-1,\frac{x\xi(x)[\xi'(1)-\xi'(x)]}{\xi'(x)[x\xi'(x)-\xi(x)](1-x)}-1\big).$

Then when $h^1_2(x,\lambda)$ is strictly decreasing, we have $\xi'(1)\xi''(x)x(1-x)-\xi'(x)[\xi'(1)-\xi'(x)]<0,$ which implies that $\frac{\xi'(1)x}{\xi'(x)}<\frac{\xi'(1)-\xi'(x)}{\xi''(x)(1-x)}.$ Moreover we notice that for $x \in (0,1),$

\begin{eqnarray*}
w_1-\frac{\xi'(1)-\xi'(x)}{\xi''(x)(1-x)}
&=&\frac{\xi'(1)-\xi'(x)}{\xi'(x)\xi''(x)[x\xi'(x)-\xi(x)](1-x)}\cdot \Big \{ x\xi(x)\xi''(x)-\xi'(x)[x\xi'(x)-\xi(x)]\Big \} \\
&=&\frac{\xi'(1)-\xi'(x)}{\xi'(x)\xi''(x)[x\xi'(x)-\xi(x)](1-x)}\cdot \Big \{ [\lambda x^{p+1}+(1-\lambda)x^{s+1}]\xi''(x)\\
&&-[p\lambda x^{p-1}+s(1-\lambda)x^{s-1}][(p-1)\lambda x^p+(s-1)(1-\lambda)x^s] \Big \}\\
&=&\frac{\xi'(1)-\xi'(x)}{\xi'(x)\xi''(x)[x\xi'(x)-\xi(x)](1-x)} \cdot  \lambda(1-\lambda)(s-p)^2x^{p+s-1}>0,
\end{eqnarray*}
which implies that when $h^1_2(x,\lambda)$ is strictly decreasing, it holds that $\frac{\xi'(1)x}{\xi'(x)}<\frac{\xi'(1)-\xi'(x)}{\xi''(x)(1-x)}<w_1.$ Since $f_1(x,z_2)$ is strictly decreasing on $(-1,w_1-1),$ we obtain that 
\[
h^1_1(x)= x^{-2} \cdot f_1(x, \frac{\xi'(1)x}{\xi'(x)}-1)>x^{-2}\cdot f_1(x, \frac{\xi'(1)-\xi'(x)}{\xi''(x)\cdot(1-x)}-1) =h_1^2(x,\lambda)
\]
Furthermore, since for fixed $x$, $f_2(x,z_2)$ is strictly decreasing with respect to $z_2,$ it yields that
\[
h^1_2(x,\lambda)=  f_2(x,\frac{\xi'(1)x}{\xi'(x)}-1)>f_2(x, \frac{\xi'(1)-\xi'(x)}{\xi''(x)\cdot(1-x)}-1)= h^2_2(x,\lambda).    
\] 

We now prove item (3). We split this item into two cases based on Lemma $\ref{tx}$. 

\textbf{Case I}: $t(x)$ has exactly 2 roots $q_1 \in (0,1)$ and $q_2 \in (1,+\infty)$. In this case, it holds that $\lim_{x \rightarrow 0^+} t(x)>0$ and $\lim_{x \rightarrow 1^-} t(x)<0. $ Thus, for $x \in (0,q_1)$, $t(x)>0$ while for $x \in (q_1,1)$, $t(x)<0,$ which implies that $\bar{q}_1=q_1$ and $\bar{q}_2=1.$

\textbf{Case II}: $t(x)$ has either 2 or 0 roots in $(0,1).$ In this case, we have $\lim_{x \rightarrow 0^+} t(x)>0$ and $\lim_{x \rightarrow 1^-} t(x)>0.$  If $t(x)$ has no root in $(0,1)$, then $t(x)>0$ for $x \in (0,1).$ If $t(x)$ has two roots $q_1$ and $q_2$ in $(0,1)$, then $t(x)<0$ for $x \in (q_1,q_2)$ and $t(x)\geq0$ for $x \in (0,q_1] \cup [q_2,+\infty).$ Then it yields that $\bar{q}_1=q_1$ and $\bar{q}_2=q_2.$

Now recall that
the sign of $\frac{d}{dx} h_1^2(x,\lambda)$ is determined by $t^2_1(x)$
and $h_1^2(x,\lambda)$ have either 2 or 0 critical points on $(0,+\infty)$.
If $h^1_2(x,\lambda)$ has two critical points, denoted by $s_1$ and $s_2$, then for $x \in (0,s_1)\cup (s_2,+\infty),$ $h_1^2(x,\lambda)$ is strictly decreasing while for $x \in (s_1,s_2),$ $h_1^2(x,\lambda)$ is strictly increasing. Thus $h_1^2(x,\lambda)$ has at most 2 roots in (0,1).

The limit $\lim_{x \rightarrow 1}h^2_1(x,\lambda)=\Psi(\lambda)$ was shown  in the proof of Theorem \ref{coro1rsb}.

Now we turn to item (5). By the proof of Lemma \ref{hs}, the sign of $\frac{d}{dx}h^2_2(x,\lambda)$ is the same as the sign of $m(x).$  Note that
\begin{eqnarray*}
m'''(1)&=& 2\xi''(1)\xi''''(1)-3\xi'''(1)^2\\
&=&s^3(s-1)^2(s-2)(1-\lambda)^2+p^3(p-1)^2(p-2)\lambda^2\\
&&-2p(p-1)s(s-1)\Big[(p-2)(p-3)+(s-2)(s-3)-3(p-2)(s-2) \Big]\lambda(1-\lambda),
\end{eqnarray*}
then $\lambda_{2 \rightarrow 1F}$ and $\lambda'_{2\rightarrow1F}$ are the two zeros of $2\xi''(1)\xi''''(1)-3\xi'''(1)^2=0,$ which implies that 
$\frac{d^2}{dx^2}h^2_2(1,\lambda)>0$ if and only if $\lambda \in (\lambda_{2 \rightarrow 1F},\lambda'_{2 \rightarrow 1F})$.

\end{proof}

By Lemmas \ref{tx} and \ref{h4}, $h^1_1(x,\lambda)$ has zeros in $(0,1)$ if and only if $\bar{q}_1,\bar{q}_2$ exists and $h^1_1(\bar{q}_2,\lambda)>0.$ At the same time, $h^1_1(x, \lambda)$ has at most 2 zeros on (0,1). We denote the smaller one by $q_1^1.$ Notice that $q_1^1 \in (\bar{q}_1,\bar{q}_2).$

Since $h^2_1(\bar{q}_1,\lambda)=h^1_1(\bar{q}_1,\lambda)<0$ and $h^2_1(\bar{q}_2,\lambda)=h^1_1(\bar{q}_2,\lambda)>0,$ there exists a zero of $h^2_1(x,\lambda)$ in $(\bar{q}_1,\bar{q}_2).$ The uniqueness of this zero is guaranteed by Lemma \ref{h4}. We denote this zero of $h^2_1(x,\lambda)$ by $q^2_1.$

Similarly, $h^1_2(x,\lambda)$ has zeros in $(0,1)$ if and only if $\bar{q}_1, \bar{q}_2$ exists and $h^1_2(\bar{q}_1,\lambda)>0.$ If $\bar{q}_2<1,$ $h^1_2(x,\lambda)$ has two zeros in $(0,1)$ and we denote the larger zero by $q^1_2.$ If $\bar{q}_2=1,$ we set $q^1_2=1.$ It holds that $q^1_2 \in [\bar{q}_1,\bar{q}_2].$

Lastly, we turn to the definition of $q^2_2(\lambda)$ and consider the following 3 cases: $(i)$ When $\bar{q}_1<\bar{q}_2<1$, $h^1_2(\bar{q}_1,\lambda)>0$ and $h^1_2(\bar{q}_2,\lambda)<0.$ Combining with the fact that $h^2_2(\bar{q}_1,\lambda)=h^1_2(\bar{q}_1,\lambda)$ and $h^1_2(\bar{q}_2,\lambda)=h^2_2(\bar{q}_2,\lambda)$, there exists zeros of $h^2_2(x,\lambda)$ in $(\bar{q}_1,\bar{q}_2).$ We denote the largest zero in $(\bar{q}_1,\bar{q}_2)$ by $q_2^2.$

$(ii)$ When $\bar{q}_1<\bar{q}_2=1$ and $\frac{d^2}{dx^2}h^2_2(x,\lambda)|_{x=1}>0,$ it still holds that $h^1_2(\bar{q}_1)>0$. Since $\frac{d^2}{dx^2}h^2_2(x,\lambda)|_{x=1}>0,$ there exists zeros of $h^2_2(x,\lambda)$ in $(\bar{q}_1,1)$ by intermediate value theorem. We denote the largest zero in $(\bar{q}_1,1)$ by $q^2_2(\lambda).$

$(iii)$ When $\bar{q}_1<\bar{q}_2=1$ and $\frac{d^2}{dx^2}h^2_2(x,\lambda)|_{x=1}>0,$ we set $q^2_2(\lambda)=1.$

Based on the definitions above and assuming the existence, we notice that it always holds that $q^2_2<q^1_2$ and $q_1^1 < q^2_1.$ Moreover, if $q^1_1$ exists, then $q^2_1$ also exists. If $q^1_2$ exists, then $q^2_2$ also exists.

When $q^1_1,q^2_1,q^1_2$ and $q^2_2$ exist, there are altogether three cases of their relation:
\begin{eqnarray*}
&(i)& q^1_2<q^1_1\\
&(ii)& q^1_2>q^1_1 \text{ and } q_2^2 <q^2_1, \\
&(iii)& q_2^2 >q^2_1.
\end{eqnarray*}

Now we are ready to prove Theorem \ref{coro12rsb}$(ii)$ and Theorem \ref{coro12frsb}$(iii)$.
\begin{proof}[Proof of Theorem \ref{coro12rsb}$(ii)$ and Theorem \ref{coro12frsb}$(iii)$]
We first prove that the model $\xi$ is 2RSB if and only if $q^1_1(\lambda),q^2_1(\lambda),q^1_2(\lambda),q^2_2(\lambda)$ exist and satisfy $q^1_2(\lambda)>q^1_1(\lambda) \text{ and } q_2^2(\lambda) <q^2_1(\lambda)$.

There are altogether 4 possibilities of their relations as follows:
\begin{eqnarray*}
&(i)& q^2_2<q^1_1<q^1_2<q^2_1, \\
&(ii)& q^1_1<q^2_2<q^1_2<q^2_1, \\
&(iii)& q^1_1<q^2_2<q^2_1<q^1_2, \\
&(iv)& q^2_2<q^1_1<q^2_1<q^1_2. 
\end{eqnarray*}
We now prove by cases. We will omit the dependence on $\lambda$ to easy the notation in the following paragraphs.

\textbf{Case $(i)$}. $q^2_2<q^1_1<q^1_2<q^2_1.$

From Figure \ref{F1}, we notice that $h_1^1(q^1_2)>0$ and $h^2_1(q^1_2)<0.$ Recall the definition $h^1_1(q^1_2)= (q^1_2)^{-2} \cdot f_1(q^1_2, \frac{\xi'(1)q^1_2}{\xi'(q^1_2)}-1) $ and $h^2_1(q^1_2)=(q^1_2)^{-2}\cdot f_1(q^1_2, \frac{\xi'(1)-\xi'(q^1_2)}{\xi''(x)\cdot(1-q^1_2)}-1) $, then there exists $z^1_1 \in \Big (\frac{\xi'(1)q^1_2}{\xi'(q^1_2)}-1, \frac{\xi'(1)-\xi'(q^1_2)}{\xi''(x)\cdot(1-q^1_2)}-1 \Big)$ such that $f_1(q^1_2,z^1_1)=0.$

Similarly, by $h^1_2(q^1_1)>0$ and $h^2_2(q^1_1)<0$, it is equivalent to say that $f_2(q^1_1,\frac{\xi'(1)q^1_1}{\xi'(q^1_1)}-1)>0$ and $f_2(q^1_1,\frac{\xi'(1)-\xi'(q^1_1)}{\xi''(q^1_1)(1-q^1_1)}-1)<0.$ Then there exists $z^1_2 \in \Big (\frac{\xi'(1)q^1_1}{\xi'(q^1_1)}-1, \frac{\xi'(1)-\xi'(q^1_1)}{\xi''(x)\cdot(1-q^1_1)}-1 \Big)$ such that $f_2(q^1_1,z^1_2)=0.$ 

Moreover, we have $h^1_2(q^1_2)=f_2(q_2^1,\frac{\xi'(1)q_2^1}{\xi'(q_2^1)}-1)=0$ and $h^1_1(q^1_1)=f_1(q_1^1,\frac{\xi'(1)q^1_1}{\xi'(q^1_1)}-1)=0$.

Based on the facts above,  $(q_1^1,\frac{\xi'(1)q^1_1}{\xi'(q^1_1)}-1)$ and $(q^1_2,z^1_1)$ the zeros of $f_1(q,z_2)$ while $(q^1_1,z^1_2)$ and $(q_2^1,\frac{\xi'(1)q_2^1}{\xi'(q_2^1)}-1)$ are the zeros of $f_2(q,z_2)$. By continuity, there exists $(q^1_*,z^1_*)\in [q^1_1,q^1_2] \times \mathbb{R}$ satisfying that $f_1(q,z_2)=f_2(q,z_2)=0.$ Since $q^1_* \in [q^1_1,q^1_2], $ we have $f_2(q^1_*,\frac{\xi'(1)q^1_*}{\xi'(q^1_*)}-1)>0$ and $f_2(q^1_*, \frac{\xi'(1)-\xi'(q^1_*)}{\xi''(q^1_*)(1-(q^1_*))}-1)<0,$ which implies there exists point $a_1 \in \Big(\frac{\xi'(1)q^1_*}{\xi'(q^1_*)}-1, \frac{\xi'(1)-\xi'(q^1_*)}{\xi''(q^1_*)(1-(q^1_*))}-1\Big)$ such that $f_2(q^1_*,a_1)=0.$ Since for any $q>0,$ the zero of $f_2(q,\cdot)=0$ exists and is unique. We then obtain that $a_1=z^1_*$ and $z^1_* \in  \Big(\frac{\xi'(1)q^1_*}{\xi'(q^1_*)}-1, \frac{\xi'(1)-\xi'(q^1_*)}{\xi''(q^1_*)(1-(q^1_*)}-1)\Big).$

Therefore $q^1_*$ and $z^1_*$ are the desired $q$ and $z_2$ satisfying the conditions in Theorem \ref{2RSB} and making the model $\xi(x)$ to be 2RSB.

\begin{figure}[H]
\centering
	\begin{subfigure}
	\centering
	\includegraphics[width=8cm,height=6cm]{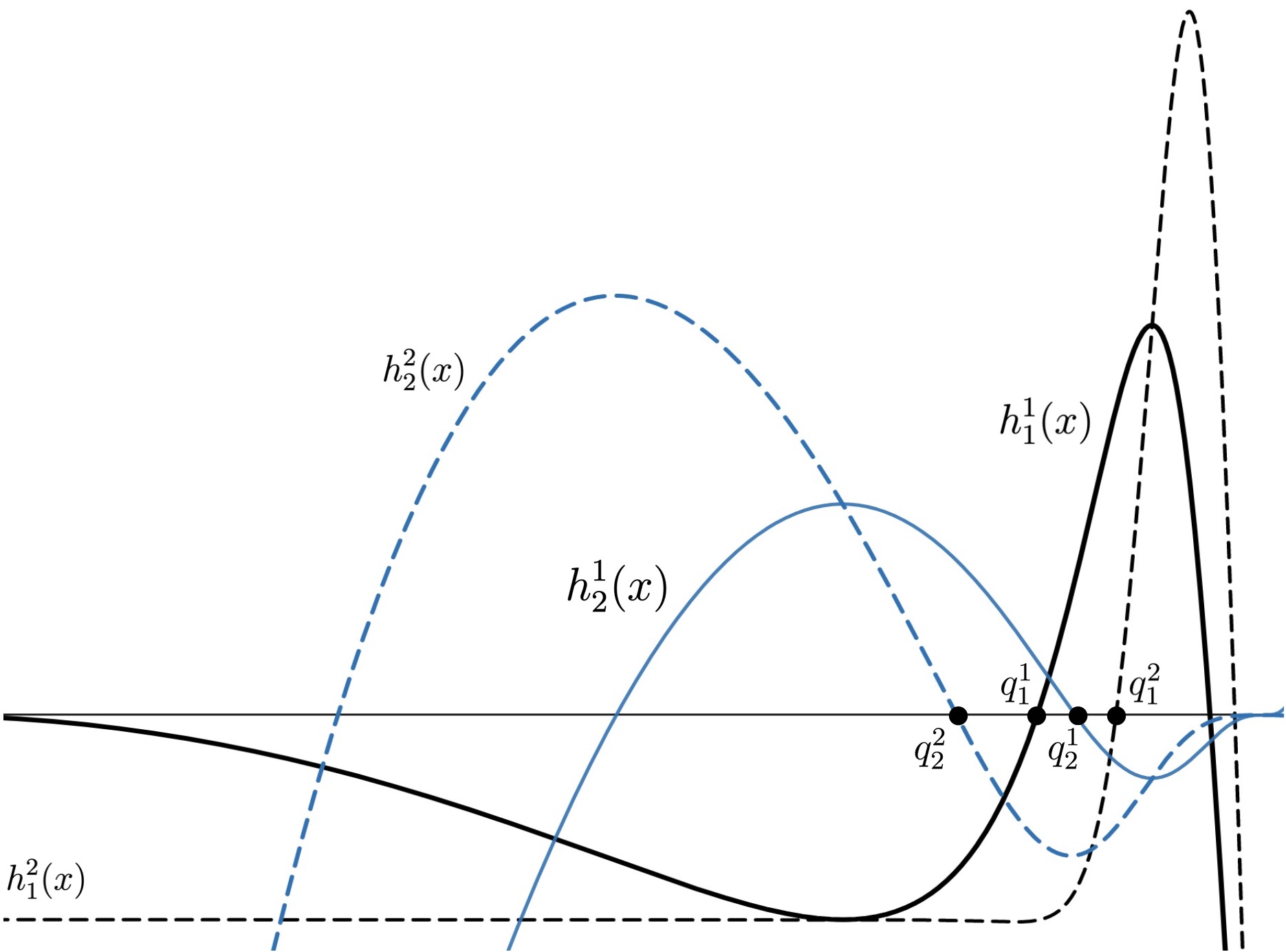} \label{F1}
\end{subfigure}
\begin{subfigure}
\centering
	\includegraphics[width=8cm,height=6cm]{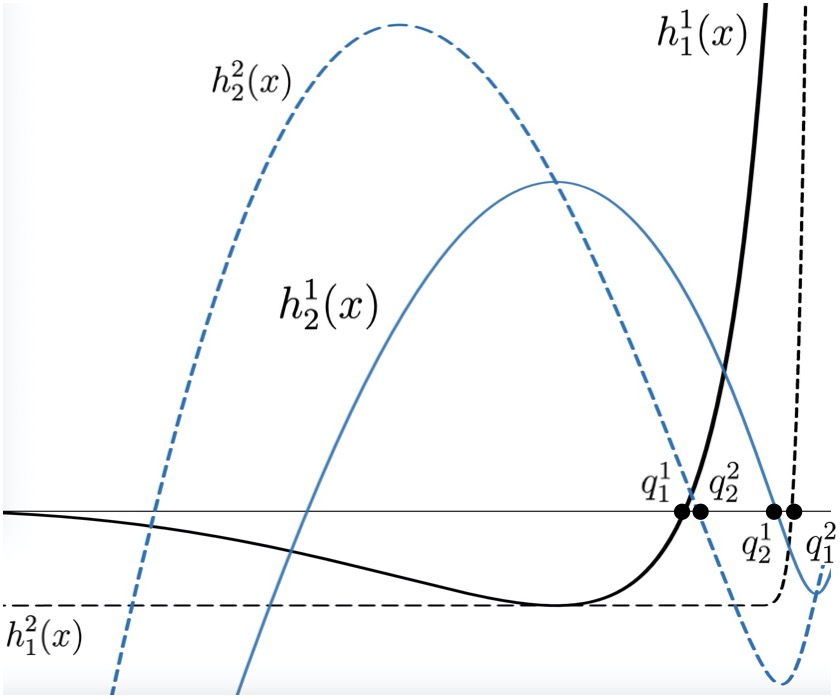} \label{F2}
	\end{subfigure}
\caption{Sketch graphs of Case $(i)$ and (ii).  The left subfigure describes Case $(i)$. The right subfigure describes Case $(ii)$. The black solid curve represents $h^1_1(x,\lambda)$.   The black dotted curve represents $h^2_1(x,\lambda).$ The blue solid curve represents $h^1_2(x,\lambda).$ The blue dotted curve represents $h^2_2(x,\lambda).$  }	
\label{figure2}
	\end {figure}

\textbf{Case $(ii)$}. $q^1_1<q^2_2<q^1_2<q^2_1$.





From Figure \ref{F2}, we note that $h^1_1(q^1_2)>0$ and $h^2_1(q^1_2)<0,$ which implies that $f_1(q^1_2, \frac{\xi'(1)q^1_2}{\xi'(q^1_2)}-1)>0$ and $f_1(q^1_2,\frac{\xi'(1)-\xi'(q^1_2)}{\xi''(q^1_2)(1-q^1_2)}-1)<0.$ Then there exists $z^2_1 \in \Big( \frac{\xi'(1)q^1_2}{\xi'(q^1_2)}-1, \frac{\xi'(1)-\xi'(q^1_2)}{\xi''(q^1_2)(1-q^1_2)}-1\Big)$ such that $f_1(q^1_2,z^2_1)=0.$

Similarly, by $h^1_1(q^2_2)>0$ and $h^2_1(q^2_2)<0.$ we have $f_1(q^2_2, \frac{\xi'(1)q^2_2}{\xi'(q^2_2)}-1)>0$ and $f_1(q^2_2,\frac{\xi'(1)-\xi'(q^2_2)}{\xi''(q^2_2)(1-q^2_2)}-1)<0$. Then there exists $z^2_2 \in \Big( \frac{\xi'(1)q^2_2}{\xi'(q^2_2)}-1, \frac{\xi'(1)-\xi'(q^2_2)}{\xi''(q^2_2)(1-q^2_2)}-1\Big)$ such that $f_1(q^2_2,z^2_2)=0.$

Moreover, we have that $h^2_2(q^2_2)=h^1_2(q^1_2)=0,$ which implies that $f_2(q_2^2,\frac{\xi'(1)-\xi'(q^2_2)}{\xi''(q^2_2)(1-q^2_2)}-1)=f_2(q^1_2, \frac{\xi'(1)q^1_2}{\xi'(q^1_2)}-1)=0$.

Based on the reasoning above, $(q_2^2,\frac{\xi'(1)-\xi'(q^2_2)}{\xi''(q^2_2)(1-q^2_2)}-1)$ and $(q^1_2, \frac{\xi'(1)q^1_2}{\xi'(q^1_2)}-1)$ are two zeros of $f_1(q,z_2)$ while $(q^2_2,z^2_2)$ and $(q^1_2,z^2_1)$ are two zeros of $f_2(q,z_2).$ By continuity, there exists $ [q^2_*,z^2_*] \in [q^2_2,q^1_2] \times \mathbb{R}$ satisfying that $f_1(q,z_2)=f_2(q,z_2)=0.$ Since $q^2_* \in [q^2_2,q^1_2],$ it yields that $f_2(q^2_*,\frac{\xi'(1)q^2_*}{\xi'(q^2_*)}-1)>0$ and $f_2(q^2_*, \frac{\xi'(1)-\xi'(q^2_*)}{\xi''(q^2_*)(1-q^2_*)}-1)<0.$ By intermediate value theorem, there exists $a_2 \in \Big(\frac{\xi'(1)q^2_*}{\xi'(q^2_*)}-1, \frac{\xi'(1)-\xi'(q^2_*)}{\xi''(q^2_*)(1-q^2_*)}-1 \Big)$ such that $f_2(q^2_*,a_2)=0.$ By uniqueness of zeros of $f_2(q,\cdot),$ we have that $a_2=z^2_*$ and $z^2_* \in \Big(\frac{\xi'(1)q^2_*}{\xi'(q^2_*)}-1, \frac{\xi'(1)-\xi'(q^2_*)}{\xi''(q^2_*)(1-q^2_*)}-1 \Big).$ 

Therefore $q^2_*$ and $z^2_*$ satisfy the conditions of $q$ and $z_2$ in Theorem \ref{2RSB} making the model $\xi(x)$ to be 2RSB.

\textbf{Case $(iii)$}. $q^1_1<q^2_2<q^2_1<q^1_2$.

From Figure \ref{F3}, we observe that $h^1_2(q^2_1)>0$ and $h^2_2(q^2_1)<0,$ which implies that $f_2(q^2_1, \frac{\xi'(1)q^2_1}{\xi'(q^2_1)}-1)>0$ and $f_2(q^2_1,\frac{\xi'(1)-\xi'(q^2_1)}{\xi''(q^2_1)(1-q^2_1)}-1)<0.$ Then by intermediate value theorem, there exists $z^3_1 \in \Big( \frac{\xi'(1)q^2_1}{\xi'(q^2_1)}-1,\frac{\xi'(1)-\xi'(q^2_1)}{\xi''(q^2_1)(1-q^2_1)}-1 \Big)$ such that $f_2(q^2_1,z^3_1)=0.$

Similarly, since $h^1_1(q^2_2)>0$ and $h^2_1(q^2_2)<0$, we have that $f_1(q^2_2,\frac{\xi'(1)q^2_2}{\xi'(q^2_2)}-1)>0$ and $f_1(q^2_2,\frac{\xi'(1)-\xi'(q^2_2)}{\xi''(q^2_2)(1-q^2_2)}-1)<0.$ Thus there exists $z^3_2 \in \Big( \frac{\xi'(1)q^2_2}{\xi'(q^2_2)}-1 , \frac{\xi'(1)-\xi'(q^2_2)}{\xi''(q^2_2)(1-q^2_2)}-1 \Big)$ such that $f_1(q^2_2,z^3_2)=0.$ Also, it holds that $f_2(q_2^2,\frac{\xi'(1)-\xi'(q^2_2)}{\xi''(q^2_2)(1-q^2_2)}-1)=f_1(q_1^2, \frac{\xi'(1)-\xi'(q^2_1)}{\xi''(q^2_1)(1-q^2_1)}-1)=0.$

Sum the reasoning above up,$(q^2_2,z^3_2),(q_1^2, \frac{\xi'(1)-\xi'(q^2_1)}{\xi''(q^2_1)(1-q^2_1)}-1)$ are two zeros of $f_1(q,z_2)$ and $(q_2^2,\frac{\xi'(1)-\xi'(q^2_2)}{\xi''(q^2_2)(1-q^2_2)}-1), (q^2_1,z^3_1)$ are two zeros of $f_2(q,z_2)$. Still by continuity, there exists $[q^3_*,z^3_*] \in [q^2_2,q^2_1] \times \mathbb{R}$  satisfying that $f_1(q,z_2)=f_2(q,z_2)=0.$ Since $q^3_* \in [q^2_2,q^2_1],$ it yields that $f_2(q^3_*,\frac{\xi'(1)q^3_*}{\xi'(q^3_*)}-1)>0$ and $f_2(q^3_*, \frac{\xi'(1)-\xi'(q^3_*)}{\xi''(q^3_*)(1-q^3_*)}-1)<0.$ By intermediate value theorem, there exists $a_3 \in \Big(\frac{\xi'(1)q^3_*}{\xi'(q^3_*)}-1, \frac{\xi'(1)-\xi'(q^3_*)}{\xi''(q^3_*)(1-q^3_*)}-1 \Big)$ such that $f_2(q^3_*,a_3)=0.$ By uniqueness of $f_2(q,\cdot),$ we have that $a_3=z^3_*$ and $z^3_* \in \Big(\frac{\xi'(1)q^3_*}{\xi'(q^3_*)}-1, \frac{\xi'(1)-\xi'(q^3_*)}{\xi''(q^3_*)(1-q^3_*)}-1 \Big).$ 

Hence $q^3_*$ and $z^3_*$ are the $q$ and $z$ making the model $\xi(x)$ to be 2RSB.

\begin{figure}[H]
\centering
\begin{subfigure}
	\centering
\includegraphics[width=8cm,height=6cm]{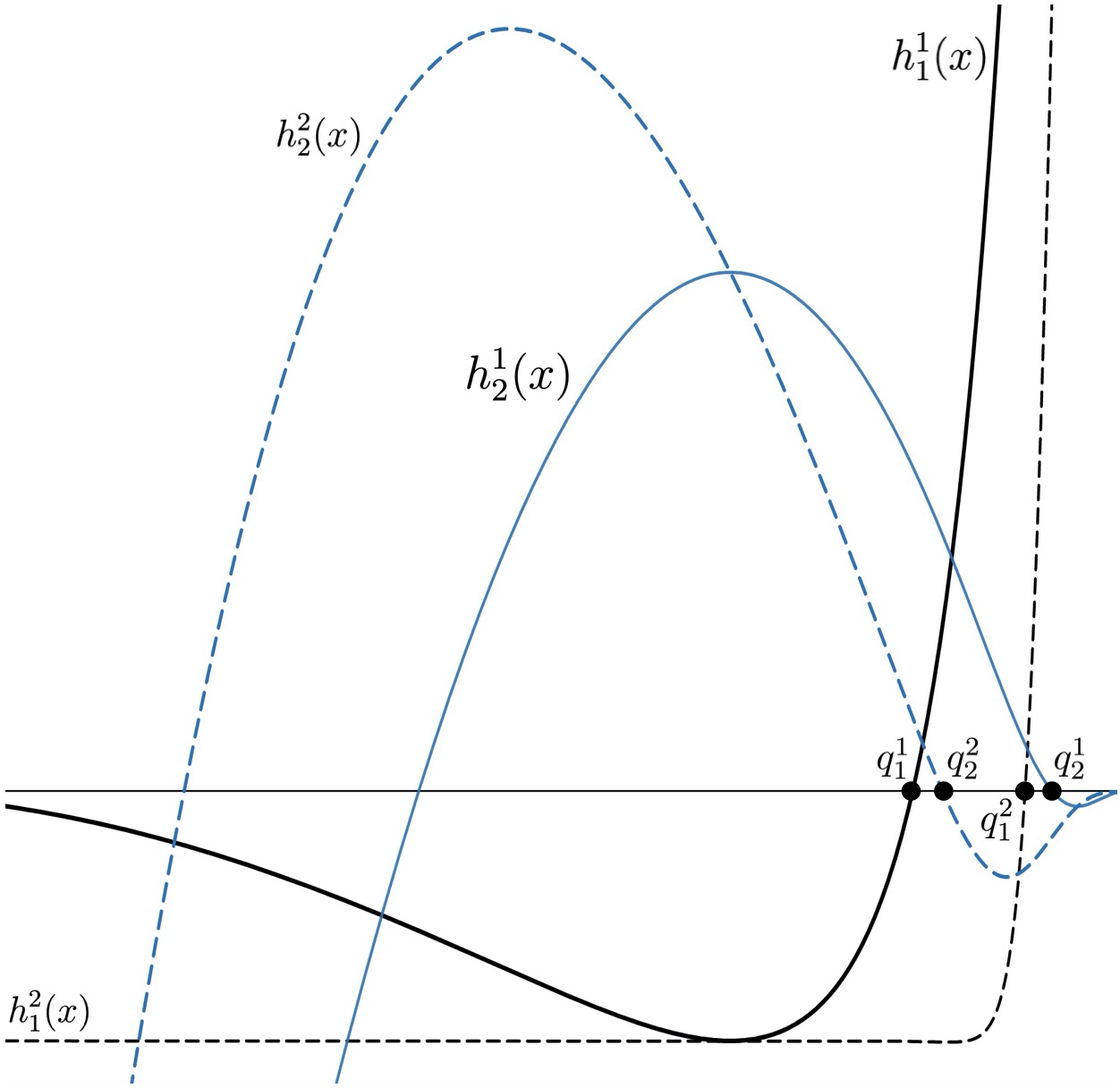} \label{F3}
	\end{subfigure}
	\begin{subfigure}
	\centering
\includegraphics[width=8cm,height=6cm]{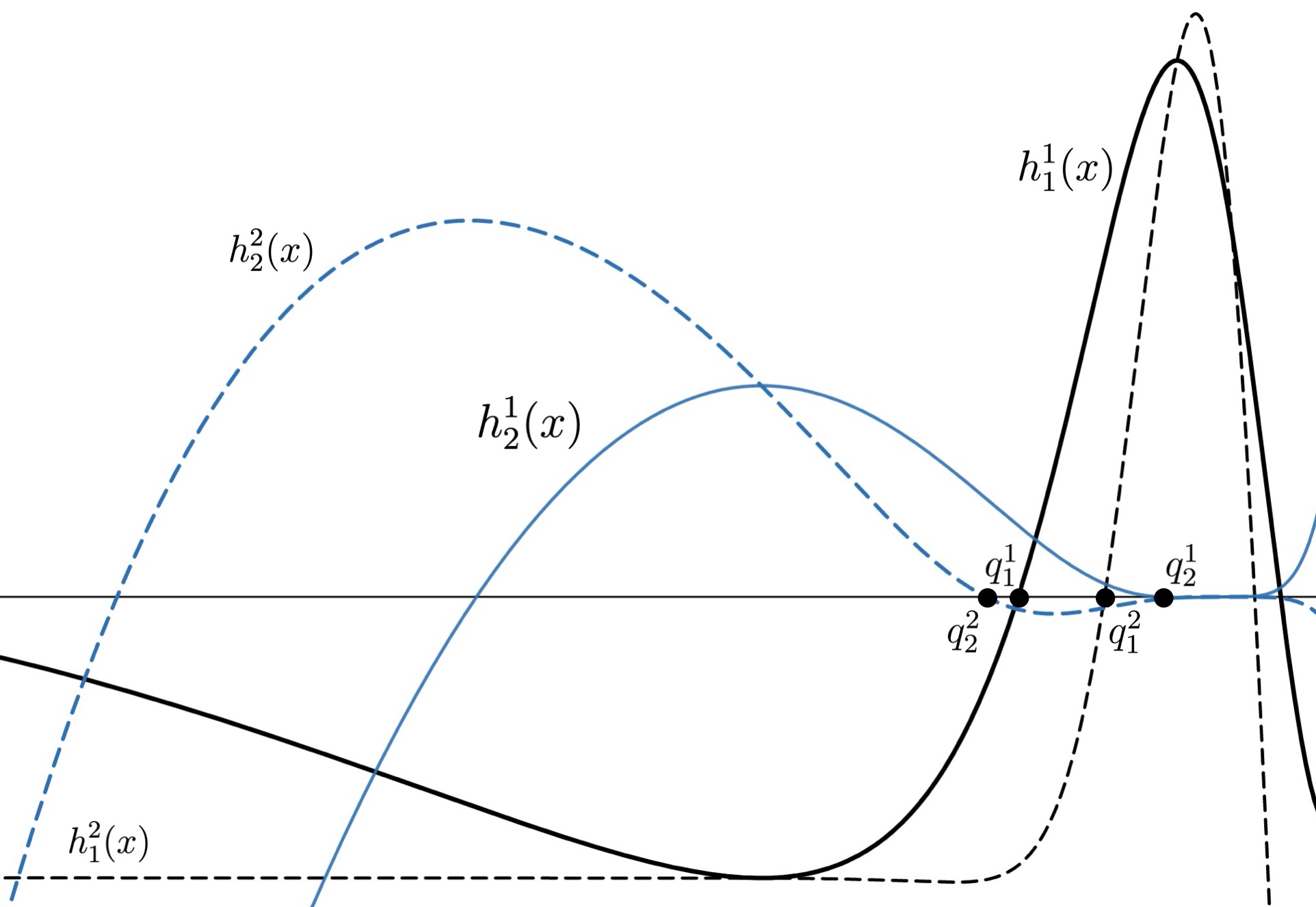} \label{F4}
\end{subfigure}
	\caption{Sketch graphs of Case $(iii)$ and $(iv)$.  The left subfigure describes Case $(iii)$. The right subfigure describes Case $(iv)$.  The coloring of the curves are the same as in Figure \ref{figure2}.  }	
\label{f3}
	\end {figure}

\textbf{Case $(iv)$}. $q^2_2<q^1_1<q^2_1<q^1_2.$

From Figure \ref{F4}, we first notice that $h^1_2(q^2_1)>0$ and $h^2_2(q^2_1)<0, $ which implies that $f_2(q^2_1, \frac{\xi'(1)q^2_1}{\xi'(q^2_1)}-1)>0$ and $f_2(q^2_1,\frac{\xi'(1)-\xi'(q^2_1)}{\xi''(q^2_1)(1-q^2_1)}-1)<0.$ Then by intermediate value theorem, there exists $z^4_1 \in \Big( \frac{\xi'(1)q^2_1}{\xi'(q^2_1)}-1,\frac{\xi'(1)-\xi'(q^2_1)}{\xi''(q^2_1)(1-q^2_1)}-1 \Big)$ such that $f_2(q^2_1,z^4_1)=0.$

Also, since $h^1_2(q^1_1)>0$ and $h^2_2(q^1_1)<0$, it is equivalent to say that $f_2(q^1_1,\frac{\xi'(1)q^1_1}{\xi'(q^1_1)}-1)>0$ and $f_2(q^1_1,\frac{\xi'(1)-\xi'(q^1_1)}{\xi''(q^1_1)(1-q^1_1)}-1)<0.$ Then there exists $z^4_2 \in \Big (\frac{\xi'(1)q^1_1}{\xi'(q^1_1)}-1, \frac{\xi'(1)-\xi'(q^1_1)}{\xi''(x)\cdot(1-q^1_1)}-1 \Big)$ such that $f_2(q^1_1,z^4_2)=0.$  Moreover, we have  $f_1(q_1^1,\frac{\xi'(1)q^1_1}{\xi'(q^1_1)}-1)=f_1(q_1^2, \frac{\xi'(1)-\xi'(q^2_1)}{\xi''(q^2_1)(1-q^2_1)}-1)=0$.

Based on the reasoning above, we obtain that there exists two continuous curves with one connecting $(q_1^1,\frac{\xi'(1)q^1_1}{\xi'(q^1_1)}-1),(q_1^2, \frac{\xi'(1)-\xi'(q^2_1)}{\xi''(q^2_1)(1-q^2_1)}-1)$ and the other connecting $(q^1_1,z^4_2),(q^2_1,z^4_1).$ By continuity, there is at least one common element $[q^4_*,z^4_*] \in [q^1_1,q^2_1] \times \mathbb{R}$ in the zero sets of $f_1(q,z_2)$ and $f_2(q,z_2).$ Since $q^4_* \in [q^1_1,q^2_1],$ we obtain that $f_2(q^4_*,\frac{\xi'(1)q^4_*}{\xi'(q^4_*)}-1)>0$ and $f_2(q^4_*, \frac{\xi'(1)-\xi'(q^4_*)}{\xi''(q^4_*)(1-q^4_*)}-1)<0.$ By intermediate value theorem, there exists $a_4 \in \Big(\frac{\xi'(1)q^4_*}{\xi'(q^4_*)}-1, \frac{\xi'(1)-\xi'(q^4_*)}{\xi''(q^4_*)(1-q^4_*)}-1 \Big)$ such that $f_2(q^4_*,a_4)=0.$ By uniqueness of $f_2(q,\cdot),$ we have that $a_4=z^4_*$ and $z^4_* \in \Big(\frac{\xi'(1)q^4_*}{\xi'(q^4_*)}-1, \frac{\xi'(1)-\xi'(q^4_*)}{\xi''(q^4_*)(1-q^4_*)}-1 \Big).$  Thus $q^4_*$ and $z^4_*$ are the $q$ and $z$ making the model $\xi(x)$ to be 2RSB.

Now we prove that if the system of equations $\eqref{eqnh2}$ has a solution, then it must be unique. By Theorem \ref{2RSB}, it is equivalent for us to show that when $q,\lambda \in (0,1)$ satisfies that $1+z_2=\frac{\xi'(1)-\xi'(q)}{\xi''(q)(1-q)},$ it holds that $\frac{d}{d\lambda} \Big( 1+z_2-\frac{\xi'(1)-\xi'(q)}{\xi''(q)(1-q)} \Big)>0.$ 

It is not difficult to see that
\begin{align*}
\frac{dz_2}{dq}=\frac{\xi''(q)(1-q)-[\xi'(1)-\xi'(q)+(1-q)\xi''(q)]\frac{1-\xi(q)-\xi'(q)(1-q)}{(1-q)[\xi'(1)-\xi'(q)]}}{(1-q)[\xi'(1)-\xi'(q)]\big( \frac{2+z_2}{z_2^3} \log(1+z_2)-\frac{2}{z^2_2}\big)},
\end{align*}
which implies that when $1+z_2=\frac{\xi'(1)-\xi'(q)}{\xi''(q)(1-q)},$
\begin{align*}
\frac{dz_2}{dq}-\Big(\frac{\xi'(1)-\xi'(q)}{\xi''(q)(1-q)}\Big)'_q &= \frac{\frac{\xi''(q)(1-q)[(1-q)\xi'(1)-1+\xi(q)]-[\xi'(1)-\xi'(q)][1-\xi(q)-\xi'(q)(1-q)]}{(1-q)[\xi'(1)-\xi'(q)]}}{\frac{\xi''(q)(1-q) \Big\{ \xi''(q)(1-q)[(1-q)\xi'(1)-1+\xi(q)]-[\xi'(1)-\xi'(q)][1-\xi(q)-\xi'(q)(1-q)] \Big \}}{-[\xi'(1)-\xi'(q)][\xi'(1)-\xi'(q)-\xi''(q)(1-q)]}} \\
&+\frac{[\xi'(1)-\xi'(q)](1-q)\xi'''(q)-\xi''(q)[\xi'(1)-\xi'(q)-\xi''(q)(1-q)]}{\xi''(q)^2(1-q)^2}\\
&=\frac{\xi'''(q)[\xi'(1)-\xi'(q)](1-q)-2\xi''(q)[\xi'(1)-\xi'(q)-\xi''(q)(1-q)]}{\xi''(q)^2(1-q)^2}\\
&=\frac{m(q)}{\xi''(q)^2(1-q)^2}.
\end{align*}

Then we compute the derivative of $\frac{d}{d\lambda} \Big( 1+z_2-\frac{\xi'(1)-\xi'(q)}{\xi''(q)(1-q)} \Big)$ with respect to $\lambda$ when  $1+z_2=\frac{\xi'(1)-\xi'(q)}{\xi''(q)(1-q)}$ to obtain:
\begin{align*}
\frac{d}{d\lambda} &\Big( 1+z_2-\frac{\xi'(1)-\xi'(q)}{\xi''(q)(1-q)} \Big) \Big|_{1+z_2=\frac{\xi'(1)-\xi'(q)}{\xi''(q)(1-q)}} =\frac{dz_2}{dq} \cdot q'_\lambda - \Big(\frac{\xi'(1)-\xi'(q)}{\xi''(q)(1-q)}\Big)'_q \cdot q'_\lambda - \Big(\frac{\xi'(1)-\xi'(q)}{\xi''(q)(1-q)}\Big)'_\lambda\\
&=\frac{m(q)}{\xi''(q)^2(1-q)^2}\cdot q'_\lambda -\frac{ps}{\xi''(q)^2(1-q)} \big[ (s-1)q^{s-2}-(s-p)q^{p+s-3}-(p-1)q^{p-2} \big].
\end{align*}

Byy \eqref{f1}, it holds that
 $ q'_\lambda=-\frac{\big( f_1(q,z_2) \big)'_{\lambda} }{\Big(\big( f_1(q,z_2) \big)'_{z_2} \cdot \frac{d z_2}{dq} + \big( f_1(q,z_2) \big)'_{q} \Big)}.$
 Also,
 \begin{align*}
 \Big(\big( f_1(q,z_2) \big)'_{z_2} \cdot \frac{d z_2}{dq} + \big( f_1(q,z_2) \big)'_{q} \Big)
 &=\Big( -\frac{[q\xi'(q)-\xi(q)](1-q)}{q^2[\xi'(1)-\xi'(q)]}+\frac{1}{1+z_2}-\frac{\xi(q)[\xi'(1)-\xi'(q)]}{1+z_2)^2\xi'(q)^2(1-q)} \Big) \cdot \frac{dz_2}{dq}\\
 &+\frac{[\xi'(q)-q\xi''(q)][q\xi'(q)^2-q\xi(q)\xi''(q)-\xi(q)\xi'(q)]}{\xi'(q)^3\xi''(q)q^3(1-q)[\xi'(1)-\xi'(q)]} \\
 &\cdot \Big[ q(1-q)\xi''(q)[\xi'(q)-2\xi'(1)]+[\xi'(1)-\xi'(q)]\xi'(q)(2-q) \Big]\\
 &=\frac{2[q\xi''(q)-\xi'(q)][q\xi'(q)^2-q\xi(q)\xi''(q)-\xi(q)\xi'(q)]}{\xi''(q)\xi'(q)^3[\xi'(1)-\xi'(q)]q^3(1-q)} \cdot \\
 &\Big[ \xi'(1)\xi''(q)q(1-q)-\xi'(q)[\xi'(1)-\xi'(q)] \Big]\\
  &=\frac{2[q\xi''(q)-\xi'(q)][q\xi'(q)^2-q\xi(q)\xi''(q)-\xi(q)\xi'(q)]}{\xi''(q)\xi'(q)^3[\xi'(1)-\xi'(q)]q^3(1-q)} \cdot t(q)
 \end{align*}
 and
 \begin{align*}
 &-\big( f_1(q,z_2) \big)'_{\lambda}=-\frac{[q\xi''(q)-\xi'(q)][q\xi'(q)^2-q\xi(q)\xi''(q)-\xi(q)\xi'(q)]\big[q\xi'(q)\xi'''(q)-2\xi''(q)[q\xi''(q)-\xi'(q)]\big]}{(s-p)\lambda(1-\lambda)q^2\xi'(q)^3\xi''(q)^2}\\
 &\qquad -\frac{ps[(p-1)q^{p-2}-(s-1)q^{s-2}+(s-p)q^{p+s-3}]}{\xi''(q)^2(1-q)} \frac{[q\xi''(q)-\xi'(q)][q\xi'(q)^2-q\xi(q)\xi''(q)-\xi(q)\xi'(q)]}{\xi'(q)^2[\xi'(1)-\xi'(q)]q^2}\\
 &=-[q\xi''(q)-\xi'(q)][q\xi'(q)^2-q\xi(q)\xi''(q)-\xi(q)\xi'(q)] \cdot\\
 & \Big\{ \frac{q\xi'(q)\xi'''(q)-2\xi''(q)[q\xi''(q)-\xi'(q)]}{(s-p)\lambda(1-\lambda)q^2\xi'(q)^3 \xi''(q)^2}+\frac{ps\big[(p-1)q^{p-2}-(s-1)q^{s-2}+(s-p)q^{p+s-3}}{\xi'(q)^2\xi''(q)^2[\xi'(1)-\xi'(1)]q^2}  \Big\}.
 \end{align*}
Then we define $p(x)=\big[x\xi'(x)\xi'''(x)-2\xi''(x)\big(x\xi''(x)-\xi'(x)\big)\big]\cdot [\xi'(1)-\xi'(x)]+ps(s-p)\lambda(1-\lambda)[(p-1)x^{p-2}-(s-1)x^{s-2}+(s-p)x^{p+s-3}]\xi'(x).$ If we combine the computations above, we obtain
\begin{align*}
-2(s-p)&\lambda(1-\lambda)\xi''(q)^2t(q)(1-q)\cdot\frac{d}{d\lambda} \Big( 1+z_2-\frac{\xi'(1)-\xi'(q)}{\xi''(q)(1-q)} \Big) \Big|_{1+z_2=\frac{\xi'(1)-\xi'(q)}{\xi''(q)(1-q)}} \\
&= q p(q)m(q)+2ps(s-p)\lambda(1-\lambda)[(s-1)q^{s-2}-(s-p)q^{p+s-3}-(p-1)q^{p-2}]\xi''(q)t(q).
\end{align*}

We recall that for $x \in (\bar{q}_1,\bar{q}_2), t(x)<0.$ If we define the largest zero of $m(x)$ in $[0,1]$ by $\tilde{q}_2$ then $\tilde{q}_2 \in (\bar{q}_1,\bar{q}_2).$
Since $1+z_2=\frac{\xi'(1)-\xi'(q)}{\xi''(q)(1-q)},$ it holds that $t(x)<0,m(x)<0$ for $x \in (\bar{q}_1,\tilde{q}_2).$ Now we claim that for $x \in (\bar{q}_1,\tilde{q}_2),$ it holds that $ x p(x)m(x)+2ps(s-p)\lambda(1-\lambda)[(s-1)x^{s-2}-(s-p)x^{p+s-3}-(p-1)x^{p-2}]\xi''(x)t(x)>0.$ Since $q  \in (\bar{q}_1,\tilde{q}_2),$ it yields that $\frac{d}{d\lambda} \Big( 1+z_2-\frac{\xi'(1)-\xi'(q)}{\xi''(q)(1-q)} \Big) \Big|_{1+z_2=\frac{\xi'(1)-\xi'(q)}{\xi''(q)(1-q)}} >0,$ which implies the desired conclusion.

Note that $ps\big[(p-1)q^{p-2}-(s-1)q^{s-2}+(s-p)q^{p+s-3}]>0,$ for $x \in (0,1)$ 
and $xm(x)+ps(s-p)\lambda(1-\lambda)\big ((p-1)x^{p-2}-(s-1)x^{s-2}+(s-p)x^{p+s-3}\big)(1-x)>0$ for $x \in (0,1).$ 
Also for $x \in (\bar{q}_1,\tilde{q}_2),$ it satisfies that $t^1_2(x)=x\xi'(x)\xi'''(x)-2\xi''(x)[x\xi''(x)-\xi'(x)]>0.$ 

Then we obtain that
\begin{align*}
& x p(x)m(x)+2ps(s-p)\lambda(1-\lambda)[(s-1)x^{s-2}-(s-p)x^{p+s-3}-(p-1)q^{p-2}]\xi''(x)t(x)\\
&=xm(x)[x\xi'(x)\xi'''(x)-2\xi''(x) \big( x\xi''(x)-\xi'(x) \big)][\xi'(1)-\xi'(x)]\\
&+ps(s-p)\lambda(1-\lambda)[(p-1)x^{p-2}-(s-1)x^{s-2}+(s-p)x^{p+s-3}][\xi'(x)xm(x)-2\xi''(x)t(x)]\\
&=[x\xi'(x)\xi'''(x)-2\xi''(x)(x\xi''(x)-\xi'(x)][\xi'(1)-\xi'(x)] \\
&\cdot \big [xm(x)+ps(s-p)\lambda(1-\lambda)\big ((p-1)x^{p-2}-(s-1)x^{s-2}+(s-p)x^{p+s-3}\big)(1-x) \big]\\
&>0,
\end{align*}
which finishes the proof of the claim.

By Lemma \ref{lemmalambda}, we obtain that $\lambda_{1 \rightarrow 2}<\lambda^*_1<\lambda_{2 \rightarrow 1F}.$ Then by the definition of $q^2_2(\lambda)$ and Lemma \ref{h4}$(v),$ if $\lambda \in (\lambda_{1 \rightarrow 2},\lambda_{2 \rightarrow 1F}),$ it holds that $q^2_2(\lambda)$ exists and $q^2_2(\lambda)<1$.
By Lemma \ref{lemmalambda} and the proposition of $\Psi(\lambda)$, we also obtain that $\Psi(\lambda_{2 \rightarrow 1F})<0 $ if and only if $\lambda_{2 \rightarrow 1F} > \lambda_{2 \rightarrow 1}.$

Recall that $\lim_{x \rightarrow 1}h^2_1(x,\lambda)=\Psi(\lambda)$ and $\Psi(\lambda_{2 \rightarrow 1})=0$, which implies that $q^2_1(\lambda_{2\rightarrow1})=1$.  Thus if $\lambda_{2 \rightarrow 1}< \lambda_{2 \rightarrow 1F}$, then it holds that $q^2_2(\lambda_{2 \rightarrow 1})<q^2_1(\lambda_{2 \rightarrow 1})$. By uniqueness of \eqref{eqnh2} and the fact that $q^2_2(\lambda_{1 \rightarrow 2})<q^2_1(\lambda_{1 \rightarrow 2})$, we obtain that \eqref{eqnh2} has no solution for $\lambda \in [\lambda_{1 \rightarrow 2}, \lambda_{2 \rightarrow 1}]$ and $x \in [0,1]$, which implies that $q^2_2(\lambda)<q^2_1(\lambda)$ for $\lambda \in (\lambda_{1 \rightarrow 2}, \lambda_{2 \rightarrow 1}).$  Since $q^1_2(\lambda)>q^1_1(\lambda)$ for $\lambda \in (\lambda_{1 \rightarrow 2}, \lambda_{2 \rightarrow 1}),$ we conclude that the model is 2RSB when $\lambda \in (\lambda_{1 \rightarrow 2}, \lambda_{2 \rightarrow 1}),$ $\Delta>0$, $\Psi(\lambda^*_1)>0$ and $\Psi(\lambda_{2 \rightarrow 1F})<0.$

Now we turn to the case when $\Delta>0$, $\Psi(\lambda^*_1)>0$ and $\lambda_{2 \rightarrow 1}> \lambda_{2 \rightarrow 1F}$. It then holds that $q^2_1(\lambda_{2 \rightarrow 1F})<q^2_2(\lambda_{2 \rightarrow 1F})=1.$ By the fact that $q^2_2(\lambda_{1 \rightarrow 2})<q^2_1(\lambda_{1 \rightarrow 2})$, the system of equations \eqref{eqnh2} has a unique solution $(\lambda_{2 \rightarrow 2F},x_{2 \rightarrow 2F})$ in $[\lambda_{1 \rightarrow 2}, \lambda_{2 \rightarrow 1F}] \times [0,1].$ Thus when $\lambda \in [\lambda_{1 \rightarrow 2}, \lambda_{2 \rightarrow 2F}],$ it holds that $q^1_2(\lambda)>q^1_1(\lambda)$ and $q^2_2(\lambda)<q^2_1(\lambda),$ which implies that the model is 2RSB and finishes the proof of Theorem \ref{coro12frsb}$(iii)$.

\end{proof}

\subsection{The 2FRSB Criterion} \label{sec2frsb}

We assume the Parisi measure $\nu_P$ has the form 
\begin{eqnarray*}
\nu(dx)=k_1 \cdot \mathbbm{1}_{[0,q_1)}(x)dx+\tilde{w}(x) \mathbbm{1}_{[q_1,q_2)}(x)+k_2 \cdot \mathbbm{1}_{[q_2,1)}(x)dx+\Delta \delta_{\{1\}}(dx),
\end{eqnarray*} 
where $\tilde{w}(s)$ is a strictly increasing continuous function on $[q_1,q_2).$

By Theorem \ref{criterion}, it holds that $g(x)=0$ for $x\in [q_1,q_2].$ Then we have $g'(x)=g''(x)=0$ for $x \in [q_1,q_2]$, which yields that
\begin{eqnarray}\label{w1}
&&\xi'(x)=\int^x_0 \frac{dr}{\nu((r,1])^2}  \nonumber \\
&& \nu((x,1])=\int^{q_2}_x \tilde{w}(s)ds +k_2(1-q_2)+\Delta=\xi''(x)^{-\frac{1}{2}} \text{ for } x \in [q_1,q_2].
\end{eqnarray}
In particular,
\begin{eqnarray}\label{w2}
\int^{q_2}_{q_1} \tilde{w}(s)ds +k_2(1-q_2)+\Delta=\xi''(q_1)^{-\frac{1}{2}} \text{ and } k_2(1-q_2)+\Delta=\xi''(q_2)^{-\frac{1}{2}},
\end{eqnarray}
and then
\begin{eqnarray}\label{w3}
\xi'(q_1)&=&\frac{1}{k_1[\int_{q_1}^{q_2} \tilde{w}(s)ds+k_2(1-q_2)+\Delta]}-\frac{1}{k_1[k_1q_1+\int^{q_2}_{q_1}\tilde{w}(s)ds+k_2(1-q_2)+\Delta]} \nonumber \\
&=&\frac{q_1}{\xi''(q_1)^{-\frac{1}{2}}[k_1q_1+\xi''(q_1)^{-\frac{1}{2}}]}.
\end{eqnarray}

By direct computation, we obtain that for $0 \leq x \leq q_1$,
 \begin{equation*} \nu((r,1])= \left\{
\begin{array}{lcl}
k_2(1-r)+\Delta, &\text{ for }& q_2 \leq r \leq 1,    \\
\int^{q_2}_r \tilde{w}(s)ds +k_2(1-q_2)+\Delta, &\text{ for }& q_1 < r < q_2, \\
k_1(q_1-r)+\int^{q_2}_{q_1}\tilde{w}(s)ds+k_2(1-q_2)+\Delta, &\text{ for }& 0\leq r \leq q_1. 
\end{array} \right. \end{equation*} 
Thus for $x \in [0,q_1],$
\begin{eqnarray*}
\int^x_0 \frac{dr}{\nu((r,1])^2}&=&\frac{1}{k_1[k_1(q_1-x)+\int_{q_1}^{q_2} \tilde{w}(s)ds+k_2(1-q_2)+\Delta]} \\
&&\qquad -\frac{1}{k_1[k_1q_1+\int^{q_2}_{q_1}\tilde{w}(s)ds+k_2(1-q_2)+\Delta]},
\end{eqnarray*}
for $x \in [q_1,q_2],$ by equation \eqref{w1},
\begin{eqnarray*}
\int^x_0 \frac{dr}{\nu((r,1])^2}&=&\frac{1}{k_1[\int_{q_1}^{q_2} \tilde{w}(s)ds+k_2(1-q_2)+\Delta]}\\
&& \qquad -\frac{1}{k_1[k_1q_1+\int^{q_2}_{q_1}\tilde{w}(s)ds+k_2(1-q_2)+\Delta]}+\xi'(x)-\xi'(q_1).
\end{eqnarray*}
Also, for
$x \in [q_2,1],$
\begin{eqnarray*}
\int^x_0 \frac{dr}{\nu((r,1])^2}&=&\frac{1}{k_1  [\int^{q_2}_{q_1} \tilde{w}(s)ds+k_2(1-q_2)+\Delta]}-\frac{1}{k_1 [k_1 q_1+\int^{q_2}_{q_1}\tilde{w}(s)ds+k_2(1-q_2)+\Delta]}\\
&&+\xi'(q_2)-\xi'(q_1)+\frac{1}{k_2[k_2(1-x)+\Delta]}-\frac{1}{k_2[k_2(1-q_2)+\Delta]}.
\end{eqnarray*}

Then it yields that for $x \in [0,q_1],$
\begin{eqnarray*}
g(x)&=&1-\xi(q_2)+\xi(q_1)-\xi(x)-[\xi'(q_2)-\xi'(q_1)](1-q_2)\\
&&+\frac{1-q_2}{k_2[k_2(1-q_2)+\Delta]}  -\frac{1}{k^2_1} \log \Big( 1+\frac{k_1(q_1-x)}{\int^{q_2}_{q_1}\tilde{w}(s)ds+k_2(1-q_2)+\Delta} \Big)\\
&&-\frac{q_1(1-q_2)}{[\int^{q_2}_{q_1}\tilde{w}(s)ds +k_2(1-q_2)+\Delta][k_1q_1+\int^{q_2}_{q_1}\tilde{w}(s)ds+k_2(1-q_2)+\Delta]} \\
&&+\frac{q_1-x}{k_1[k_1q_1+\int^{q_2}_{q_1} \tilde{w}(s)ds+k_2(1-q_2)+\Delta]} - \frac{1}{k^2_2} \log \Big (1+\frac{k_2(1-q_2)}{\Delta} \Big) ,
\end{eqnarray*}

for $q_2 \leq x \leq 1,$ 
\begin{eqnarray*}
g(x)&=&1-\xi(x)-[\xi'(q_2)-\xi'(q_1)](1-x)+\frac{1-x}{k_2[k_2(1-q_2)+\Delta]} -\frac{1}{k_2^2} \log \Big( 1+ \frac{k_2(1-x)}{\Delta}\Big)\\
&&-\frac{q_1(1-x)}{[\int^{q_2}_{q_1}\tilde{w}(s)ds+k_2(1-q_2)+\Delta][k_1q_1+\int^{q_1}_{q_2}\tilde{w}(s)ds + k_2(1-q_2)+\Delta]}.
\end{eqnarray*}

We then simplify the expressions of $g(x)$ as follows: for $x \in [0,q_1],$
\begin{eqnarray*}
g(x)&=&1-\xi(q_2)+\xi(q_1)-\xi(x)-[\xi'(q_2)-\xi'(q_1)](1-q_2)+\frac{1-q_2}{k_2[k_2(1-q_2)+\Delta]} \\
&&-\frac{q_1(1-q_2)}{\xi''(q_1)^{-\frac{1}{2}}[k_1q_1+\xi''(q_1)^{-\frac{1}{2}}]} - \frac{1}{k^2_2} \log \Big (1+\frac{k_2(1-q_2)}{\Delta} \Big) \\
&&+\frac{q_1-x}{k_1[k_1q_1+\xi''(q_1)^{-\frac{1}{2}}]} -\frac{1}{k^2_1} \log \Big (1+\frac{k_1(q_1-x)}{\xi''(q_1)^{-\frac{1}{2}}} \Big)
\end{eqnarray*}
and for $x \in [q_2,1],$
\begin{eqnarray*}
g(x)&=&1-\xi(x)-[\xi'(q_2)-\xi'(q_1)](1-x)+\frac{1-x}{k_2[k_2(1-q_2)+\Delta]}-\frac{1}{k_2^2} \log \Big (1+\frac{k_2(1-x)}{\Delta} \Big) \\
&&-\frac{q_1(1-x)}{\xi''(q_1)^{-\frac{1}{2}}[k_1q_1+\xi''(q_1)^{-\frac{1}{2}}]}.
\end{eqnarray*}
Also, $g(x)=g(q_2)$ for $x \in [q_1,q_2].$

Since $\xi'(1)=\int^1_0 \frac{dr}{\nu((r,1])^2}$, we obtain that
\begin{eqnarray*}
\xi'(1)&=&\frac{1}{k_1  [\int^{q_2}_{q_1} \tilde{w}(s)ds+k_2(1-q_2)+\Delta]}-\frac{1}{k_1 [k_1 q_1+\int^{q_2}_{q_1}\tilde{w}(s)ds+k_2(1-q_2)+\Delta]}\\
&&+\xi'(q_2)-\xi'(q_1)+\frac{1}{k_2\Delta}-\frac{1}{k_2[k_2(1-q_2)+\Delta]}\\
&=& \frac{q_1}{\xi''(q_1)^{-\frac{1}{2}}[k_1q_1+\xi''(q_1)^{-\frac{1}{2}}]}+\xi'(q_2)-\xi'(q_1)+\frac{1-q_2}{\Delta[k_2(1-q_2)+\Delta]}.
\end{eqnarray*}
Combining the equation \eqref{w2},\eqref{w3}, we obtain that
\begin{eqnarray}\label{delta}
\xi'(1)-\xi'(q_2)=\frac{1-q_2}{\Delta[k_2(1-q_2)+\Delta]} , \text{ then } \Delta=\frac{1-q_2}{\xi''(q_2)^{-\frac{1}{2}}[\xi'(1)-\xi'(q_2)]}
\end{eqnarray}
and
\begin{eqnarray}\label{k1}
k_1q_1=\frac{q_1\xi''(q_1)^{\frac{1}{2}}}{\xi'(q_1)}-\xi''(q_1)^{-\frac{1}{2}},\text{ then } k_1=\frac{q_1\xi''(q_1)-\xi'(q_1)}{q_1\xi'(q_1)\xi''(q_1)^{\frac{1}{2}}}.
\end{eqnarray}

By \eqref{w2}, we also have
\begin{eqnarray}\label{k2}
k_2(1-q_2)=\xi''(q_2)^{-\frac{1}{2}}-\frac{1-q_2}{\xi''(q_2)^{-\frac{1}{2}}[\xi'(1)-\xi'(q_2)]} = \frac{\xi'(1)-\xi'(q_2)-\xi''(q_2)(1-q_2)}{\xi''(q_2)^{\frac{1}{2}}[\xi'(1)-\xi'(q_2)]}.
\end{eqnarray}

Also by the condition $g(0)=g(q_1)=0$, it yields that
\begin{eqnarray}\label{g1}
\xi(q_1)+\frac{q_1}{k_1[k_1q_1+\xi''(q_1)^{-\frac{1}{2}}]} -\frac{1}{k^2_1} \log \Big (1+\frac{k_1q_1}{\xi''(q_1)^{-\frac{1}{2}}} \Big) =0,
\end{eqnarray}
and 
\begin{eqnarray}\label{g2}
&&1-\xi(q_2)-[\xi'(q_2)-\xi'(q_1)](1-q_2)+\frac{1-q_2}{k_2[k_2(1-q_2)+\Delta]} \nonumber\\
&&-\frac{q_1(1-q_2)}{\xi''(q_1)^{-\frac{1}{2}}[k_1q_1+\xi''(q_1)^{-\frac{1}{2}}]} - \frac{1}{k^2_2} \log \Big (1+\frac{k_2(1-q_2)}{\Delta} \Big)=0.
\end{eqnarray}

Thus by \eqref{delta},\eqref{k1} and \eqref{k2}, we can rewrite \eqref{g1} and \eqref{g2} as:
\begin{eqnarray*}
\xi(q_1)+\frac{q_1\xi'(q_1)^2}{q_1\xi''(q_1)-\xi'(q_1)}-\frac{(q_1)^2\xi'(q_1)^2\xi''(q_1)}{[q_1\xi''(q_1)-\xi'(q_1)]^2} \log \Big( \frac{q_1\xi''(q_1)}{\xi'(q_1)} \Big) =0,
\end{eqnarray*}
and 
\begin{align*}
&&-1+\xi(q_2)+\xi'(q_2)(1-q_2) -\frac{\xi''(q_2)[\xi'(1)-\xi'(q_2)](1-q_2)^2}{\xi'(1)-\xi'(q_2)-\xi''(q_2)(1-q_2)}\\
&&+ \frac{\xi''(q_2)[\xi'(1)-\xi'(q_2)]^2(1-q_2)^2}{[\xi'(1)-\xi'(q_2)-\xi''(q_2)(1-q_2)]^2} \log \Big (\frac{\xi'(1)-\xi'(q_2)}{\xi''(q_2)(1-q_2)} \Big)=0.
\end{align*}

Recall the definition of $h_1^2(x,\lambda)$ and $h^2_2(x,\lambda)$, we notice that $g(0)=0,g(q_1)=0$ is equivalent to $h_1^2(q_1,\lambda)=0,h_2^2(q_2,\lambda)=0.$

Now we turn to the proof of Theorem \ref{coro12frsb}$(iii)$ as follows:
\begin{proof}[Proof of Theorem \ref{coro12frsb}$(iii)$]
We first claim that
the model $\xi(x)$ is 2FRSB if there exists $0<q_1<q_2<1$ such that 
\begin{eqnarray}\label{cond2frsb1}
h^2_1(q_1,\lambda)=h^2_2(q_2,\lambda)=0
\end{eqnarray}
and
\begin{eqnarray}\label{cond2frsb2}
\frac{\xi'(1)-\xi'(q_2)}{\xi''(q_2)(1-q_2)}>\frac{\xi'(1)q_2}{\xi'(q_2)} \text{ and } \frac{\xi'(1)-\xi'(q_1)}{\xi''(q_1)(1-q_1)}>\frac{\xi'(1)q_1}{\xi'(q_1)}.
\end{eqnarray}

Now for $x\in[0,1],$ we consider $g_1(x)$ whose expression is given by $g(x)$ when $x \in [0,q_1]$. We also consider $g_2(x)$, whose expression is given by $g(x)$, for $x \in [q_2,1].$

It suffices for us to show that 
$g_1(x)>0$ when $x \in (0,q_1)$ if $ \frac{\xi'(1)-\xi'(q_1)}{\xi''(q_1)(1-q_1)}>\frac{\xi'(1)q_1}{\xi'(q_1)}$
and $g_2(x)>0$ when $x \in (q_2,1)$ if $\frac{\xi'(1)-\xi'(q_2)}{\xi''(q_2)(1-q_2)}>\frac{\xi'(1)q_2}{\xi'(q_2)} .$

We first prove the first part.
By computation, 
\begin{align*}
g_1'(x)&=-\xi'(x)-\frac{1}{k_1[k_1q_1+\xi''(q_1)^{-\frac{1}{2}}]}+\frac{1}{k_1[\xi''(q_1)^{-\frac{1}{2}}+k_1(q_1-x)]} \\
&=-\xi'(x)-\frac{\xi'(q_1)^2}{[q_1\xi''(q_1)-\xi'(q_1)]}+\frac{(q_1)^2\xi'(q_1)^2\xi''(q_1)}{[q_1\xi''(q_q)-\xi'(q_1)][q_1\xi'(q_1)+\big (q_1\xi''(q_1)-\xi'(q_1)\big)(q_1-x)]}\\
&=\frac{1}{[q_1\xi'(q_1)+\big(q_1\xi''(q_1)-\xi'(q_1)\big)(q_1-x)]} \cdot \\
& \qquad \Big \{ -(q_1)^2\xi''(q_1)\xi'(x)+[q_1\xi''(q_1)-\xi'(q_1)]x\xi'(x)+\xi'(q_1)^2x \Big \} \\
&:=\frac{1}{[q_1\xi'(q_1)+\big(q_1\xi''(q_1)-\xi'(q_1)\big)(q_1-x)]}  G_1(x).
\end{align*}
Then $g'_1(x)$ has the same sign as $G_1(x).$

Notice that 
\begin{align*}
G_1(x) =& \xi'(q_1)^2x-(q_1)^2\xi''(q_1) p \lambda x^{p-1} +[q_1\xi''(q_1)-\xi'(q_1)] \lambda px^p\\
&-(q_1)^2\xi''(q_1)s(1-\lambda)x^{s-1}+[q_1\xi''(q_1)-\xi'(q_1)](1-\lambda)sx^s.
\end{align*}
Then $G_1(x)$ has at most 4 changes of signs. By Descartes' rule of signs, $G_1(x)$ has at most 4 strictly positive zeros. Through direct computation, we obtain \begin{eqnarray*}
G'_1(q_1)=-(q_1)^2\xi''(q_1)^2+[q_1\xi''(q_1)-\xi'(q_1)]\xi'(q_1) +[q_1\xi''(q_1)-\xi'(q_1)]q_1\xi''(q_1)+\xi'(q_1)^2=0
\end{eqnarray*}
Since $G_1(q_1)=0$, it yields that $x=q_1$ is a zero with at most multiplicity 2.
Moreover,
\begin{eqnarray*}
G(1)=-(q_1)^2\xi''(q_1)\xi'(1)+[q_1\xi''(q_1)-\xi'(q_1)]\xi'(1)+\xi'(q_1)^2<0.
\end{eqnarray*}
Here we use the assumption $\frac{\xi'(1)-\xi'(q_1)}{\xi''(q_1)(1-q_1)}>\frac{\xi'(1)q_1}{\xi'(q_1)}.$
Also since $\lim_{x \rightarrow \infty}G_1(x)=+\infty,$ there exists $q_1'' \in (1,+\infty),$ such that $G_1(q_1'')=0.$ 

Recall that $g_1(0)=g_1(q_1)=0$, by mean value theorem, there exists $q_1' \in (0,q_1)$ such that $G_1(q_1')=0.$
Therefore the zeros of $G_1(x)$ are exactly  $q_1',q_1,q_1''$, which implies that $G_1(x)$ has only one zero in $(0,q_1).$

Notice that $G_1(x)$ is positive in a small neighborhood $(0,\delta_1)$ of 0, we conclude that $g_1(x)>0$ for $x \in (0,q_1).$

We then prove the second part of the claim. By computation, using \eqref{w3},\eqref{delta} and \eqref{k1},
\begin{align*}
g_2'(x)&=-\xi'(x)+\xi'(q_2)-\xi'(q_1)-\frac{1}{k_2[k_2(1-q_2)+\Delta]}\\
& \qquad +\frac{q_1}{\xi''(q_1)^{-\frac{1}{2}}[\xi''(q_1)^{-\frac{1}{2}}+k_1q_1]}+\frac{1}{k_2[\Delta+k_2(1-x)]} \\
&=-\xi'(x)+\xi'(q_2) +\frac{\xi''(q_2)^{\frac{1}{2}}(x-q_2)}{[\Delta+k_2(1-x)]} \\
&=-\xi'(x)+\xi'(q_2) +\frac{\xi''(q_2)[\xi'(1)-\xi'(q_2)](1-q_2)(x-q_2)}{\xi''(q_2)(1-q_2)^2+[\xi'(1)-\xi'(q_2)-\xi''(q_2)(1-q_2)](1-x)}\\
&=\frac{1}{\xi''(q_2)(1-q_2)^2+[\xi'(1)-\xi'(q_2)-\xi''(q_2)(1-q_2)](1-x)} \cdot \\
&\qquad \Big \{  \xi''(q_2)[\xi'(1)-\xi'(q_2)](1-q_2)(x-q_2) \\
&\qquad \qquad +\Big [[-\xi'(x)+\xi'(q_2)][\xi''(q_2)(1-q_2)^2+\Big(\xi'(1)-\xi'(q_2)\xi''(q_2)(1-q_2)\Big)(1-x) \Big] \Big \} \\
&:=\frac{1}{\xi''(q_2)(1-q_2)^2+[\xi'(1)-\xi'(q_2)-\xi''(q_2)(1-q_2)](1-x)} \cdot G_2(x)
\end{align*}

Note that $G_2(x)$ has 4 changes of signs since we assumed that $\frac{\xi'(1)-\xi'(q_2)}{\xi''(q_2)(1-q_2)}>\frac{\xi'(1)q_2}{\xi'(q_2)} .$ Then by Descartes' rule of signs, $G_2(x)$ has at most 4 strictly positive roots. Since $g_2(q_2)=g_2(1)=0$, there exists $q_2' \in (q_2,1)$ such that $g_2'(q_2')=G_2(q_2')=0$.

Moreover, a direct computation yields that
\begin{eqnarray*}
G_2(1)= \xi''(q_2)[\xi'(1)-\xi'(q_2)](1-q_2)^2+ \xi''(q_2)(1-q_2)^2[-\xi'(1)+\xi'(q_2)]   =0
\end{eqnarray*}
and
\begin{eqnarray*}
&&G_2(q_2)=0, \\
&&G_2'(q_2)= -\xi'(q_2)[\xi'(1)-\xi'(q_2)]+\xi'(1)\xi''(q_2)(1-q_2)\\
&&\qquad \qquad -[\xi'(1)-\xi'(q_2)]\xi''(q_2) +\xi''(q_2)^2q_2(1-q_2) -\xi''(q_2)^2q_2(1-q_2) \\
&&\qquad \qquad +[\xi'(1)-\xi'(q_2)]\xi'(q_2)-\xi'(q_1)\xi''(q_2)(1-q_2)+q\xi''(q_2)[\xi'(1)-\xi'(q_2)]=0,
\end{eqnarray*}
which implies that $x=1$ is a zero of $G_2(x)$ with multiplicity at least 1 and $x=q_2$ is a zero with multiplicity at least 2.

Thus $q_2,q_2',1$ are exactly all of the positive zeros of $G_2(x)$. It yields that there is only one critical point of $g_2(x)$ in $(q_2,1).$
 By assumption, the constant coefficient of $G_2(x)$ is strictly positive. Also since $q_2$ is a zero of $G_2(x)$ with multiplicity 2, there exists a neighborhood $(q_2-\delta_2,q_2)\cup (q_2,q_2+\delta_2)$ such that $G_2(x)>0$ for $x$ in it. It follows that $g_2(x)>0$ in $(q_2,q_2+\delta_2).$  Since $g_2(x)$ has only one critical point in $(q_2,1)$, we conclude that $g_2(x)>0$ for $x \in (q_2,1).$

Now it remains to show that $\xi''(x)^{-\frac{1}{2}}$ is concave for $x \in [q_1,q_2].$ Indeed,
\begin{align*}
\frac{d^2}{dx^2} \Big ( \xi''(x)^{-\frac{1}{2}} \Big)&=-\frac{1}{4}\xi''(x)^{-\frac{5}{2}}\Big (2\xi''(x)\xi''''(x)-3\xi'''(x)^2 \Big) \\
&=\frac{1}{4}\xi''(x)^{-\frac{5}{2}} x^{2p-6} \cdot \Big (s^3(s-1)^2(s-2)(1-\lambda)^2x^{2s-2p}+p^3(p-1)^2(p-2)\lambda^2\\
&-p(p-1)s(s-1)(2(p-2)(p-3)+2(s-2)(s-3)-6(p-2)(s-2)))\lambda(1-\lambda)x^{s-p}\Big)
\end{align*}
Since the right-hand side is quadratic with respect to $x^{s-p},$ it implies the range of $\xi''(x)^{-\frac{1}{2}}$ to be concave is an interval. Now we define the local minimum of $h^2_1(x,\lambda)$ and $h^2_2(x,\lambda)$ by $q'_1$ and $q'_2$. By Lemma \ref{h4}, it holds that $q'_1<q_1<q_2<q_2'.$ Then it suffices for us to show that $q'_1$ and $q'_2$ are in this interval.

We first show that $2\xi''(q_1')\xi''''(q'_1)-3\xi'''(q'_1)^2<0.$ Recall that the sign of $h^2_1(x,\lambda)$ is determined by $t^2_1(x).$
Since $q'_1$ is a local minimum of $h^2_1,$ it yields that
\begin{eqnarray*}
&&t^2_1(q'_1)=q'_1\xi'(q'_1)\xi'''(q'_1)-2\xi''(q'_1)[q'_1\xi''(q'_1)-\xi'(q'_1)] =0\\
&&q'_1\xi'(q'_1)\xi''''(q'_1)-3\xi'''(q'_1)[q'_1\xi''(q'_1)-\xi'(q'_1)]>0,
\end{eqnarray*}
which implies that
\begin{eqnarray*}
2\xi''(q_1')\xi''''(q'_1)-3\xi'''(q'_1)^2 &=&\frac{2\xi''(q_1')}{q_1'\xi'(q_1')} \cdot \Big( q'_1\xi'(q'_1)\xi''''(q'_1) -\frac{3q'_1\xi'(q'_1)\xi'''(q'_1)^2}{2\xi''(q'_1)}\Big) \\
&=&\frac{2\xi''(q_1')}{q_1'\xi'(q_1')} \cdot \Big( q'_1\xi'(q'_1)\xi''''(q'_1) -3[q'_1\xi''(q'_1)-\xi'(q'_1)]\xi'''(q'_1)\Big) \\
&>&0.
\end{eqnarray*}

It suffices to prove that $2\xi''(q_2')\xi''''(q'_2)-3\xi'''(q'_2)^2<0.$ Since $q'_2$ is a local minimum of $h^2_2,$ it satisfies that
\begin{eqnarray*}
&&m(q_2')=\xi'''(q_2')[\xi'(1)-\xi'(q_2')](1-q_2')-2\xi''(q_2')[\xi'(1)-\xi'(q_2')-\xi''(q_2')(1-q_2')] =0,\\
&&\xi''''(q_2')[\xi'(1)-\xi'(q_2')](1-q_2')-3\xi'''(q_2')[\xi'(1)-\xi'(q_2')-\xi''(q_2')(1-q_2')]>0.
\end{eqnarray*}
which implies that 
\begin{eqnarray*} 
&&2\xi''(q_2')\xi''''(q'_2)-3\xi'''(q'_2)^2\\
&=& \frac{2\xi''(q_2' )}{[\xi'(1)-\xi'(q_2' )](1-q_2' )} \Big ( \xi''''(q'_2)[\xi'(1)-\xi'(q_2' )](1-q_2' )-\frac{3\xi'''(q'_2)^2[\xi'(1)-\xi'(q_2' )](1-q_2' )}{2\xi''(q_2' )} \Big)\\
&=& \frac{2\xi''(q_2' )}{[\xi'(1)-\xi'(q_2' )](1-q_2' )} \Big ( \xi''''(q'_2)[\xi'(1)-\xi'(q_2' )](1-q_2' )-3\xi'''(q_2')[\xi'(1)-\xi'(q_2')-\xi''(q_2')(1-q_2')]\Big)\\
&>&0,
\end{eqnarray*}
as desired.

Now we prove that when $s>p>2$ satisfy that $\Delta>0,\Psi(\lambda^*_1)>0$ and $\Psi(\lambda_{2 \rightarrow 1})>0$, the model $\xi(x)$ is 2FRSB if $\lambda \in [\lambda_{2 \rightarrow 2F},\lambda_{2 \rightarrow 1F}].$ We claim that if we set $q_1=q^2_2(\lambda)$ and $q_2=q^2_1(\lambda),$ then the conditions \eqref{cond2frsb1} and \eqref{cond2frsb2} are satisfied, which implies that the model $\xi$ is 2FRSB. By the definition of $q^2_2(\lambda)$ and $q^2_1(\lambda)$, the condition \eqref{cond2frsb1} holds automatically. It then suffices for us to show that the condition \eqref{cond2frsb2} holds and that $q^2_1(\lambda) < q_2^2(\lambda) <1$ if  $\lambda \in (\lambda_{2 \rightarrow 2F},\lambda_{2 \rightarrow 1F}).$ 

Indeed, since the system of equations \eqref{eqnh2} has a unique solution $(\lambda_{2 \rightarrow 2F},x_{2 \rightarrow 2F})$ and $q^2_1(\lambda_{2 \rightarrow 1F})<q^2_2(\lambda_{2 \rightarrow 1F})=1$, it holds that $q^2_1(\lambda)<q^2_2(\lambda)<1$ for $\lambda \in (\lambda_{2 \rightarrow 2F}, \lambda_{2 \rightarrow 1F}).$ Since $\bar{q}_1<q^2_2<q^2_1<\bar{q}_2,$ the condition \eqref{cond2frsb2} is satisfied.
Thus the model $\xi(x)$ is 2FRSB for $\lambda \in (\lambda_{2 \rightarrow 2F}, \lambda_{2 \rightarrow 1F}).$

\end{proof}

\subsection{The 1FRSB Criterion} \label{sec1frsb}
We assume the Parisi measure $\nu_P$ has the form 
\begin{eqnarray*}
\nu(ds)=k_1 \cdot \mathbbm{1}_{[0,q_1)}(s)ds+\tilde{w}(s) \mathbbm{1}_{[q_1,1)}(s)+\Delta \delta_{\{1\}}(ds),
\end{eqnarray*} 
where $\tilde{w}(s)$ is a nondecreasing continuous function on $[q_1,q_2).$

Let
\begin{eqnarray*}
g(x)=\xi(q_1)-\xi(x)-[\xi'(1)-\xi'(q_1)](1-q_2)  +\frac{q_1-x}{k_1[k_1q_1+\xi''(q_1)^{-\frac{1}{2}}]} -\frac{1}{k^2_1} \log \Big (1+\frac{k_1(q_1-x)}{\xi''(q_1)^{-\frac{1}{2}}} \Big).
\end{eqnarray*}

By Theorem \ref{criterion}, 
it holds that
\begin{eqnarray*}
\xi'(x)=\int^x_0 \frac{dr}{\nu((r,1])^2} \text{ and } \nu((x,1])=\int^{1}_x \tilde{w}(s)ds +\Delta=\xi''(x)^{-\frac{1}{2}} \text{ for } x \in [q_1,1].
\end{eqnarray*}
In particular,
\begin{eqnarray*}
\int^{1}_{q_1} \tilde{w}(s)ds +\Delta=\xi''(q_1)^{-\frac{1}{2}} \text{ and } \Delta=\xi''(1)^{-\frac{1}{2}}.
\end{eqnarray*}
and then
\begin{eqnarray*}
\xi'(q_1)=\frac{q_1}{\xi''(q_1)^{-\frac{1}{2}}[k_1q_1+\xi''(q_1)^{-\frac{1}{2}}]}.
\end{eqnarray*}

Thus the equations $g(0)=g(q_1)$ can be rewritten as:
\begin{eqnarray*}
h_1^2(x,\lambda)=\xi(q_1)+\frac{q_1\xi'(q_1)^2}{q_1\xi''(q_1)-\xi'(q_1)}-\frac{(q_1)^2\xi'(q_1)^2\xi''(q_1)}{[q_1\xi''(q_1)-\xi'(q_1)]^2} \log \Big( \frac{q_1\xi''(q_1)}{\xi'(q_1)} \Big) =0,
\end{eqnarray*}

We are now ready to prove Theorem \ref{coro12frsb}$(iv)$.
\begin{proof}[Proof of Theorem \ref{coro12frsb}$(iv)$]
The proof of Theorem \ref{coro12frsb}$(iv)$ is similar to the proof of Theorem \ref{coro12frsb}$(iii)$. We omit the details of the calculation and present only the key steps.

We first claim that the model $\xi(x)$ is 1FRSB if 
(a) there exists $q_1 \in (0,1)$ such that 
\begin{eqnarray*}
h^2_1(q_1)=0 \text{ and } \frac{\xi'(1)-\xi'(q_1)}{\xi''(q_1)(1-q_1)}>\frac{\xi'(1)q_1}{\xi'(q_1)}.
\end{eqnarray*}
(b) for $q \in (q_1,1), h^2_2(x,\lambda) \neq 0.$

Since $ \frac{\xi'(1)-\xi'(q_1)}{\xi''(q_1)(1-q_1)}>\frac{\xi'(1)q_1}{\xi'(q_1)},$ it holds that $g(x)>0$ for $x \in (0,q_1).$ The proof is the same as the one in the proof of Theorem \ref{coro12frsb}$(iv)$. By the same reasoning as the proof of Theorem \ref{coro12frsb}$(iv)$, it holds that $\xi''(x)^{-\frac{1}{2}}$ is concave in $(q_1,1)$.
For $\lambda \in (\lambda_{2 \rightarrow 1F},\lambda_{2 \rightarrow 1})$, it holds that $q^1_1,q^1_2,q^1_2,q^2_2$ exists and $q^2_1<q^2_2=1$. Then by the definition of $q^2_1(\lambda)$ and $q^2_2(\lambda),$ we set $q_1=q^2_1(\lambda)$, which implies that the model $\xi(x)$ is 1FRSB.
\end{proof}

\begin{proof}[\bf Proof of Corollary \ref{cor:finiteTemp}]
The proof of Corollary \ref{cor:finiteTemp} is similar to the proof of Corollary 1 in \cite{auffinger2017sk}.
\end{proof}

\bibliographystyle{abbrv}
\bibliography{biblio3}

\end{document}